\documentclass[reqno,12p]{amsart}
\usepackage[T1]{fontenc}

\usepackage{geometry}
\usepackage{mathtools,amssymb,amsthm,mathrsfs,color,lineno,paralist,graphicx,float}

\definecolor{orange}{rgb}{1,0.5,0}

\definecolor{green}{rgb}{0.5,0.7,0.3}

\usepackage[colorlinks,
linkcolor=red,
anchorcolor=green,
citecolor=blue,
]{hyperref}
\usepackage{amsmath}
\usepackage{pdfpages}
\usepackage{amstext}
\usepackage{stmaryrd}
\usepackage{bbm}
\usepackage{comment}
\usepackage{faktor,url}
\usepackage{xfrac}
\usepackage{mathabx}
\usepackage{enumitem}

\usepackage[normalem]{ulem}

\numberwithin{equation}{section}
\numberwithin{figure}{section}
\theoremstyle{plain}
\newtheorem{thm}
{\protect\theoremname}[section]
  \theoremstyle{definition}
  \newtheorem{defn}[thm]{\protect\definitionname}
  \theoremstyle{rem}
  \newtheorem{rem}{\protect\remarkname}
  \theoremstyle{plain}
  \newtheorem{lem}[thm]{\protect\lemmaname}
  \theoremstyle{plain}
  \newtheorem{prop}[thm]{\protect\propositionname}
  \theoremstyle{plain}

\newtheorem{thmx}{Theorem}
\newtheorem{thmxprime}{Theorem}

  \providecommand{\definitionname}{Definition}
  \providecommand{\lemmaname}{Lemma}
  \providecommand{\propositionname}{Proposition}
  \providecommand{\remarkname}{Remark}
\providecommand{\theoremname}{Theorem}
\providecommand{\corollaryname}{Corollary}

\newcommand{\bs}{\boldsymbol}

\title{Non-local ergodic measures in symplectic dynamics}

\author[J. Paradela]{J. Paradela}
\address[JP]{Instituto de Ciencias Matemáticas, Consejo Superior de Investigaciones Científicas, 28049
Madrid, Spain}
\email{jaimeparadela@gmail.com}

\begin{document}

\begin{abstract}
    We consider the direct product of two symplectomorphisms, one of which exhibits a basic set and the other one acts on a compact surface and has a non-degenerate elliptic equilibrium. Under a local domination condition we show that a typical \textit{real-analytic}  perturbation of this system displays a large  normally hyperbolic lamination on which the natural invariant measure is \textit{ergodic}. 

    The proof, which introduces a novel geometric source of semi-local ergodicity, combines tools ranging from KAM theory to geometric measure theory and is robust enough to be of interest in more general configurations.
\end{abstract}

\maketitle
\tableofcontents

\section{Introduction}

A longstanding problem in classical mechanics is to understand the topological and measure theoretic
properties of typical symplectic dynamics.  After Oxtoby and Ulam established genericity of ergodicity
among conservative homeomorphisms \cite{MR5803}, the groundbreaking work of Kolmogorov,  Arnold and Moser,
showed that for sufficiently regular topologies, there exist open sets of non-ergodic symplectomorphisms
(see \cite{MR209436,MR2269239} and the references therein). Later,  Herman constructed open
sets of symplectomorphisms which are not transitive (see the survey article by Yoccoz \cite{MR1206072}). These results show that the typical symplectic dynamics varies drastically from low to high regularity categories. 

\medskip

On the low regularity side of the spectrum, genericity of transitivity was established in the $C^1$ symplectic setting
by Arnaud, Bonatti and Crovisier \cite{MR2173426} (following earlier work of Bonatti and Crovisier \cite{MR2090361} in
the conservative scenario).  On the measure theoretic counterpart, the typicality of ergodicity has been intensively studied under extra assumptions such as the existence of a global partially hyperbolic splitting (i.e. the tangent bundle decomposes into three invariant bundles: an expanding one, a contracting one and a center one, the latter being ``dominated'' by the other two). Crucially this additional structure excludes the obstructions to ergodicity described above and several important results have been obtained,  including also some which apply in regular topologies (see \cite{MR3682778,MR4422613,MR4422214} and the references therein).
\medskip

For very regular categories,  in the absence of a global partially hyperbolic splitting, much less is known.  The phenomenon which better highlights the subtleties hindering ergodicity in the symplectic smooth scenario is the robust persistence of elliptic behavior. On one hand, under mild non-degeneracy conditions and suitable regularity assumptions, it is well known that elliptic equilibria are
robustly accumulated by invariant tori whose union forms a set of positive measure and on which the dynamics is conjugated to a Diophantine translation (we refer the interested reader to \cite{MR3146593}). We call this union the KAM set, as the existence of these tori is deduced by means of the machinery developed by Kolmogorov, Arnold and Moser. Note that this implies the existence of uncountably many ergodic components.
\medskip

Still, one can hope to find  measures with large support and wild statistical properties in the complement of the KAM set. This is the main goal of the present work.  Before stating our main result, we first describe the class of dynamics that we will  study (real-analytic Hamiltonian diffeomorphisms), identify a class of structures that under mild assumptions can be found in the complement of the KAM set (normally hyperbolic laiminations) and discuss heuristically why these objects are natural candidates to be the support of ergodic measures.

\subsection*{The Liouville measure}
Given a smooth  manifold $X$ we will denote by $\mathcal M(X)$ the set of Borel  measures on $X$. In the case where $(X,\varpi)$ is a symplectic manifold, we define the \textit{Liouville measure} $\mathrm{vol}\in \mathcal M(X)$ as the unique Borel measure satisfying 
\[
\mathrm{vol}(A)=\int_A \varpi^{d},\qquad\qquad d=\frac 12\mathrm{dim}(X)
\]
for every Borel measurable set $A\subset X$.
\subsection{Hamiltonian diffeomorphisms}
 Given a real-analytic   symplectic manifold $(X,\varpi)$  we we denote by 
\begin{equation}\label{eq:spaceofHams}
\mathcal H^\omega(X)=\{H:X\times [0,1]\to \mathbb R\colon H\text{ is real-analytic}\}
\end{equation}
the space of real-analytic  time-dependent Hamiltonians equipped with the real-analytic topology (see Appendix \ref{sec:realanalytictop} for the precise definition of this topology).
We write $H_t(\cdot)=H(\cdot,t)$ and recall that to any $H\in \mathcal H^\omega(X)$ we can associate a time-dependent Hamiltonian vector field $X_{H_t}: X\to TX$ defined (implicitly) by
\begin{equation}\label{eq:correspondenceformsvectors}
\varpi(\cdot,X_{H_t})=\mathrm dH_t(\cdot).
\end{equation}
We denote by 
\begin{align*}
\Phi_H:X&\to X
\end{align*}
the time-one map associated to the flow of this vector field and recall that it is a symplectic transformation, i.e. $(\Phi_H)^*\ \varpi=\varpi$. In particular, $\Phi_{H}$ preserves the \textit{Liouville measure} $\mathrm {vol}\in \mathcal M(X)$ defined by $\varpi$. Finally, we introduce the group of real-analytic Hamiltonian diffeomorphisms
\[
\mathrm{Ham}^\omega(X)=\{\Phi_H:X\to X\colon H\in \mathcal H^\omega(X)\}.
\]
The interested reader is referred to  \cite{MR3674984} for an introduction to Hamiltonian and symplectic dynamics.\medskip

\subsection{Direct products of Hamiltonian diffeomorphisms}

Consider two real-analytic symplectic manifolds $(N,\varpi_N)$ and $(M,\varpi_M)$ together with Hamiltonian diffeomorphisms $F_N\in\mathrm{Ham}^\omega(N,\varpi_N)$ and $F_M\in \mathrm{Ham}^\omega(M,\varpi_M)$. Then, we define their \textit{direct product} 
    \begin{equation}\label{eq:directproductmap}
\begin{aligned}
  F_0:M\times N&\to M\times N\\
  (x_M,x_N)&\mapsto (F_M(x_M),F_N(x_N)).
\end{aligned}
\end{equation}
In Theorem \ref{thm:mainham} below we will show that ergodic measures with large support can be found in typical real-analytic unfoldings of direct products  satisfying certain assumptions.
    Note that the direct product \eqref{eq:directproductmap} is again a real-analytic Hamiltonian diffeomorphism (i.e. $F_0\in \mathrm{Ham}^\omega(M\times N)$). Indeed,  there exist $H_N\in \mathcal H^\omega(N)$ and $H_M\in \mathcal H^\omega(M)$ such that $F_N=\Phi_{H_N}$ and $F_M=\Phi_{H_M}$. Then, an equivalent definition of the direct product \eqref{eq:directproductmap} is $F_0=\Phi_{H_0}$ where $H_0=H_M+H_N\in \mathcal H^\omega(M\times N)$.
\medskip

\subsubsection*{Non-degenerate elliptic equilibria}
Our first assumption (see \ref{it:Hyp1} below) requires a mild amount of elliptic behavior for the dynamics on $N$. Recall that a  fixed point $q\in N$ of a  symplectic map $F:N\to N$ is \textit{elliptic} if the spectrum of the linearization of $F$ at $q$ is simple and it is contained  in the unit circle in $\mathbb C$. In addition, we say that an elliptic fixed point $q$ is  \textit{non-degenerate} if there exists a neighbourhood  $U\subset N$ of $q$ and a local coordinate system $\psi:z\in \mathbb B^d\subset\mathbb C^{d}\to U\subset N$ such that $\psi^*\varpi_N=i \mathrm d\bar z\wedge \mathrm dz$ and conjugates $F$ to 
    \[
    \psi^{-1}\circ F\circ\psi:z\mapsto  z(\exp (i \Omega(|z|^2))+O(|z|^3))
    \]
    where $|z|^2$ is shorthand notation for $(|z_1|^2,\dots,|z_d|^2)\in \mathbb R^d$ and  $\Omega(|z|^2)=\mathrm{diag}(\Omega_1(|z|^2),\dots,\Omega_d(|z|^2))$ is a diagonal matrix satisfying
    \[
    \mathrm{det} D\Omega(0)\neq 0.
    \]
Namely, up to higher order terms, the map $F$ is conjugated to a product of rotations with frequencies $\{\Omega_i(|z|^2)\}_i$ and the frequency map $|z|^2\mapsto \Omega(|z|^2)$ is a local diffeomorphism. 

\medskip

Given a direct product of Hamiltonian diffeomorphisms $F_0=F_M\times F_N$ as in \eqref{eq:directproductmap} we say that it satisfies the property:
\begin{enumerate}[label={[H\arabic*]}]
    \item \label{it:Hyp1} if the map $F_N:N\to N$ exhibits a non-degenerate elliptic equilibrium.
    \end{enumerate}

\subsubsection*{Normally hyperbolic submanifolds of direct products}

Our second assumption (\ref{it:Hyp2} below) asks for a mild amount of hyperbolicity for the dynamics on $M$. To state it precisely we need the following definition.

\begin{defn}\label{defn:bunchingpinching}
   Let $F_0=F_M\times F_N$ be the direct product of two diffeomorphisms $F_M:M\to M$ and $F_N:N\to N$. Let $a\in M$ be a hyperbolic fixed point for $F_M$ and let $T_a M=E^s_a\oplus E^u_a$ be the corresponding splitting of the tangent space of $M$ at $a$. For $\kappa\geq 0$ we say that the invariant submanifold 
   \[
   \Lambda=\{a\}\times N\subset M\times N
   \]
   is $\kappa$-normally hyperbolic  if there exist metrics $\lVert\cdot\rVert_M$ on $M$ and $\lVert\cdot\rVert_N$ on  $N$ such that
       \begin{equation}\label{eq:bunchingdefn}
         \lVert  DF_N^{-1}\rVert_N^{\kappa+1} \lVert DF_M|_{E^s_a}\rVert_M < 1\qquad\qquad\text{and}\qquad\qquad \lVert DF_N\rVert^{\kappa+1}_N \lVert DF_M^{-1}|_{E^u_a}\rVert_M< 1.
       \end{equation}
\end{defn}

In words, the stable  forward contraction (unstable backwards contraction) dominates strongly the center backwards (forward) behavior. It is well known that  the exponent $\kappa$ in this property can be related to the regularity of certain invariant foliations associated to (the continuation) of the submanifold $\Lambda$ and which are important to  understand the homoclinic picture associated to  $\Lambda$ (see  \cite{MR501173} or Section 3 of \cite{paradela2026robustlytransitivebehaviorsymplectic}). 

\begin{defn}
    In the setting of Definition \ref{defn:bunchingpinching} we say that $\Lambda$ is \textit{$\kappa$-pinched} if \begin{equation}\label{eq:pinchingdefn}
        \bar \alpha^{\kappa+1} \lambda<1
       \end{equation}
       where 
       \[
       \bar\alpha=\frac{\max\{
\lVert DF_M|_{E_a^s}\rVert_M,\,
\lVert DF_M^{-1}|_{E_a^u}\rVert _M
\}}{\min\{
\lVert DF_M^{-1}|_{E_a^s}\rVert^{-1}_M,\,
\lVert DF_M|_{E_a^u}\rVert^{-1} _M
\}}\qquad\text{and}\qquad \lambda=\max\{
\lVert DF_M|_{E_a^s}\rVert_M,\,
\lVert DF_M^{-1}|_{E_a^u}\rVert _M
\}.
       \]
\end{defn}

Pinching compares the maximal rate of normal expansion/contraction with the minimal rate of normal expansion/contraction. Observe that in the case $\mathrm{dim}(M)=2$ the assumption is trivial.
\medskip

Given a direct product of Hamiltonian diffeomorphisms $F_0=F_M\times F_N$ as in \eqref{eq:directproductmap} we say that it satisfies property:
 \begin{enumerate}[label={[H\arabic*]},start=2]
    \item \label{it:Hyp2} 
    if the map $F_M:M\to M$ displays a hyperbolic fixed point $a\in M$ whose stable and unstable manifolds intersect transversally and, moreover, $\Lambda=\{a\}\times N$ is $\kappa$-normally hyperbolic and $\kappa$-pinched for  $F_0$. 
    \end{enumerate}

 \subsubsection*{Emergence of ergodic laminations} 
   Note that under assumption \ref{it:Hyp2}, the classical Smale horseshoe Theorem (see \cite{MR1326374}) shows that for any alphabet $A\subset \mathbb N$, the map $F_M:M\to M$ displays an open set $Q\subset M$ and a subset   $\Sigma\subset Q$ homeomorphic to $A^\mathbb Z$ such that $\Sigma$ is a basic set for the return map $R$ induced on $Q$. Hence, if we let $\mathcal Q_0=Q\times N$, the return map 
   \begin{equation}\label{eq:returnmap unperturbed heuristic}
   \mathcal R_0:\breve{\mathcal Q}_0\subset\mathcal Q_0\to \mathcal Q_0
   \end{equation}
   induced by $F_0$ on $\mathcal Q_0$,  exhibits an invariant set 
   \begin{equation}\label{eq:lamin unperturbed heuristic}
   \mathcal X_0=\Sigma\times N\subset M\times N.
   \end{equation}
Moreover, thanks to the domination assumption in \ref{it:Hyp2}, $\mathcal X_0$ is robust under small perturbations  \cite{MR501173}. Although for $\mathcal R_0|_{\mathcal X_0}$ the fiber dynamics (along $N$) is decoupled from the chaotic base dynamics (along $\Sigma$) driven by the shift, one expects that after an arbitrarily small perturbation, the chaoticity of the base ``contaminates'' the fiber dynamics, generating rich dynamics in the continuation of $\mathcal X_0$. This intuition has been confirmed in a number of works (see for instance \cite{Moeckeldrift,MR2104598,NassiriPujalsTransitivity,GelfreichTuraev,MR3917797,guardia2025partiallyhyperbolicdynamics3body,paradela2026robustlytransitivebehaviorsymplectic} and the references therein). 
\medskip

From the topological viewpoint, the pioneering work of Nassiri and Pujals \cite{NassiriPujalsTransitivity} shows that any $F_0$ satisfying \ref{it:Hyp1} can be $C^\infty$ approximated by a symplectomorphism for which the continuation $\mathcal X$ of the lamination $\mathcal X_0$ is \textit{robustly transitive}. Later, the main result in \cite{paradela2026robustlytransitivebehaviorsymplectic}  shows that the same conclusion indeed holds for a locally dense set of real-analytic perturbations of a Hamiltonian diffeomorphism $F_0$ satisfying \ref{it:Hyp1}-\ref{it:Hyp2}.  The main tool behind these constructions is a \textit{symplectic blender}, a robust source of semi-local transitivity introduced by Nassiri and Pujals in \cite{NassiriPujalsTransitivity} which extends to the symplectic setting the classical notion of blender introduced by Bonatti and Díaz \cite{BonattiDiazblenders,BeyondUH}. To upgrade the construction of symplectic blenders in \cite{NassiriPujalsTransitivity} (which relied heavily on compactly supported perturbations) to the real-analytic setting, \cite{paradela2026robustlytransitivebehaviorsymplectic} extends previous ideas developed by the author and Guàrdia in \cite{guardia2025partiallyhyperbolicdynamics3body}\footnote{A closely related mechanism for the creation of symplectic blenders in the real-analytic category was previously developed by Li and Turaev \cite{LiTuraevPreprint}. A comparison between the two mechanisms is given in \cite{paradela2026robustlytransitivebehaviorsymplectic}.}. 
\medskip
    \begin{figure}
    \centering
\includegraphics[scale=0.55]{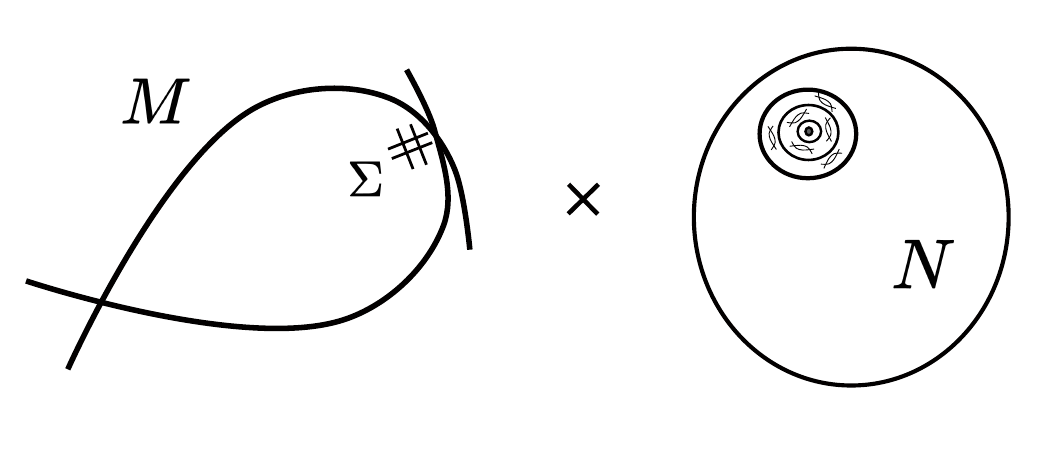}
    \caption{The basic set for $F_M$ and the elliptic equilibrium for $F_N$.}
    \label{fig:fig1}
\end{figure}

Besides the topological viewpoint, a natural next question is the study of statistical properties for  measures supported on the lamination $\mathcal X$. In particular, in Theorem \ref{thm:mainham} below,  we will be interested in the ergodicity of the following natural class of invariant measures.

\begin{defn}[Natural measure on a symplectic lamination]\label{defn:natural measure}
Let $\mathrm{vol}\in \mathcal M(M\times N)$ stand for the Liouville measure on $M\times N$ defined by the product symplectic form $\varpi=\varpi_M\times \varpi_N$. Let $A\subset \mathbb N$ (possibly equal to $\mathbb N$) and denote by $\Omega_A=A^\mathbb Z$. We say that a set $\mathcal X$ is a \textit{symplectic lamination} if there exists  a homeomorphism $\Psi:\Omega_A\times N\to \mathcal X\subset M\times N$ such that for each $\omega\in \Omega_A$:
\begin{itemize}
    \item  the leaf $N_\omega=\Psi(\{\omega\}\times N)$ is a $C^1$ manifold,
    \item  the two-form $\varpi_\omega=(\iota_\omega)^*   \varpi$ defined in terms of the inclusion map $\iota_\omega:N_\omega\to M\times N$ is non-degenerate, hence a symplectic form.
\end{itemize} 
Then, we say that a measure $ \mathfrak m$ on $\mathcal X$ is \textit{natural} if its pushforward under $\Psi$ is of the form 
    \[
    \mathrm d\mathfrak m (\omega,z)=\mathrm d \mu_{\mathrm{Ber}}(\omega)\  \mathrm d \mathrm{vol}_\omega(z)
    \]
    where $\mu_{\mathrm{Ber}}$ is any Bernoulli probability measure  on $\Omega_A$ with strictly positive weights and $
    \mathrm{vol}_\omega\in \mathcal M(N_\omega)$ is Liouville measure associated to the symplectic form $\varpi_\omega$.
\end{defn}

Before proceeding it is worth discussing the invariance of the measure $\mathfrak m$. Denote by $\sigma:\Omega_A\to \Omega_A$ the full shift acting on the space of bi-infinite sequences.  Then,  given a homeomorphism $\Psi:\Omega_A\times N\to \mathcal X\subset M\times N$ as above, an open set $\mathcal Q\supset \mathcal X$ and a symplectic map $\mathcal R:\mathcal Q\to M\times N$ satisfying $\mathcal R(N_\omega)=N_{\sigma\omega}$,  we must have $
(\mathcal R|_{N_\omega})^*\ \varpi_{\sigma\omega}=\varpi_\omega$. Consequently, $
    (\mathcal R|_{N_\omega})_*\  \mathrm{vol}_\omega=\mathrm{vol}_{\sigma\omega}$. 
    Since any Bernoulli measure is invariant under the shift we conclude that 
    \[
    (\mathcal R|_{\mathcal X})_*\  \mathfrak m=\mathfrak m.
    \]
    Namely, the measure $\mathfrak m$ is invariant for $\mathcal R|_{\mathcal X}:\mathcal X\to \mathcal X$. 
\medskip

    For our purposes in this work  it will be enough to think of $\mathcal R$ as a return map (induced by symplectic a map $F$) to an open set $\mathcal Q$ (so symplecticity is automatic) in the neighbourhood of a homoclinic manifold to a normally hyperbolic manifold. Indeed, in the direct product case one can take $\mathcal R_0$ as in \eqref{eq:returnmap unperturbed heuristic} and $\mathcal X_0$ as in \eqref{eq:lamin unperturbed heuristic}.
\medskip

\subsection{Real-analytic Hamiltonian  unfoldings of direct products}
Given $r\in\mathbb N$ define the open subset 
\begin{equation}\label{eq:Ballperturbations}
\mathcal B^\omega_r(X)=\{H\in \mathcal H^\omega(X)\colon |H|_{C^{r}(X)}< 1\}
\end{equation}
of real-analytic Hamiltonians. 
In Theorem \ref{thm:mainham} below, given a direct product of Hamiltonian diffeomorphisms $F_0=F_M\times F_N$ as in \eqref{eq:directproductmap} associated to a Hamiltonian $H_0=H_M+H_N$ we consider one-parametric unfoldings of the form
\begin{equation}\label{eq:deformationmap}
F_{\varepsilon,f}=\Phi_{H_{\varepsilon, f}} \qquad \text{where}\qquad H_{\varepsilon, f}= H_0+\varepsilon f\qquad\text{and}\qquad f\in \mathcal B^\omega_r(M\times N).
\end{equation}
    Observe that for any $f\in \mathcal B^\omega(X)$ and any $0\leq |\varepsilon|\leq 1$, the  time-one map
    \begin{equation}\label{eq:unfolding}\Phi_{\varepsilon f}:M\times N\to M\times N
    \end{equation}
      is uniformly $C^{r-1}$-close to the identity. In particular, $F_{\varepsilon,f}$ is $C^{r}$-close to $F_0=\Phi_{H_0}$.
\medskip

\subsection{Main results}

We recall that a measure is \textit{ergodic} if the only invariant subsets have measure zero or one.

\subsection*{Ergodic measures with large support for Hamiltonian diffeomorphisms}

The following is our first main result.

\begin{thmx}\label{thm:mainham}
Let  $(M,\varpi_M)$ and $(N,\varpi_N)$ be real-analytic symplectic manifolds with $N$ being a compact surface. There exists  $\kappa_0>0$ and $r_0\in\mathbb N$ such that the following hods. Let $F_0\in \mathrm{Ham}^\omega(M\times N)$ be a direct product of Hamiltonian diffeomorphisms  which satisfies the assumptions \ref{it:Hyp1}-\ref{it:Hyp2} with $\kappa\geq\kappa_0$. 

Then, for any $f$ belonging to an open and dense subset of $\mathcal B^\omega_{r_0}(M\times N)$ there exists 
$\varepsilon_0(F_0,f)>0$ such that for any $-\varepsilon_0\leq \varepsilon\leq \varepsilon_0$, the Hamiltonian diffeomorphism $
   F_{\varepsilon,f}:M\times N\to M\times N$ in \eqref{eq:deformationmap} displays an open set $\mathcal Q_\varepsilon\subset M\times N$ and a subset $\mathcal X_\varepsilon\subset \mathcal Q_\varepsilon$  homeomorphic to $\mathbb N^\mathbb Z\times N$ such that: 
    \begin{enumerate}
        \item $\mathcal X_\varepsilon$ is a locally maximal normally hyperbolic symplectic lamination for the return map $\mathcal R_\varepsilon$ induced on $\mathcal Q_\varepsilon$. Moreover, its projection onto $N$ is the whole manifold.

        \item Any natural measure $\mathfrak m_\varepsilon$ on $\mathcal X_\varepsilon$ is invariant for $\mathcal R_\varepsilon$ and,
        \item 
        Any natural measure  $\mathfrak m_\varepsilon$ on $\mathcal X_\varepsilon$ is  ergodic of Lebesgue almost every $\varepsilon \in (-\varepsilon_0,\varepsilon_0)$.
    \end{enumerate}
\end{thmx}

We present the proof of Theorem \ref{thm:mainham} in Section \ref{sec:outlineThmA} by reducing it to a number of more technical statements which will be proved in the remaining sections of the paper.  Before continuing, let us discuss briefly the assumptions and main novelties of our result. We later review the relationship with previous works (Section \ref{sec:previousworks}) and discuss some extensions of our main result (Section \ref{sec:futurework}).

\begin{rem}
    Although Theorem \ref{thm:mainham} is presented only for $\mathrm{dim}N=2$, a higher dimensional version can be obtained with the same techniques used in this paper. We have opted for focusing here on the two-dimensional version since the geometric mechanism behind the proof is more transparent in this setting. The higher dimensional version will appear in a separate note.
\end{rem}
\medskip

\begin{figure}
    \centering
\includegraphics[scale=0.45]{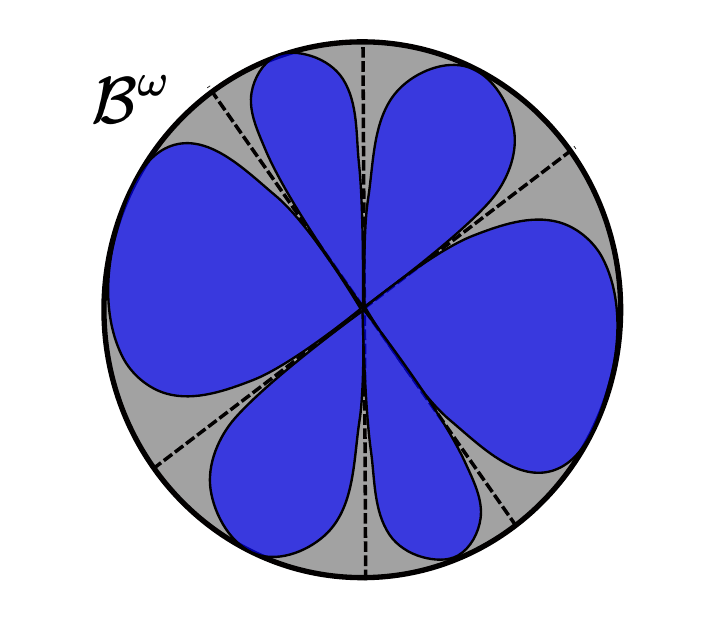}
    \caption{We depict in blue the shape of the set of perturbations within $\mathcal H^\omega(M\times N)$ for which Theorem \ref{thm:mainham} holds. The excluded directions correspond to the dashed lines. Along any ray within the blue region ergodicity of the natural measure holds almost everywhere.}
    \label{fig:fig2}
\end{figure}

\subsubsection*{Domination   assumption:}  We have not tried to optimize the lower value of $\kappa_0$ for which Theorem \ref{thm:mainham} holds. In any case it is worth remarking that in some cases of interest, the domination assumption in \ref{it:Hyp2} is automatically satisfied. This is the case if, for instance, the unperturbed dynamics on $N$ is given by a completely integrable map all whose singularities are elliptic and $\mathrm{dim}(M)=2$. A natural example is to consider $N=\mathbb C\mathbb P$   equipped with the Fubini-Study form $\varpi_{\mathrm FS}$ and let $F_N$ be the time-one map associated to a Hamiltonian given by a non-linear function of the moment map coordinates (see \cite{MR1853077}).
\medskip

\subsubsection*{Strong form of topological instability}
In the case where the unperturbed dynamics on $N$ is integrable, classical results show that, after arbitrarily small perturbations, there exist orbits which leak through the remnants of the invariant toric foliation (of the integrable system) and connect regions separated by a distance independent of the size of the perturbation. Some recent contributions are  \cite{MR3646879,GelfreichTuraev,MR4298716,MR4033892,MR4729212,MR4913967}. 

Under stronger assumptions (direct product structure of the unperturbed map), our result implies that a typical real-analytic perturbation, gives rise to the existence of ergodic measures with large support  in the sense that  its projection to $N$ does not depend on the size of the perturbation (indeed it coincides with   $N$).
\medskip

\noindent\textit{Explicit assumptions:} Not only \ref{it:Hyp1}-\ref{it:Hyp2} but also the conditions defining the open and dense subset of $\mathcal B^\omega(M\times N)$ for which Theorem \ref{thm:mainham} holds are completely explicit. Notably, they only involve first order (in $\varepsilon$) conditions on a explicit collection of \textit{integrals depending linearly on $f$} and can be verified in given one-parametric families making use of perturbative computations or via computer assisted technology. 

We refer the reader to Section \ref{sec:outline} and in particular to Proposition \ref{prop:Reduction Skew-prod to IFS} for a precise description of the set of assumptions.
\medskip

\noindent\textit{Ergodicity of natural measure:} One of the corollaries of the main results in \cite{NassiriPujalsTransitivity} is that a direct product  $F_0$ satisfying \ref{it:Hyp1} can be approximated in the $C^\infty$ topology by symplectomorphisms exhibiting ergodic measures with large support. Namely, the support of these measures contains a lamination (a set  homeomorphic to $\{1,\dots,k\}^\mathbb Z\times N$ for some $k\geq 2$) whose projection onto $N$ is large. However, the existence of these measures is deduced by an abstract argument (obtained in \cite{MR2811152}) which gives little additional description of the measure.

In Theorem \ref{thm:mainham}, we are able to show ergodicity of the natural measure on the lamination $\mathcal X_\varepsilon$, giving a partial affirmative answer to Problem 6.3 in \cite{NassiriPujalsTransitivity}.
\medskip

\noindent\textit{A novel semi-local source of ergodicity:} The main ingredient in the proof of Theorem \ref{thm:mainham} is  a geometric construction which implies \textit{local ergodicity} for a suitable class of iterated function systems on surfaces. We refer the reader to Section \ref{sec:IFSintro} and in particular to Theorem \ref{thm:main}, for a detailed discussion of this mechanism.

The main drawback of this local construction, tailored for one-parametric families, is that it only produces almost everywhere results in the parameter (hence the $\varepsilon$-almost everywhere conclusion in Theorem \ref{thm:mainham}). In particular, with our current approach we cannot guarantee that the ergodic measures in Theorem \ref{thm:mainham} are stably ergodic.

\medskip

\subsection*{Ergodic measures with large support for symplectomorphisms}
   Theorem \ref{thm:mainham} has an analogue for symplectomorphisms which we now present. To do so we will use the  well known fact  that any  neighbourhood of the identity in  $\mathrm{Symp}^\omega(X,\varpi)$ (the group of real-analytic symplectomorphisms, see Appendix \ref{sec:realanalytictop}) can be identified with a neighbourhood of zero in the vector space of real-analytic closed one-forms (see \cite{MR3674984}). More precisely, we denote by $\mathfrak H^\omega$ the space of real-analytic closed one-forms $\theta:X\to T^*X$ equipped with the real-analytic topology  and let
\begin{equation}\label{eq:...}
\mathfrak B^\omega_r(X)=\{\theta\in \mathfrak h^\omega(X)\colon |\theta|_{C^{r}}\leq 1\}.
\end{equation}
To  any $\theta\in \mathfrak B^\omega_r(X)$ we can associate a symplectic vector field $X_\theta: X\to TX$ defined (implicitly) as
\begin{equation}\label{eq:correspondenceformsvectors}
\varpi(\cdot,X_\theta)=\theta(\cdot).
\end{equation}
In this way, we obtain a one-parametric family of real-analytic symplectomorphisms $\{\Phi^\varepsilon_\theta\}_\varepsilon:X\to X$ defined implicitly by the initial value problem
\begin{equation}\label{eq:deformationmap one form}
\begin{cases}
   \frac{\mathrm{d}\Phi^\varepsilon_\theta}{\mathrm d\varepsilon}=X_\theta\circ\Phi^\varepsilon_\theta\\
  \Phi^0_\theta=\mathrm{id}.
\end{cases}
\end{equation}

\begin{rem}
    Observe that for any $\theta\in \mathfrak B^\omega_r(X)$ and any $0\leq |\varepsilon|\leq 1$, the real-analytic symplectomorphism $\Phi^\varepsilon_\theta:X\to X$ defined by \eqref{eq:correspondenceformsvectors}-\eqref{eq:deformationmap} is uniformly $C^{r-1}$-close to identity.
\end{rem}

\begin{thmxprime}\label{thm:mainsymp}
Let $(N,\varpi_N)$ and $(M,\varpi_M)$ be  real-analytic symplectic manifolds and let $N$ be a compact surface. There exists $\kappa>0$ and $r_0\in\mathbb N$  such that the following holds. Let $F_0\in \mathrm{Symp}^\omega(M\times N)$ be a direct product of symplectomorphisms  which satisfies assumptions \ref{it:Hyp1}-\ref{it:Hyp2} with $\kappa\geq \kappa_0$.
 
Then, for any $\theta$ belonging to an open and dense subset of $\mathfrak B^\omega_{r_0}(M\times N)$  there exists $\varepsilon_0(F_0,\theta)>0$ such that for  any $-\varepsilon_0<\varepsilon\leq \varepsilon_0$, the map
\begin{equation}\label{eq:unfoldingforsymplectic}
   F_{\varepsilon,\theta}=F_0\circ\Phi_{\theta}^\varepsilon\in \mathrm{Symp}^\omega(M\times N)
   \end{equation}
   with $\Phi_{\theta}^\varepsilon$ as in \eqref{eq:deformationmap one form} displays an open set $\mathcal Q_\varepsilon\subset M\times N$ and a subset $\mathcal X_\varepsilon\subset \mathcal Q_\varepsilon$  homeomorphic to $\mathbb N^\mathbb Z\times N$ such that: 
       \begin{enumerate}
        \item $\mathcal X_\varepsilon$ is a locally maximal  normally hyperbolic symplectic lamination for the return map $\mathcal R_\varepsilon$ induced on $\mathcal Q_\varepsilon$. Moreover, its projection onto $N$ is the whole manifold.

        \item Any natural measure $\mathfrak m_\varepsilon$ on $\mathcal X_\varepsilon$ is invariant for $\mathcal R_\varepsilon$ and,
        \item 
        Any natural measure  $\mathfrak m_\varepsilon$ on $\mathcal X_\varepsilon$ is  ergodic of Lebesgue almost every $\varepsilon \in (-\varepsilon_0,\varepsilon_0)$.
    \end{enumerate}
\end{thmxprime}

\medskip

\subsection{Comparison with previous literature}\label{sec:previousworks}

We now compare our result and techniques with other works on related topics.

\subsubsection*{Works on stable ergodicity of partially hyperbolic symplectomorphisms}

We will see below the dynamics on $\mathcal X_\varepsilon$ is driven by a partially hyperbolic skew-product.  Hence, one may wonder if the techniques developed in \cite{MR3682778,MR4422613,MR4422214} (and related works on stable ergodicity)  can be adapted to deduce Theorem \ref{thm:mainham}.  In all these works,  ergodicity is achieved by showing that the system has the so-called \textit{accessibility} property.  Namely,  for (almost every) pair of points $z,z'$ in the manifold there exists a $su$-path, i.e. a piecewise defined path which alternates segments of strong stable and strong unstable manifolds, which connects them.  
\medskip

To the best of our understanding,  the route followed to establish accessibility in these works relies on the possibility of generating a ``continuum'' of $su$-loops (su-paths to a fixed reference center leaf).  However,  due to the fractal nature of the lamination $\mathcal X_\varepsilon$, this does not seem possible in the setting of Theorem \ref{thm:mainham}.  Therefore, it is not clear how to adapt the techniques in these works to obtain a proof of Theorem \ref{thm:mainham}.

\subsubsection*{Works on statistical properties of random compositions}

DeWitt and Dolgopyat have recently  obtained a number of remarkable results concerning the statistical properties of a broad class of random dynamical systems. Let us describe these results in connection with those in Theorem \ref{thm:mainham}.  Consider a smooth surface $N$ equipped with a  smooth probability measure $\mathrm{vol}$,  a Bernoulli probability measure $\mu_{\mathrm{Ber}}$ on the sequence space $\Omega=\{1,\dots, k\}^\mathbb N$, and denote by 
\[
\frak m=\mu_{\mathrm{Ber}}\otimes \mathrm{vol}.
\]
Given a tuple of $C^r$ ($r\geq 1$) diffeomorphisms $\{F_i\}_{i=1}^k:N\to N$ preserving $\mathrm{vol}$ define the (locally constant) skew-product 
\[
\mathcal F:(\omega,z)\in \Omega\times N\mapsto (\sigma\omega, F_{\omega_0}(z))
\]
with $\sigma:\Omega\to \Omega $ being the right shift. Then,  a particular case of the main result in \cite{dewitt2024expandingaveragediffeomorphismssurfaces} shows that the system $(\mathcal F,\Omega\times N, \mathfrak m)$   \textit{mixes exponentially} provided it is  \textit{expanding on average}. A spectacular generalization of this result to higher dimensions can be found in \cite{dewitt2025conservativecoexpandingaveragediffeomorphisms}. Even more recently  \cite{bénard2026rigidityequidistributionrandomwalks,brown2026absolutecontinuityinvariantmeasures} have established strong rigidity and equidistribution results for  stationary measures of random diffeomorphisms satisfying certain variations/generalizations of the expanding on average property. Notably, the results in \cite{brown2026absolutecontinuityinvariantmeasures} apply to random dynamics on non-compact manifolds.
\medskip

Provided one could establish that in the setting of Theorem \ref{thm:mainham} the corresponding skew-product driving the dynamics on the lamination is expanding on average, a small adaptation of their result to bypass the fact that the skew-product is not locally-constant (although the coupling is weak), would yield exponential mixing for the natural measure. The abundance of expanding on average tuples has been studied in a number of works \cite{MR4477288,MR4656974}.  We believe that the expanding on average property should hold typically in the setting of Theorem \ref{thm:mainham} but we have not been able to verify this assertion. See Section \ref{sec:nihimintro} (in particular, Proposition \ref{prop:mainnhim}) for a detailed description of the skew-product driving the dynamics on the lamination $\mathcal X_\varepsilon$ in Theorem \ref{thm:mainham}.
\medskip

An alternative route to make use of the results in \cite{dewitt2024expandingaveragediffeomorphismssurfaces} to deduce Theorem \ref{thm:mainham} would be the following. We will see in the proof of Theorem \ref{thm:mainham} that it can be reduced to study ergodicty of a suitable iterated function system on a fixed reference leaf. This iterated function system is generated by a countable collection of maps and contains, in particular, a countable family of $su$-loops to the fixed reference leaf. It seems plausible that one could find words formed by compositions of these $su$-loops and a probability measure supported on these words for which the corresponding random dynamics is expanding on average and hence exponential mixing.  Translated to the original problem, i.e. the statistical properties of the measure $\mathfrak m_\varepsilon$, this would imply ergodicity. However,  one cannot deduce exponential mixing for $\mathfrak m_\varepsilon$ with this approach.

Although this approach gives much more flexibility (compared to trying to establish the expanding on average property for $\mathfrak m_\varepsilon$) it also presents some technical  challenges and it would be worth investigating it, as it would provide an alternative proof of Theorem \ref{thm:mainham} at the expense of a more involved description of the open and dense set of perturbations for which ergodicty of $\mathfrak m_\varepsilon$ holds.

\subsection{Extensions of Theorem \ref{thm:mainham}}\label{sec:futurework}

We now discuss briefly some possible extensions of Theorem \ref{thm:mainham}. A more technical discussion can be found in Section \ref{sec: extensions section}.

\subsubsection*{Statistical properties of the natural measure}
Another interesting extension is to study statistical properties of the natural measure besides ergodicity, such as establishing a Central Limit Theorem, exponential mixing, etc...  With our current approach, although it would require some non-trivial work,  it seems plausible to upgrade the conclusion of Theorem \ref{thm:mainham} from ergodicity of the natural measure to $K$-property which, in particular, implies mixing of all orders (although with no rates), see \cite{MR209436}.

\subsubsection*{Beyond the direct product setting} It also seems possible to  obtain an analogue of Theorem \ref{thm:mainham} just assuming the existence of a normally hyperbolic lamination with compact leaves and suitable bunching/pinching conditions. A priori, this extension does not seem limited to the two-dimensional center setting and could be carried directly in the higher dimensional case. Similarly, the natural strategy  for proving the $K$-property does not rely on working close to a direct product.

\medskip

\subsubsection*{Ergodic measures with support of maximal Hausdorff dimension} Suppose now that $\mathrm{dim}(M)=2$ and  that for $\varepsilon=0$ the lamination $\mathcal X_0$ in Theorem \ref{thm:mainham} has Hausdorff dimension $0<d\leq 4=\mathrm{dim}(M)+\mathrm{dim}(N)$. Under some additional technical work it seems possible to show that there exists $\theta=\theta(\varepsilon,f)\to 1$ as $\varepsilon \to 0$ such that $\mathrm{HD}(\mathcal X_\varepsilon)\geq \theta d$. In particular, by working close to a homoclinic tangency, one could obtain ergodic measures whose support has Hausdorff dimension arbitrarily close to maximal (see also Section 1.5 in \cite{NassiriPujalsTransitivity}). A more demanding problem would be to construct ergodic measures whose support has maximal Hausdorff dimension.

\subsection{Acknowledgements}
The author is supported by the ICMAT–Severo Ochoa grant CEX2023-001347-S funded by MICIU/AEI/10.13039/501100011033.

\subsection{Declaration of AI use}
All of the mathematical ideas of this work are due to the author except for the trick in Lemma \ref{lem:weigh integration estimate growth affine}, which was suggested by ChatGPT.  During the final preparation of the article 
the author used ChatGPT for proofreading. The author takes full responsibility for the mathematical arguments
and the final content of the article.

\section{Definitions, strategy and proof of Theorem \ref{thm:mainham}}\label{sec:outlineThmA}

In this section we describe the strategy of the proof of Theorem \ref{thm:mainham} and introduce the most relevant definitions. 
\medskip
 
 First, in Section \ref{sec:nihimintro} we recall the definition of \textit{normally hyperbolic manifold} and prepare ourselves to study \textit{homoclinic phenomena} associated to normally hyperbolic manifolds such as the existence of \textit{normally hyperbolic laminations}.  To that end we review the construction of the so-called \textit{scattering map} introduced by Delshams, de la Llave and Seara \cite{DelaLLavegaps,DelaLLaveScattmap}. As we will see below, the scattering map, defined as a composition of stable and unstable holonomies, encodes the dynamics along the set of homoclinic orbits to a normally hyperbolic manifold. Having in mind the setting in Theorem \ref{thm:mainham}, in Proposition \ref{prop:mainnhim} we give a number of  sufficient conditions  for the existence of a  normally hyperbolic lamination and relate the dynamics on the lamination to the scattering map dynamics. 
 \medskip

Sections \ref{sec:IFSintro}-\ref{sec:globalergodicity} are devoted to the study of certain \textit{iterated function systems} (IFS's) on compact manifolds. The reason for this (apparent) detour is that the dynamics in the lamination in Proposition \ref{prop:mainnhim} is given by a \textit{skew-product} over the shift. Hence one might expect to be able to translate some dynamical constructions from the IFS setting to the skew-product setting. 

In particular, having in mind the structure of the skew-product dynamics in Proposition \ref{prop:mainnhim}, our goal in Sections \ref{sec:IFSintro}-\ref{sec:globalergodicity} is to present sufficient conditions for a suitable class of IFS's to be \textit{ergodic}. This is done in two steps. 

The first step (which occupies Section \ref{sec:IFSintro})  focuses on the creation of  a \textit{semi-local source} of \textit{ergodicity} which we have called \textit{Hopf boxes}. Sufficient conditions for their construction are given in  Proposition \ref{prop:TT and HB}.

The second step (in Section \ref{sec:globalergodicity}) involves \textit{global transport} in the IFS, namely, wether the IFS has the property of connecting any two points up to an small (but finite) error. Exploiting the fact  that the IFS's of interest for us include a couple of close to identity maps,  a simple argument borrowed from the \textit{control theory} literature shows that this is the case. Combined, local ergodicity and global transport force global ergodicity of the corresponding IFS (see Theorem \ref{thm:ergodicity IFS final}).

\medskip

In Section \ref{sec:skewprodintro},  building on the results   obtained in Sections \ref{sec:IFSintro}-\ref{sec:globalergodicity}, we give sufficient conditions which guarantee that the skew-product in  Proposition \ref{prop:mainnhim} is ergodic with respect to any natural measure. This is the content of Proposition \ref{prop:Reduction Skew-prod to IFS}.
\medskip

Finally, in Section \ref{sec:genericityintro} we recall a result obtained in \cite{paradela2026robustlytransitivebehaviorsymplectic} which shows that for an open and dense set of perturbations in $\mathcal B^\omega(M\times N)$ (recall its definition in \eqref{eq:Ballperturbations}), the corresponding one-parametric family of Hamiltonian diffeomorphisms 
\[
F_{\varepsilon,f}\in \mathrm{Ham}^\omega(M\times N)
\]
as in \eqref{eq:deformationmap}: 
\begin{enumerate}
    \item satisfies the asumptions in Proposition \ref{prop:mainnhim} and,
    \item  the corresponding skew-product dynamics  can be understood in terms of an IFS which verifies the assumptions introduced in Section \ref{sec:IFSintro} (namely the ones needed to apply Propositions \ref{prop:TT and HB} and \ref{prop:IFS property R}). In particular, Proposition \ref{prop:Reduction Skew-prod to IFS} applies.
\end{enumerate}

\begin{rem}\label{rem:remarksofttechniques}
    The presentation in  Sections \ref{sec:nihimintro}-\ref{sec:skewprodintro} does not rely on the particular scenario of Theorem \ref{thm:mainham} and is presented in a more abstract setting. In particular, throghout these sections we will consider $C^r$  diffeomorphisms  $F:X\to X$ of a compact Riemannian manifold with $r\geq 2$ and do not assume any product structure for $X$ nor that $F$ is a Hamiltonian diffeomorphism.  Only in Section \ref{sec:genericityintro} we will come back to the setting of Theorem \ref{thm:mainham}.
\end{rem}

\medskip

\subsection{Normally hyperbolic cylinders and homoclinic channels}\label{sec:nihimintro}
Consider a $C^r$ diffeomorphism $F:X\to X$ of a  Riemannain manifold $X$ with $r\geq 2$.
\begin{defn}\label{defn:nhim}
    Let $\Lambda\subset X$ be an invariant submanifold for $F$, i.e. $F(\Lambda)=\Lambda$. We say that $\Lambda$ is \textit{normally hyperbolic} if there exists $C>0$ and constants
    \[
    0<\lambda<\alpha^{-1}<1
    \]
    such that for any $x\in \Lambda$ there exists a splitting $T_x X=E^s_x\oplus T_x\Lambda \oplus E^u_x$ verifying
    \begin{align*}
        v\in E^s_x &\Longleftrightarrow \lVert DF^n v\rVert\leq C\lambda^n \lVert v\rVert,\qquad n\geq 0\\
         v\in T_x\Lambda &\Longleftrightarrow \lVert DF^n v\rVert\leq C\alpha^{|n|} \lVert v\rVert,\qquad n\in\mathbb Z\\
          v\in E^u_x &\Longleftrightarrow \lVert DF^n v\rVert\leq C\lambda^{-n} \lVert v\rVert,\qquad n\leq 0.
    \end{align*}
\end{defn}

\begin{rem}
    It is not difficult to check that the particular definition we gave for the direct product case (see Definition \ref{defn:bunchingpinching}) satisfies the requirements in Definition \ref{defn:nhim} (see for instance \cite{MR501173}).
\end{rem}

Given a normally hyperbolic manifold we define its \textit{stable and unstable sets}
\begin{equation}\label{eq:stablemanifoldsdefn}
\begin{aligned}
W^s(\Lambda)=&\{x_*\in X\colon \exists x\in \Lambda \text{ such that }\mathrm{dist}(F^n(x),F^n(x_*))\to 0\text{ as } n\to +\infty\}\\
W^u(\Lambda)=&\{x_*\in X\colon \exists x\in \Lambda \text{ such that }\mathrm{dist}(F^n(x),F^n(x_*))\to 0\text{ as } n\to -\infty\}
\end{aligned}
\end{equation}
and, given $x\in \Lambda$ introduce
\begin{equation}\label{eq:stablemanifoldsfibersdefn}
\begin{aligned}
W^s_x=&\{x_*\in X\colon \mathrm{dist}(F^n(x),F^n(x_*))\to 0\text{ as } n\to +\infty\}\\
W^u_x=&\{x_*\in X\colon \mathrm{dist}(F^n(x),F^n(x_*))\to 0\text{ as } n\to -\infty\}.
\end{aligned}
\end{equation}
The following facts are standard (see \cite{MR501173,DelaLLaveScattmap})
\begin{enumerate}
   
 \item The manifold $\Lambda$ is $C^l$ with 
    \begin{equation}\label{eq:spectralgap}
        1\leq l<\mathrm{min}\left(r,\frac{|\log\lambda|}{\log \alpha}\right),
    \end{equation}

    \item The maps  $x\mapsto W^{u,s}_x$ are $C^{l-1-j}$ when $W^{u,s}_x$ are given the $C^j$ topology and at any $x\in \Lambda$ we have $T_x W^{u,s}_x=E_x^{u,s}$,

    \item $W^{u,s}_x$ are $C^{r}$  manifolds for any $x\in \Lambda$.    
\end{enumerate}
In particular, it follows that $\{W^{u,s}_x\}_{x\in\Lambda}$ give a $C^{l-1}$ foliation of $W^{u,s}(\Lambda)$.
\medskip

\subsubsection*{Homoclinic channels and the scattering map}

We now introduce some tools to study homoclinic phenomena associated to normally hyperbolic manifolds. Our presentation follows closely that in \cite{DelaLLaveScattmap} and we refer the interested reader there for a much more complete exposition. 

We start by introducing the \textit{wave maps}
\begin{align*}
\Omega^{u,s}:W^{u,s}(\Lambda)&\mapsto \Lambda\\
x&\mapsto x^{u,s}
\end{align*}
where $x^{u,s}\in \Lambda$ is the unique point such that $x\in W^{u,s}_{x^{u,s}}$. Notice that these maps are $C^{l-1}$.
\begin{defn}[Homoclinic channel] Let $\Lambda$ be a normally hyperbolic manifold of a $C^r$ diffeomorphism $F:X\to X$ and let  $\Omega^{u,s}:x\in W^{u,s}(\Lambda)\to x^{u,s}\in\Lambda$ be the corresponding wave maps. We say that a submanifold $\Gamma\subset W^s(\Lambda)\cap W^u(\Lambda)$ is a \textit{homoclinic channel} to $\Lambda$ if at any $x\in \Gamma$
\[
T_x\Gamma=T_x W^s(\Lambda)\cap T_x W^u(\Lambda)\qquad\qquad T_x \Gamma \oplus T_x W^{u,s}_{x^{u,s}}=T_x W^{u,s}(\Lambda)
\]
\end{defn}

By analogy with their lower dimensional counterpart, i.e. hyperbolic periodic orbits displaying transverse homoclinic orbits to them, one usually expects that the existence of homoclinic channels to a normally hyperbolic submanifold leads to some sort of ``chaotic behavior''. 

One of the more relevant tools for the study of the dynamics around homoclinic channels is the so-called scattering map, introduced by Delshams, de la Llave and Seara in \cite{DelaLLavegaps,DelaLLaveScattmap}.

\begin{defn}[Scattering map]
    Let $\Lambda$ be a normally hyperbolic manifold of a $C^r$ diffeomorphism  $F:X\to X$, let $\Omega^{u,s}:x\in W^{u,s}(\Lambda)\to x^{u,s}\in \Lambda$ be the corresponding wave maps and let $\Gamma\in W^s(\Lambda)\cap W^u(\Lambda)$ be a homoclinic channel. Then, we define the \textit{scattering map}
    \begin{align*}
    S_\Gamma: D_\Gamma\subset \Lambda\to \Lambda\qquad\qquad \text{ where }\qquad D_\Gamma=\Omega^u|_\Gamma(\Gamma)
    \end{align*}
    by 
    \[
S_\Gamma=\Omega^s|_\Gamma\circ (\Omega^u|_\Gamma)^{-1}: x^u\mapsto x^s
    \]
\end{defn}
Observe that, by construction, the scattering map is a $C^{l-1}$ diffeomorphism onto their image.

\medskip

\subsubsection*{Existence of normally hyperbolic laminations}
 Throughout this section we will consider a one-parametric family of $C^r$ $(r\geq  2$) diffeomorphisms  $\{F_\varepsilon\}_{\varepsilon}$ of a Riemaniann manifold $X$ (see Remark \ref{rem:remarksofttechniques}). We shall suppose that there exists $\varepsilon_0(F_0)>0$ such that for any $0\leq \varepsilon\leq\varepsilon_0(F_0)$:
\begin{enumerate}[label={[A\arabic*]}]
    \item \label{it:H1} \textit{$\kappa$-normal hyperbolicity:} $F_\varepsilon:X\to X$ has a compact normally hyperbolic invariant manifold $\Lambda_\varepsilon$ which associated rates $0<\lambda_\varepsilon<\alpha_\varepsilon^{-1}<1$ satisfying
    \begin{equation}\label{eq:spectralgap}
        \alpha_\varepsilon^{\kappa+1}\lambda_\varepsilon<1
    \end{equation}
    for some $\kappa>0$. Moreover, we assume that  there exists a reference manifold $N$ and a parametrization $\psi_\varepsilon:N\to \Lambda_\varepsilon$. We denote by 
    \begin{equation}\label{eq:innerdynamicsdefn}
        T_\varepsilon=\psi_\varepsilon^{-1}\circ F_\varepsilon|_{\Lambda_\varepsilon}\circ\psi_\varepsilon:N\to N
    \end{equation}
the map describing the restriction of the dynamics to $\Lambda_\varepsilon$.

\item\label{it:H1+1/2} \textit{$\kappa$-pinching:} Let 
\[
\underline\lambda_\varepsilon=\min\{\lVert D F_\varepsilon^{-1}|_{E^s}\rVert^{-1},\lVert D F_\varepsilon|_{E^u}\rVert^{-1}\}\qquad\text{and}\qquad \bar\alpha_\varepsilon=\frac{\lambda_\varepsilon}{\underline\lambda_\varepsilon}.
\]
    Then, 
    \begin{equation}\label{eq:pinching}
    \bar\alpha_\varepsilon^{\kappa+1}\lambda_\varepsilon <1
    \end{equation}
    for some $\kappa>0$.
    \item \label{it:H2} \textit{Full homoclinic intersections:} There exist $m$ ($m\geq 1$)  homoclinic manifolds 
    \[\{\Gamma^i_\varepsilon\}_{i=1}^m\subset W^s(\Lambda_\varepsilon)\pitchfork W^u(\Lambda_\varepsilon)
    \]
    such that for every $i=1,\dots,m$ the wave maps    $\Omega^{u,s}|_{\Gamma^i_\varepsilon}:\Gamma^i_\varepsilon\to \Lambda_\varepsilon$  are global diffeomorphisms and 
    \[
\bigg( \bigcup_{n\in\mathbb Z} F_\varepsilon^n(\Gamma^i_\varepsilon)\bigg)\cap \bigg( \bigcup_{n\in\mathbb Z} F_\varepsilon^n(\Gamma^j_\varepsilon)\bigg)=\emptyset.
    \]
    for any pair $i\neq j$.We denote by $S_{\Gamma^i_\varepsilon}:\Lambda\to \Lambda$ the corresponding scattering maps and let 
    \begin{equation}\label{eq:scattmapdynamicsdefn}S_{i,\varepsilon}=\psi_\varepsilon^{-1}\circ S_{\Gamma^i_{\varepsilon}}\circ\psi_\varepsilon:N\to N.
    \end{equation}
\end{enumerate}
\medskip

Before continuing, let us make some comments on the assumptions. First of all, one should notice that since both conditions are open, it follows from regular dependence on parameters (see \cite{MR501173,DelaLLaveScattmap}), that it is only necessary to verify Assumptions \ref{it:H1}-\ref{it:H2} for $F_0$ and then, they automatically hold for all $F_\varepsilon$ with $\varepsilon>0$ sufficiently small (a more detailed discussion can be found in Section 3 of \cite{paradela2026robustlytransitivebehaviorsymplectic}). We prefer to state the assumptions already for the family to avoid innecessary repetitions. 

Concerning the nature of the assumptions,   \ref{it:H1} guarantees that the stable and unstable foliations associated to $\Lambda_0$ are at least $C^{\kappa}$ if $r\geq\kappa$ while \ref{it:H2} will be  of use to show that the homoclinic picture associated to $\Lambda_\varepsilon$ is sufficiently rich.   The assumption \ref{it:H1+1/2} is of more technical nature: it measures the ratio between maximal and minimal normal contraction or expansion. 
\medskip

The following result describes the set in which, under suitable assumptions to be made precise later, supports a natural ergodic measure.

\begin{prop}[Existence of a NHIL]\label{prop:mainnhim}
    There exists $r_0\in\mathbb N$ and $\kappa_0>0$ such that the following hold. Let $F_\varepsilon:X\to X$ be a parametric family of $C^r$ maps $(r\geq r_0)$   which satisfy assumptions \ref{it:H1}-\ref{it:H2} with $\kappa\geq \kappa_0$. Let $m\geq 2$ be as in \ref{it:H2}, let  $A=\{1,\dots,m\}$ and let $\Omega=(\mathbb N\times A)^\mathbb Z$.

    Then, there exists \[
    \bar\lambda\in (0,1),\qquad C(F_0)>0\qquad\text{and}\qquad \varepsilon_0(F_0)>0
    \]
    such that  for any $0\leq \varepsilon\leq \varepsilon_0$ there exists a set $\mathcal X_\varepsilon\subset X$, an open neighbourhood $\mathcal Q_\varepsilon\supset \mathcal X_\varepsilon$   and a homeomorphism $\Psi_\varepsilon:\Omega\times N\to \mathcal X_\varepsilon$ such that the following holds:
    \begin{enumerate}
        \item \textit{(Normal hyperbolicity):} The set $\mathcal X_\varepsilon$ is $\kappa$-normally hyperbolic and locally  for the return map $\mathcal R_\varepsilon$ induced by $F_\varepsilon$ on $\mathcal Q_\varepsilon$,
        \item \textit{(Skew-product dynamics):} The homeomorphism $\Psi_\varepsilon$ conjugates the dynamics of $\mathcal R_\varepsilon$ to the skew-product
        \begin{equation}\label{eq:skewprodnhim}
        \begin{aligned}
        \mathcal F_\varepsilon=\Psi^{-1}_\varepsilon\circ\mathcal R_\varepsilon\circ\Psi_\varepsilon:\Omega\times N&\to \Omega\times N\\
        (\bs\omega,z)&\mapsto (\sigma(\bs\omega), F_{\bs\omega,\varepsilon}(z))
        \end{aligned}
        \end{equation}
        where $\sigma$ is the full right shift and, for every  $\bs\omega=(\omega,\iota)\in\Omega$ and $z\in N$
        \begin{equation}\label{eq:skewprodformula}
        F_{\bs\omega,\varepsilon}=T^{\omega_0}_\varepsilon\circ (S_{\iota_0,\varepsilon}+R_{\bs\omega,\varepsilon})
        \end{equation}
    with $T_\varepsilon$ as in \eqref{eq:innerdynamicsdefn} being of class $ C^{\min\{r,\kappa+1\}}$, $\{S_{i,\varepsilon}\}_i$ as in \eqref{eq:scattmapdynamicsdefn} being of class $ C^{\min\{r,\kappa\}}$ and $R_{\bs\omega,\varepsilon}: N\to N$ satisfying 
    \begin{equation}\label{eq:Romega estimate C1}
    |R_{\bs \omega,\varepsilon}|_{C^3}\leq  C  \bar\lambda^{\min\{\omega_{-1},\omega_{0}\}}.
    \end{equation}
    \item \textit{(Weak dependence):} For any $n\geq 0$, given any  $\bs\omega,\bs\omega'\in \Omega$ with $\bs\omega_k=\bs\omega_k'$ for $k<n$ and $\omega_n\leq \omega_n'$ (resp. $\bs\omega_k=\bs\omega_k'$ for $k\geq -n$ and $\omega_{-n-1}\leq \omega'_{-n-1}$)  we have 
    \begin{equation}\label{eq: comparing Romega ad inv manifolds}
    |R_{\bs\omega,\varepsilon}-R_{\bs\omega',\varepsilon}|_{C^3}\leq C  \bar\lambda^{\omega_0+\dots+\omega_{n}}\qquad\qquad (\text{resp. }|R_{\bs\omega,\varepsilon}-R_{\bs\omega',\varepsilon}|_{C^3}\leq C  \bar\lambda^{\omega_{-1}+\dots+\omega_{-n-1}}).
    \end{equation}
    \item \textit{(Invariant natural measure):} In the case where $X$ is a symplectic manifold and $F_\varepsilon$ is a parametric family of symplectomorphisms, the set $\mathcal X_\varepsilon$ is a symplectic lamination and any natural invariant measure $\mathfrak m_\varepsilon$ is invariant for the induced dynamics $\mathcal R_{\varepsilon}|_{\mathcal X_\varepsilon}:\mathcal X_\varepsilon\to \mathcal X_\varepsilon$.
    
    \end{enumerate}
\end{prop}

Although rather technical, the proof of Proposition \ref{prop:mainnhim} is  standard and can be found, for instance, in \cite{paradela2026robustlytransitivebehaviorsymplectic} (see Proposition 3.8. therein). Before proceeding, a couple of remarks are in order.
\medskip

First, notice that the alphabet $\Omega$ is countable (not finite). To obtain a lamination in which the base is homeomorphic to $\Omega$, we study the return map to an open set which contains the homoclinic channels $\{\Gamma_i\}$ in its closure. The point of doing that instead of considering a lamination close but uniformly away from the channels (and hence with base defined over a sequence space on a finite alphabet), is the following. By allowing the symbols of $\omega$ to be arbitrarily large, the corresponding orbit comes arbitrarily close to the normally hyperbolic manifold $\Lambda_\varepsilon$. As we will see below, this is crucial for shadowing certain pseudo-orbits of the iterated function system generated by the map $T_\varepsilon$ and the maps $\{S_{i,\varepsilon}\}_{i=1}^m$.
\medskip

Second, the skew-product \eqref{eq:skewprodnhim} is not locally-constant due to the presence of the term $R_{\bs \omega,\varepsilon}$ in \eqref{eq:skewprodformula}. However, as $\min\{\omega_{-1},\omega_0\}\to \infty$ (i.e. as we approach a homoclinic channel), the size of this term decreases exponentially fast. In \eqref{eq: comparing Romega ad inv manifolds} we compare the remainder along a unstable and stable manifolds respectively. These estimates will prove crucial to study certain holonomy maps between different leaves in $\mathcal X_\varepsilon$ (see Section \ref{sec: from IFS to skewprod}). \medskip

Finally, let us make some observations about the  bunching plus pinching hypothesis (i.e. taking the exponent $\kappa$ in \ref{it:H1}-\ref{it:H1+1/2} sufficiently large). The existence of $\mathcal X_\varepsilon$ can be deduced simply taking $\kappa=0$ (see \cite{GelfreichTuraev}). It is only in establishing the asymptotic formula \eqref{eq:skewprodformula} and the estimates for $R_{\bs\omega,\varepsilon}$ that we need to increase its value. We have not tried to optimize this lower threshold.

\begin{rem}\label{rem: barlambda}
    The quantifier $\bar\lambda$ in Proposition \ref{prop:mainnhim} is of the form 
    \[
    \bar\lambda=\max_{0\leq \varepsilon\leq \varepsilon_0} \{\alpha^{\kappa_0}_\varepsilon \lambda_\varepsilon, \bar\alpha^{\kappa_0}_\varepsilon\lambda_\varepsilon\}
    \]
    with $\alpha_\varepsilon,\lambda_\varepsilon$ as in \ref{it:H1}, $\bar\alpha_\varepsilon$ as in \ref{it:H1+1/2} and $\kappa_0$ being some positive universal exponent (it does not depend on any of the inputs in Proposition \ref{prop:mainnhim}) which can be computed explicitly from the proof of Proposition 3.8 in \cite{paradela2026robustlytransitivebehaviorsymplectic}.
\end{rem}

\begin{rem}\label{rem:N two dimensional}
    In Proposition \ref{prop:mainnhim} we have not made any assumption about the dimension of the normally hyperbolic manifold $\Lambda$ in \ref{it:H1}. In the following we will restrict to the case where $\Lambda$ (and hence the manifold $N$ in its parametrization) is two-dimensional.
\end{rem}

\subsection{Local ergodicity of iterated function systems}\label{sec:IFSintro}

In order to prove that the skew-product dynamics \eqref{eq:skewprodnhim} is ergodic, the first natural step is to study the IFS generated by the symplectic diffeomorphisms $T_\varepsilon$ in \eqref{eq:innerdynamicsdefn} and $\{S_{i,\varepsilon}\}_{i=1}^m$ in \eqref{eq:scattmapdynamicsdefn} of a smooth compact Riemannian surface $N$ (see Remark \ref{rem:N two dimensional}).  We will denote by 
\begin{itemize}
    \item $\mathrm{dist}_N(\cdot,\cdot):N\times N\to \mathbb R$ the distance,
    \item $\rho\in \mathcal M(N)$ the normalized area form, i.e. $\rho(N)=1$,
\end{itemize}
induced by the Riemannian metric on $N$.
\medskip

This section is organized as follows. First,  we introduce a semi-local source of ergodicity which we call Hopf boxes. Theorem \ref{thm:main} shows that Hopf boxes (and hence local ergodicity) can be found in families of IFS's which we call $\epsilon$-good and whose distinguishing property is that they display a small region for which one can find two elements of the IFS which are $O_{C^3}(\epsilon)$ close to a suitably chosen pair of hyperbolic affine maps which exhibit good covering properties (see \eqref{eq:affinemaps}). 
\medskip

Second, we give a collection of (very) easily verifiable conditions (see Assumptions \ref{it:itemB1}-\ref{it:itemB5}) which guarantee that a certain IFS generated by a \textit{given} one-parameter family $\{T_\varepsilon,S_\varepsilon\}_\varepsilon$ is $\epsilon$-good, hence locally ergodic (on some small subset of positive measure).
\medskip

\subsection*{Notation for IFS}
In the following, given a collection $\{G_i\}_{i\in I}$ of $C^2$ diffeomorphisms of a compact Riemannian surface $N$ we denote by 
\[
\langle \{G_i\}_{i\in I}\rangle \subset \mathrm{Diff}^3(N)
\]
the group generated by $\{G_i\}_{i\in I}$. Let $n\in\mathbb N$. Given a word $w=(w_0,\dots,w_{n-1})\in I^n$ we write 
\[
G_w^{n}=G_{w_{n-1}}\circ \cdots \circ G_{w_0}\in\langle \{G_i\}_{i\in I}\rangle\qquad\qquad\text{and}\qquad\qquad G_w^{-n}=(G_w^n)^{-1}\in \langle \{G_i\}_{i\in I}\rangle.
\]
In addition, given a sequence $\omega=(\omega_1,\dots,\omega_n,\dots)\in I^\mathbb N$ we denote by 
\[
G_\omega^n=G_{\omega_{n-1}}\circ \cdots \circ G_{\omega_0}\in\langle \{G_i\}_{i\in I}\rangle\qquad\qquad\text{and}\qquad\qquad G_\omega^{-n}=(G_\omega^n)^{-1}\in \langle \{G_i\}_{i\in I}\rangle.
\]
Given a measurable set $E\subset N$ we say that it is $\langle \{G_i\}_i\rangle$-invariant provided 
\[
g(E)\subset E\qquad\qquad \text{for all}\quad g\in \langle \{G_i\}_i\rangle.
\]
\subsection*{Hopf boxes and local ergodicity for IFS}\label{sec:Hopf boxes}
For $i=1,2$ and $\lambda\in(1/2,1)$ consider the parametric family of affine maps 
\begin{equation}\label{eq:affinemaps}
F_{\mathrm{aff},i,\lambda}:z=(x,y)\in \mathbb R^2 \mapsto (f_{x,i,\lambda}(x),f_{y,i,\lambda}(y))\in \mathbb R^2
\end{equation}
with 
\[
f_{x,1,\lambda}(s)=\lambda s=f_{y,1,\lambda}^{-1}(s)
\]
and
\[
f_{x,2,\lambda}(s)=1-\lambda(1-s)=f_{y,2,\lambda}^{-1}(s).
\]
\begin{figure}
    \centering
\includegraphics[scale=0.65]{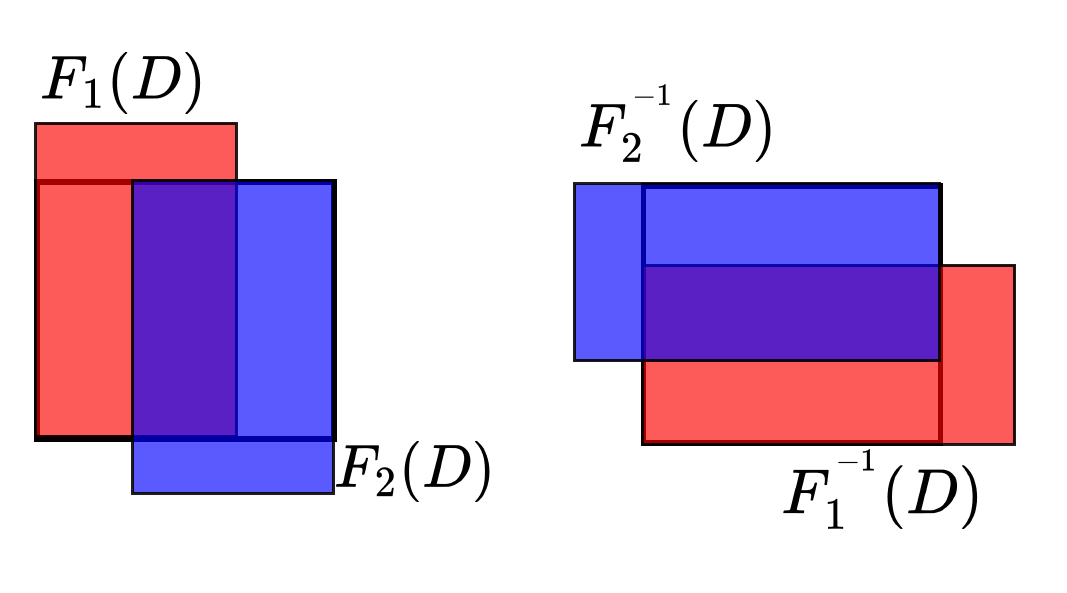}
    \caption{The action of the maps $F_{1,\lambda},F_{2,\lambda}$ on $D=[0,1]^2$ for some $\lambda\in(1/2,1)$.}
    \label{fig:fig3}
\end{figure}
\noindent A sketch  of the action of these maps on the unit square can be found in Figure \ref{fig:fig3}. For later convenience we introduce the following notation\footnote{The domain $\mathcal D$ is (arbitrarily) chosen among the set of neighbourhoods of $D$.}
\[
D=[0,1]^2\subset\mathbb R^2\qquad\qquad \mathcal D=[-2,2]^2\subset\mathbb R^2
\]
and highlight the bi-covering property
\begin{equation}\label{eq:bicoveringproperty}
D\subset \bigcup_i F_{\mathrm{aff},i,\lambda}(D)\qquad\qquad\text{and}\qquad\qquad D\subset \bigcup_i F_{\mathrm{aff},i,\lambda}^{-1}(D).
\end{equation}
The following class of maps will be the central object of our study.

\begin{defn}[$\epsilon$-good family]\label{defn:epsilon good}
    Let $\epsilon>0$ and let $J\subset (1/2,1)$ be a closed interval. We say that a  parametric family  $\{F_{1,\lambda},F_{2,\lambda}\}_{\lambda\in J}\subset \mathrm{Diff}^3(N)$ acting on a smooth surface  is $\epsilon$-\textit{good} if for all $\lambda\in J$ there exists a $C^3$ local coordinate chart $\phi_\lambda:\mathcal D\subset\mathbb R^2\to N$ such that 
    \begin{equation}\label{eq:epsilon}
\max_{i=1,2}\  \lVert \widetilde F_{i,\lambda}-F_{\mathrm{aff},i,\lambda} \rVert_{C^3(\mathcal D\times J)}\leq \epsilon,
\end{equation}
where $\widetilde F_{i,\lambda}=\phi^{-1}_\lambda \circ F_{i,\lambda}\circ\phi_\lambda$.
\end{defn}

For an $\epsilon$-good family with $0\ll \epsilon\ll 1$ much of the local affine geometry can be continued smoothly into the nonlinear dynamics. In particular, if we fix any closed interval $J\subset(1/2,1)$,  a straightforward application of the stable manifold theorem (see \cite{MR1326374,MR1963683}) shows that there exists $\epsilon_0(J)>0$ such that for any $0\leq \epsilon\leq \epsilon_0$ the map $\widetilde F_{i,\lambda}$ possesses a unique hyperbolic fixed point $z_{i,\lambda}\in \mathcal D$ with local  un/stable manifolds $W^{u,s}(z_{i,\lambda})\subset \mathcal D$. These are, moreover,  $O_{C^3}(\epsilon)$-close to the corresponding un/stable manifolds of the affine system \eqref{eq:affinemaps}.
\medskip

\begin{figure}
    \centering
\includegraphics[scale=0.65]
{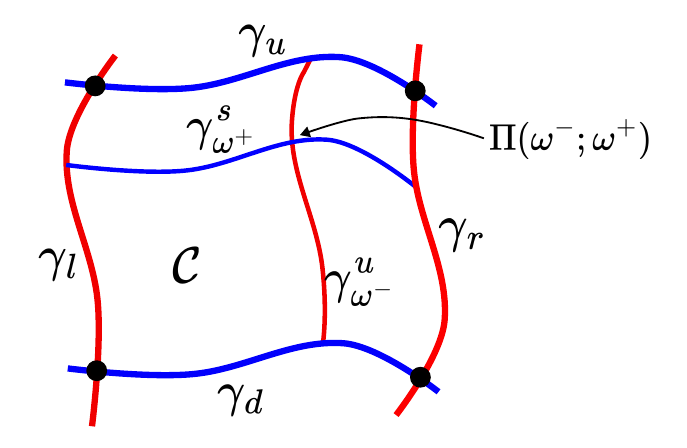}
    \caption{The Hopf box $\mathcal C$ and the coding map $\Pi$.}
    \label{fig:fig4}
\end{figure}

We denote by
\begin{itemize}
    \item $\gamma_{d,\lambda}\subset W^s(z_{1,\lambda})$ the right branch of the stable manifold of $z_{1,\lambda}$.
     \item $\gamma_{r,\lambda}\subset W^u(z_{2,\lambda})$ the lower branch of the unstable manifold of $z_{2,\lambda}$.
      \item $\gamma_{u,\lambda}\subset W^s(z_{2,\lambda})$ the left branch of the stable manifold of $z_{2,\lambda}$.
       \item $\gamma_{l,\lambda}\subset W^u(z_{1,\lambda})$ the upper branch of the unstable manifold of $z_{1,\lambda}$.
\end{itemize}
After decreasing (if needed) the value of $\epsilon_0(J)$, the curves $\gamma_{d,\lambda}$ and $\gamma_{r,\lambda}$ intersect transversally. The same holds for   $\gamma_{u,\lambda}$ and $\gamma_{l,\lambda}$ (see Figure \ref{fig:fig4}).

\begin{defn}[Hopf box]\label{defn:hopfbox}
    Fix any closed subset $J\subset (1/2,1)$ and let $0< \epsilon\leq \epsilon_0(J)$. Then, for a $\epsilon$-good family we call \textit{Hopf box} to the subset 
    \[
    \mathcal C_\lambda=\phi_\lambda(D_\lambda)\subset N
    \]
    where $D_\lambda\subset \mathbb R^2$ is the curvilinear rectangle bounded by the four curves $\gamma_{d,\lambda},\gamma_{r,\lambda},\gamma_{u,\lambda}$ and $\gamma_{l,\lambda}$.
\end{defn}

The following  result  motivates the above definition.

\begin{thmx}[Local ergodicity]\label{thm:main}
  Let $N$ be a smooth compact Riemannian surface and denote by  $\rho$ the associated area form. Consider  a $\epsilon$-good family of $C^3$  maps $\{F_{1,\lambda},F_{2,\lambda}\}_{\lambda}:N\to N$ and denote by 
  \begin{equation}\label{eq:the group action}
\mathcal G_\lambda=\langle F_{1,\lambda},F_{2,\lambda}\rangle\subset \mathrm{Diff}^3(N).
  \end{equation}
  Then, there exists a non-empty open interval $J_0\subset (1/2,1)$ such that for any $\epsilon$ small enough and  for $\mathrm{Leb}$-a.e. $\lambda\in J_0$, the group action $ (\mathcal G_\lambda,N,\rho)$
is locally ergodic on its corresponding Hopf box $\mathcal C_\lambda\subset N$. That is, for every measurable set $E$ which is $\mathcal G_\lambda$-invariant we have 
    \[
    \rho_{\mathcal C_\lambda}(E\cap \mathcal C_\lambda)\in\{0,1\}.
    \]
\end{thmx}
\medskip

Theorem \ref{thm:main} is the most relevant novelty in our construction. Its proof  occupies Sections \ref{sec:geometric preli}-\ref{sec:equivalenceLebesgue} (an outline is presented in Section \ref{sec:outline}). Note that we do not assume the maps $F_{1,\lambda},F_{2,\lambda}\in \mathrm{Diff}^3(N)$ to preserve the area form $\rho$ or to display a smooth common invariant measure. This will be important in our application of Theorem \ref{thm:main} to the proof of Theorem \ref{thm:mainham}. Before continuing our roadmap of the proof of Theorem \ref{thm:mainham} it is worth to first give a quick overview on the proof of Theorem \ref{thm:main} and also explain why, with our current approach we only obtain an almost everywhere (in $\lambda$) result.
\medskip

A very rough idea of the mechanism is the following.  We adopt the notation in Definition \ref{defn:epsilon good}.   We will see below that the mappings $\{ F_{i,\lambda}\}_i$  preserve (locally on $\mathcal C_\lambda $) uniform stable  and unstable  cones. Hence, for each $\omega\in \{1,2\}^\mathbb N$, a standard graph transform argument shows that
the sets (recall the notation for compositions introduced above)
\[
V_{\omega|n,\lambda}=\mathcal C_\lambda\cap F^{n}_{\omega,\lambda}(\mathcal C_\lambda) \qquad\qquad H_{\omega|n,\lambda}= F^{-n}_{\omega,\lambda}(\mathcal C_\lambda) \cap \mathcal C_\lambda
\]
constitute a nested sequence of unstable and stable rectangles satisfying 
\[
F_{\omega,\lambda}^n(H_{\omega|n,\lambda})=V_{\omega|n,\lambda}\qquad\qquad F_{\omega,\lambda}^{-n}(V_{\omega|n,\lambda})=H_{\omega|n,\lambda}.
\]
In particular, for each $\omega^-=(\dots,\omega_{-n},\dots,\omega_{-1})\in\{1,2\}^\mathbb N$  (resp. $\omega^+=(\omega_0,\dots,\omega_n,\dots)\in\{1,2\}^\mathbb N$) 
\[
\gamma_{\omega^-,\lambda}^u=\bigcap_n V_{(\omega_{-n},\dots,\omega_{-1}),\lambda}\qquad\qquad (\text{resp. } \gamma_{\omega^+,\lambda}^s=\bigcap_n H_{(\omega_{0},\dots,\omega_{n-1}),\lambda})
\]
is a unstable (resp. stable) $C^2$ curve which fully crosses $\mathcal C_\lambda$ (see Figure \ref{fig:fig4}). Thanks to the bi-covering property \eqref{eq:bicoveringproperty}
\[
\bigcup_{\omega^-}\gamma_{\omega^-,\lambda}^u=\mathcal C_\lambda=\bigcup_{\omega^+}\gamma_{\omega^+,\lambda}^s
\]
Moreover, for each pair $\omega^-,\omega^+$ the intersection $\gamma_{\omega^-,\lambda}^u\cap \gamma_{\omega^+,\lambda}^s$ consists of a single point. We however note that neither $\{\gamma_{\omega^-,\lambda}^u\}_{\omega^-}$ nor $\{\gamma_{\omega^+,\lambda}^s\}_{\omega^+}$ form a foliation of $\mathcal C_\lambda$ as different leaves might have nonempty intersection. Still, we can run a modified version of Hopf's argument. To do so we proceed as follows. 
\medskip

Denote by $\Omega=\{1,2\}^\mathbb Z$. We first introduce the coding map 
    \begin{equation}\label{eq:affinetwosidedcodingmap}
    \Pi_\lambda:\omega=(\omega^-;\omega^+)\in \Omega \mapsto \gamma_{\omega^-,\lambda}^u\cap \gamma_{\omega^+,\lambda}^s\in \mathcal C_\lambda
    \end{equation}
and observe that if $E$ is an invariant subset for $\mathcal G_\lambda$, then 
\[
\widehat E=\Pi^{-1}_\lambda (E\cap \mathcal C_\lambda)
\]
is invariant for the full shift $\sigma:\Omega\mapsto \Omega$.  In view of this observation, we would like to construct a measure $\mu$ on the probability space $\Omega$ (equipped with the Borel sigma algebra) which satisfies:
\begin{enumerate}
    \item \label{it:1}\textit{(Projection equivalence):} The pushforward of $\mu$ under the coding map  is equivalent to the area measure $\rho$ restricted to $D$.
    \item \label{it:2} \textit{(Ergodicity):} $\mu$ is $\sigma$-invariant and ergodic for the full shift.
\end{enumerate}
\medskip

Roughly speaking, \eqref{it:1} ensures that ergodicity of the group action $(\mathcal G_\lambda, D, \rho_D)$ is implied by ergodicity of $(\sigma,\Omega,\mu)$ (which we assumue in \eqref{it:2}).
\medskip

While the shift admits plenty of invariant ergodic measures,  the more demanding step in this program is  to guarantee that the  projection of this measure is equivalent to the area form $\rho$. To do so, we deploy ideas coming from the theory of Bernoulli convolutions which, at the current stage, only allow us to deduce a.e. statements (see Section \ref{sec:typical AC affine} for a more detailed discussion).

\subsection*{Transversality-torsion and Hopf boxes}

We now introduce a suitable class of parametric families of diffeomorphisms of a $2$-dimensional manifold  and give sufficient conditions which guarantee $\epsilon$-goodness. The conditions are local (i.e. only require control of the dynamics on a small region of the phase space) and,  as we will see below, are general enough so that, in the context of Theorem \ref{thm:mainham}, we can check that the corresponding inner map \eqref{eq:innerdynamicsdefn} and system of scattering maps \eqref{eq:scattmapdynamicsdefn} verify these hypothesis.
\medskip

To be precise, we consider a  pair of (parametric families of) $C^3$ maps $T_\varepsilon,S_\varepsilon$ of a compact surface $N$ for which there exists a local coordinate patch 
\begin{equation}\label{eq:local coordinate patch annulus}
\phi_\varepsilon:\mathbb A\to \mathcal A:=\phi_\varepsilon(\mathbb A)\subset N
\end{equation}
 on which the maps take the form 
\begin{equation}\label{eq:innermapsymbolic}
\widetilde T_\varepsilon:= \phi_\varepsilon^{-1}\circ T_\varepsilon\circ\phi_\varepsilon:\binom{\varphi}{J}\mapsto \binom{\varphi+\beta+a_\varepsilon J+R_{\varepsilon,\varphi}(\varphi,J)}{J+R_{\varepsilon,J}(\varphi,J)}
\end{equation}
and
\begin{equation}\label{eq:scattmapsymbolic}
\widetilde S_\varepsilon:= \phi_\varepsilon^{-1}\circ S_\varepsilon\circ\phi_\varepsilon:\binom{\varphi}{J}\mapsto \binom{\varphi+\varepsilon P_{\varepsilon,\varphi}(\varphi,J)}{J+\varepsilon P_{\varepsilon,J}(\varphi,J)}.
\end{equation}
We will suppose moreover that
\begin{enumerate}[label={[B\arabic*]}]
\item \label{it:itemB1} \textit{Rotation vector of the map $T_\varepsilon$:} The vector $\beta\in\mathbb R$ is of constant type, i.e. there exists $\gamma>0$ such that 
\begin{equation}\label{eq:constanttype}
\beta\in \mathrm{DC}_\gamma=\{\alpha\in\mathbb R\colon |\alpha q-p|\geq \gamma |q|^{-1}\ \text{ for all }q\in\mathbb Z\setminus\{0\}\text{ and all }p\in\{0,\dots,|q|\}\}.
\end{equation}
   \item \label{it:itemB2} \textit{Small nonlinearity of the map $T_\varepsilon$}: The functions $R_{\varepsilon,\varphi},R_{\varepsilon,J}$ satisfy 
   \[
   R_{\varepsilon,\varphi}(\varphi,J)= J^2 U_{\varepsilon}(\varphi,J),\qquad\qquad R_{\varepsilon,\varphi}(\varphi,J)= J^3 V_{\varepsilon}(\varphi,J),
   \]
   for some functions $U_\varepsilon,V_\varepsilon\in C^3(\mathbb A)$. 
   \item \label{it:itemB3} \textit{Transversality-torsion}: The function $P_{\varepsilon, J}$ satisfies
\[
P_{0,J}(0,0)=0\qquad\qquad \text{and}\qquad\qquad b=\partial_\varphi P_{0,J}\neq 0.
\] 
Moreover, 
\[
ab>0\qquad\text{where}\qquad a=a_0.
\]
\item \textit{Regularity with respect parameters:} \label{it:itemB5} There exists  $\varepsilon_0>0$ such that $R_{\varphi,\star}$ and $P_{\varphi,\star}$  ($\star=\varphi,J$) depend $C^3$ on $\varepsilon$ for $-\varepsilon_0\leq \varepsilon\leq \varepsilon_0$.
\end{enumerate}

\begin{rem}
In the following we will sometimes abuse notation and simply say that the maps $T_\varepsilon,S_\varepsilon$  satisfy assumptions  \ref{it:itemB1}-\ref{it:itemB5} without referring to the local coordinate patch \eqref{eq:local coordinate patch annulus}.
\end{rem}

\begin{figure}
    \centering
\includegraphics[scale=0.65]
{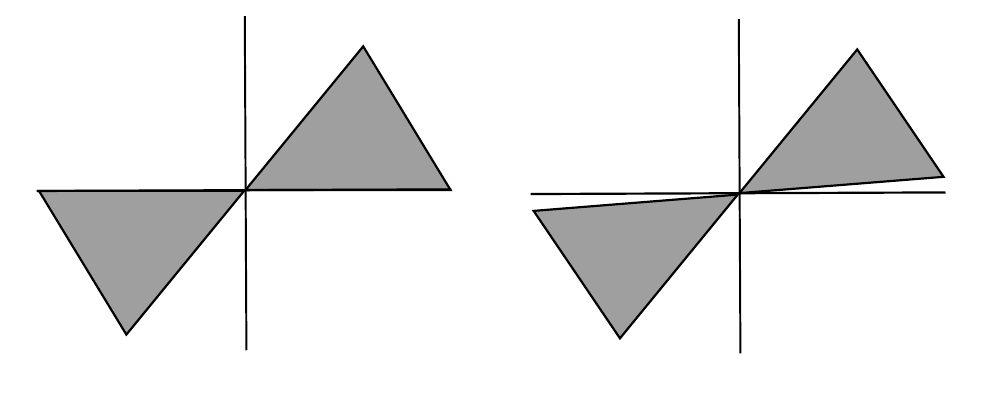}
    \caption{Here we sketch how the transversality torsion works. On the left we depict in grey the image of the cone formed by the first and second quadrants under the differential of the map $T_\varepsilon^n$ (for any $n\in \mathbb N$) evaluated at $\varphi=0=J$. Notice that this action is not hyperbolic because the image is not strictly contained inside the cone. On the right we show the image of the cone under the differential of the map $T_\varepsilon^n\circ S_\varepsilon$ evaluated at $\varphi=0=J$. For $\varepsilon>0$ this action is hyperbolic as the cone gets mapped to its interior.}
    \label{fig:fig5}
\end{figure}

We now relate this set of assumptions to $\epsilon$-goodness. To do so we denote by 
\[
K=\max_{0\leq\varepsilon\leq \varepsilon_0}\{|\partial^2_\varepsilon a_\varepsilon|,\lVert U_{\varepsilon}\rVert_{C^3}, \lVert V_{\varepsilon}\rVert_{C^3},\lVert P_{\varepsilon}\rVert_{C^3}, \lVert \partial_\varepsilon R_{\varepsilon}\rVert_{C^2},\lVert \partial_\varepsilon P_{\varepsilon}\rVert_{C^2}, \lVert \partial^2_{\varepsilon^2} R_{\varepsilon}\rVert_{C^1}, \lVert \partial^2_{\varepsilon^2} P_{\varepsilon}\rVert_{C^1} \}.
\]
\begin{prop}\label{prop:TT and HB}
   Let $T_\varepsilon,S_\varepsilon\in \mathrm{Diff}^3(N)$ satisfy assumptions \ref{it:itemB1}-\ref{it:itemB5} and let  $J\subset(1/2,1)$ be any fixed closed interval.  There exists
\[
\epsilon_0(\gamma,ab,K)>0
\]
 and a  $C^3$ function 
\[
\lambda\in J\mapsto s(\lambda)\in (0,1)
\]
such that for any $0<\epsilon\leq \epsilon_0$, the group $\langle \{T_{\varepsilon(\epsilon, \lambda)},S_{\varepsilon(\epsilon, \lambda)}\}\rangle_{\lambda \in J}$ with 
\[
\varepsilon(\epsilon,\lambda)=\epsilon^3 s(\lambda)
\]
contains a pair of words $\{F_{1,\lambda},F_{2,\lambda}\}_{\lambda\in J}$ which are $\epsilon$-good.  Moreover, the corresponding Hopf box $\mathcal C_\lambda\subset N $ (recall the definition \eqref{defn:hopfbox}) contains an open  ball of radius greater than $\varepsilon^{3/2}$.
    \end{prop}
The proof of Proposition \ref{prop:TT and HB} follows closely the construction in Section 3 of \cite{paradela2026robustlytransitivebehaviorsymplectic} (see also \cite{guardia2025partiallyhyperbolicdynamics3body}) and is presented in Section \ref{sec:proof TT }. The basic ideas are the following. On one hand, under \ref{it:itemB3} (and \ref{it:itemB2}), it is possible to create hyperbolicity by combining the torsion of the map \eqref{eq:innermapsymbolic} with the transversality of the map \eqref{eq:scattmapsymbolic} (see figure \ref{fig:fig5}).  This is the so-called transversality-torsion mechanism  (see for instance \cite{MR1949441}). On the other hand, under \ref{it:itemB1}, the strong irrationality of the rotation vector $\beta$, guarantees that among the hyperbolic words reproducing the transversality-torsion mechanism it is possible to choose a finite collection whose fixed points approximate that of the affine model \eqref{eq:affinemaps}.
\medskip

\subsection{Global ergodicity of IFS}\label{sec:globalergodicity}
In some applications (as the one which concerns us at the moment, i.e. proving Theorem \ref{thm:mainham}), the maps $F_{1,\lambda},F_{2,\lambda}$ form part of a  possibly larger collection $\{F_{i,\lambda}\}_{i=1}^k:N\to N$, $k\geq 2$ acting on a compact surface. In this case,  if the Hopf box $\mathcal C_\lambda$ associated to $F_{1,\lambda},F_{2,\lambda}$ is \textit{reachable} we can upgrade the result in Theorem \ref{thm:main} to global ergodicity. Namely,  we say that  $\{F_{i,\lambda}\}_{i=1}^k:N\to N$ satisfies Property \ref{it:Assumption1} if:

\begin{enumerate}[label={[R]}]
    \item \label{it:Assumption1} There exists a subset $R_\lambda\subset N$ of full $\rho$-measure such that for any $z\in R_\lambda$ there exists 
    \[
    g\in \langle\{F_{i,\lambda}\}_{i=1}^k\rangle
    \]
    such that $g(z)\in \mathcal C_\lambda$ where $\mathcal C_\lambda\subset N$ is the corresponding Hopf box for the maps $\{F_{i,\lambda}\}_{i=1}^2$.
    \end{enumerate}

\begin{prop}\label{prop:global ergodicity}
    Consider the setting of Theorem \ref{thm:main}  and assume moreover that the parametric family  $\{F_{1,\lambda},F_{2,\lambda}\}_\lambda\subset \mathrm {Diff}^3(N)$ therein forms part of a (not necessarily larger) collection of maps $\{F_{i,\lambda}\}_{i=1}^{k}\in \mathrm{Diff}^3(N)$ which satisfies Property \ref{it:Assumption1}. Denote by 
    \[
    \mathcal G_\lambda=\langle \{F_{i,\lambda}\}_{i=1}^{k}\rangle\subset \mathrm{Diff}^3(N).
    \]
    Then, for $\mathrm{Leb}$-a.e. $\lambda\in J_0$ the group action $(\mathcal G_\lambda,N,\rho)$ is ergodic.  That is, for every measurable set $E$ which is $\mathcal G_\lambda$-invariant we have 
    \[
    \rho (E)\in\{0,1\}.
    \]
\end{prop}

\begin{proof}
    Let $E$ be $\mathcal G_\lambda$-invariant. From Theorem \ref{thm:main} we have 
    \[
    \rho(E\cap \mathcal C_\lambda)=0\qquad\qquad\text{or}\qquad\qquad \rho( E\cap  \mathcal C_\lambda)=\rho( \mathcal C_\lambda).
    \]
    On one hand, $\rho (E\cap  \mathcal C_\lambda)=0$ clearly forces $\rho(E)=0$ since 
    \[
    E\subset \bigcup_g g^{-1}(E\cap  \mathcal C_\lambda) \qquad\qquad \mathrm{mod }\  \rho
    \]
    and each $g$ is non-singular with respect to $\rho$. On the other hand, $\rho(E\cap  \mathcal C_\lambda)=\rho( \mathcal C_\lambda)$ implies that $\rho(E^c\cap  \mathcal C_\lambda)=0$. Hence, applying the same argument we deduce that $\rho(E)=1$.
\end{proof}
\medskip

\subsection*{Property \ref{it:Assumption1}  for the scattering maps}
Propositions \ref{prop:TT and HB} and \ref{prop:global ergodicity} reduce the problem of global ergodicity of the IFS generated by $\{T_\varepsilon,\{S_{i,\varepsilon}\}_{i=1}^m\}$ (recall that these are the maps in \eqref{eq:innerdynamicsdefn}-\eqref{eq:scattmapdynamicsdefn}) to the question of whether the Hopf box $\mathcal C_\lambda$ in Proposition \ref{prop:TT and HB} is ``reachable'' from anywhere in (a full measure subset of) the manifold $N$. Our next result shows that this is the case under a mild condition on the system of scattering maps. 

As we did above, we state this condition in abstract terms and then relate its assumptions to the system of scattering maps \eqref{eq:scattmapdynamicsdefn} which arise in the setting of Proposition \ref{prop:mainnhim}. More precisely, given $m\geq 2$, and a collection of (parametric families of) maps  $\{S_{i,\varepsilon}\}_{i=1}^m:N\to N$ satisfies:
\begin{enumerate}[label={[B\arabic*]},start=5]
    \item \label{it:B6} \textit{(Hörmander's  condition):} Fix some $n \geq 0$ and $\tilde r\geq n$. There exists $C>0$,  $\varepsilon_0>0$ and a covering $\{U_j\}_j$ of $N$ by charts $\psi_j:U_j\subset \mathbb R^{2} \to N$ such that, for any $0\leq \varepsilon\leq \varepsilon_0$, on each $z\in U_j$
    \[
    \psi_j^{-1}\circ S_{i,\varepsilon}\circ\psi_j(z)=z+\varepsilon X_i^{(j)}(z)+\varepsilon^2 R_i^{(j)}(z;\varepsilon)\qquad\qquad i=1,\dots,m
    \]
    for some $X_i^{(j)},R_i^{(j)}:U_j\to \mathbb R^{2}$ such that $\lVert X_i^{(j)}\rVert_{C^{\tilde r}},\lVert R_i^{(j)}\rVert _{C^{\tilde r}}\leq C$ and at any $z\in U_j$
    \begin{equation}\label{eq:assumptionLiedistribution}
        \mathrm{Lie}^{(n,\tilde r)}_z\{X_1^{(j)},\dots, X_m^{(j)}\}=\mathbb R^{2}
    \end{equation}
    where, for each $1\leq s\leq \tilde r$, given $m$ vector fields $\{Y_i\}_{i=1}^m:U_j\to \mathbb R^{2}$ 
    \begin{equation}\label{eq:Liedistributions}
    \mathfrak L^{(s)}\{Y_1,\dots, Y_m\}
=
\{[Y_i,Y]:1\leq i\leq m,\ Y\in\mathfrak L^{(s-1)}\{Y_1,\dots, Y_m\}\}
    \end{equation}
   $\mathfrak L^{(0)}\{Y_1,\dots, Y_m\}=\{Y_1,\dots, Y_m\}$  and 
    \[
    \mathrm{Lie}^{(n,\tilde r)}_z\{Y_1,\dots, Y_m\}= \mathrm{span}\{Y(z)\in \bigcup_{s=n}^{\tilde r}  \mathfrak L^{(s)}_z\{Y_1,\dots, Y_m\}\}
    \]
\end{enumerate}

Hypothesis \ref{it:B6} with $n=0$ is nothing but a suitably adapted reformulation of those in the classical Ball-Box Theorem in sub-Riemannian geometry. It simply says that the Lie distribution of order $\tilde r\in \mathbb N$ spans $\mathbb R^2$. From an abstract viewpoint it  can be understood as a sufficient condition on a collection of close to identity diffeomorphisms of a compact surface $N$ which guarantees $\varepsilon$-transitivity of the corresponding IFS (see \cite{paradela2026robustlytransitivebehaviorsymplectic} for a more detailed discussion). For our purposes (giving a proof of Theorem \ref{thm:mainham}), assumption \ref{it:B6} will be understood jointly with assumptions \ref{it:H1}-\ref{it:H2}. Hence, the regularity class $\tilde r$ in assumption \ref{it:B6} will be determined by that in the strong spectral gap condition in assumption \ref{it:H1}.
\medskip

For $n\geq 1$ instead of looking at the whole distribution we only consider brackets of size larger than $n$. The advantage of doing so is that it implies $\varepsilon^{n+1}$-transitivity of the corresponding IFS.

\begin{prop}\label{prop:IFS property R}
    Suppose  that $\{S_{i,\varepsilon}\}_{i=1}^m$ satisfy \ref{it:B6} for some $n\geq 0$. Then, there exists $K>0$ and  
\[
\varepsilon_0(\{X_i^{(j)}\}, \lVert R_i^{(j)}\rVert_{C^{\tilde r}})>0
\]
 such that for every $z,z_*\in N$ there exists $g\in \langle \{S_{i,\varepsilon}\}_{i=1}^m\rangle$ such that 
 \[
 \mathrm{dist}_N(g(z),z_*)\leq K\varepsilon^{n+1}.
 \]
\end{prop}

Choosing now $n\geq 1$ we obtain global ergodicity of the IFS generated by $T_\varepsilon$ and $\{S_{i,\varepsilon}\}_{i=1}^m$. Although we will not be able to use this result directly to establish ergodicity of the skew-product in \ref{prop:mainnhim}, we state it here since we believe it is of independent interest.

\begin{thm}\label{thm:ergodicity IFS final}
    Let $N$ be a smooth compact Riemannian surface and denote by  $\rho$ the associated area form. Let $T_\varepsilon,S_{1,\varepsilon}\in \mathrm{Diff}^3(N)$ satisfy assumptions \ref{it:itemB1}-\ref{it:itemB5} and $\{S_{i,\varepsilon}\}_{i=1}^m\subset \mathrm{Diff}^3(N)$ satisfy \ref{it:B6} with $n\geq 1$. Then, there exists $\varepsilon_0>0$ such that  for any $0<\varepsilon\leq \varepsilon_0$ the collection $\{T_\varepsilon,\{S_{i,\varepsilon}\}_{i=1}^m\}\subset \mathrm{Diff}^3(N)$ satisfies property \ref{it:Assumption1}. In particular, if we denote by
    \[
    \mathcal G_\varepsilon=\langle \{T_\varepsilon,\{S_{i,\varepsilon}\}_{i=1}^m\}\rangle \subset \mathrm{Diff}^3(N)
    \]
    then the group action $(\mathcal G_\varepsilon, N,\rho)$ is ergodic for almost every $\varepsilon\in (0,\varepsilon_0)$.
\end{thm}
\begin{proof}
    The fact that the collection $\{T_\varepsilon,\{S_{i,\varepsilon}\}_{i=1}^m\}\subset \mathrm{Diff}^3(N)$ satisfies \ref{it:Assumption1} is straightforward from Proposition \ref{prop:IFS property R} and the fact that the associated Hopf box $\mathcal B_{\lambda(\varepsilon)}$ has size $\varepsilon^{3/2}$.  Then, ergodicity for almost every $\varepsilon$ follows from Proposition \ref{prop:global ergodicity}. 
\end{proof}

\subsection{Ergodicity of the skew-product}\label{sec:skewprodintro}

Finally, we upgrade ergodicity of the IFS generated by the collection of maps $\{T_\varepsilon,\{S_{i,\varepsilon}\}_{i=1}^m\}$ satisfying \ref{it:itemB1}-\ref{it:B6}  to ergodicity of the associated skew-product in Proposition \ref{prop:mainnhim}.

  \begin{prop}\label{prop:Reduction Skew-prod to IFS}
Assume that $r>0$ and $\kappa>0$ are large enough. Let $X$ be a symplectic manifold,  let $F_\varepsilon:X\to X$ be a parametric family of $C^r$ symplectomorphisms  which satisfy \ref{it:H1}-\ref{it:H2} and let  $T_\varepsilon:N\to N$ and  $\{S_{i,\varepsilon}\}_{i=1}^m:N\to N$ be given by \eqref{eq:innerdynamicsdefn} and \eqref{eq:scattmapdynamicsdefn}. Let $\mathcal R_\varepsilon$ be the return map and $\mathcal X_\varepsilon$ be the symplectic lamination constructed in Proposition \ref{prop:mainnhim} and let $\mathfrak m_\varepsilon$ be any  natural invariant measure for $\mathcal R_\varepsilon|_{\mathcal X_\varepsilon}:\mathcal X_\varepsilon\to \mathcal X_\varepsilon$.  Suppose that there exists $\tilde\varepsilon_0>0$ such that for all $0<\varepsilon\leq \tilde\varepsilon_0$:
\begin{enumerate}
\item Each of the maps $T_\varepsilon, \{S_{i,\varepsilon}\}_{i=1}^m$ preserves a smooth measure on $N$.
    \item The maps $\{T_\varepsilon, S_{1,\varepsilon}\}$  verify hypothesis \ref{it:itemB1}-\ref{it:itemB5} for some $\gamma$ independent of $\varepsilon$.   \item The scattering maps $\{S_{i,\varepsilon}\}_{i=1}^m$ satisfy assumption \ref{it:B6} with $n\geq 1$ and some $1\leq \tilde r\leq r$. 
    
\end{enumerate}
Then, there exists $\varepsilon_0>0$ such that for Lebesgue almost every $\varepsilon\in (0,\varepsilon_0)$ the measure $(\Psi_\varepsilon^{-1})_* \mathfrak m_\varepsilon$ on $\Sigma\times N$  is invariant and ergodic for  the corresponding skew-product $\mathcal F_\varepsilon:\Sigma\times N\to \Omega\times N$ in \eqref{eq:skewprodnhim}. In particular, for Lebesgue almost every $\varepsilon\in (0,\varepsilon_0)$, the measure  $\mathfrak m_\varepsilon$ on $\mathcal X_\varepsilon$ is ergodic for $\mathcal R_\varepsilon|_{\mathcal X_\varepsilon}$.
\end{prop}

The proof of Proposition \ref{prop:Reduction Skew-prod to IFS} is given in Section \ref{sec: from IFS to skewprod}. The strategy of the proof combines elements from that in Theorem \ref{thm:ergodicity IFS final} and a suitable shadowing argument which allows us to reproduce orbits of the IFS generated by the system of scattering maps within the skew-product dynamics.

\begin{rem}
    The avid reader might have noticed that we stated  Proposition \ref{prop:Reduction Skew-prod to IFS} for symplectomoprhisms  instead of just considering diffeomorphisms of a Riemannian manifold as we did in the previous sections. The reason for doing so is that in the symplectic setting, invariance of the natural measure $\mathfrak m_\varepsilon$ is automatic (recall the discussion after Definition \ref{defn:natural measure}).

    We also note that from a  symplectic form one can always produce a Riemannian metric. Moreover, the  volume form associated to this metric coincides with that coming from the symplectic form.
\end{rem}
\medskip

\subsection{Typicality of the assumptions}\label{sec:genericityintro}
Finally, we observe that hypothesis \ref{it:H1}-\ref{it:H2} and \ref{it:itemB1}-\ref{it:B6} are typically satisfied for  one-parametric families $\{F_{\varepsilon,f}\}_\varepsilon:M\times N\to M\times N$ of Hamiltonian diffeomorphisms unfolding from a direct product $F_0$ as in Theorem \ref{thm:mainham}.

\begin{prop}\label{prop:genericityofassumptions}
  Let  $(M,\varpi_M)$ and $(N,\varpi_N)$ be real-analytic symplectic manifolds with $N$ being a compact surface. Let $F_0\in \mathrm{Ham}^\omega(M\times N)$ be a direct product of Hamiltonian diffeomorphisms  which satisfies the assumptions \ref{it:Hyp1}-\ref{it:Hyp2} with $\kappa$ sufficiently large (independent of $F_0$) and  displays a collection of full disjoint homoclinic channels $\{\Gamma_0^i\}_{i=1}^m\subset W^s(\Lambda_0)\pitchfork W^u(\Lambda_0)$ for some $m$ large enough.  
  
  Then, there exists $\gamma(F_0)>0$ and   $\varepsilon_0(F_0)>0$  such that for any  $f$ belonging to an open and dense subset of $\mathcal B^\omega(M\times N)$, the corresponding family $F_{\varepsilon,f}\in\mathrm{Ham}^\omega(M\times N)$ in \eqref{eq:deformationmap} satisfies the following for any $0<\varepsilon\leq \varepsilon_0$:
    \begin{enumerate}
        \item \label{it:Item1genericity} $F_{\varepsilon,f}$ displays a normally hyperbolic invariant manifold $\Lambda_\varepsilon$ for which assumptions \ref{it:H1}-\ref{it:H2} are verified. The corresponding inner map  $T_\varepsilon:N\to N$ in \eqref{eq:innerdynamicsdefn} and scattering maps (associated to the continuation of $\{\Gamma_0^i\}_{i=1}^m$)  $\{S_{i,\varepsilon}\}_{i=1}^m:N\to N$ in \eqref{eq:scattmapdynamicsdefn}, preserve a common smooth measure on $N$.
       \item  \label{it:Item2genericity}  There exists a subset $\mathcal A\subset N$ diffeomorphic to the  annulus $\mathbb A=\mathbb T\times[-1,1]$  on which the maps $\{T_\varepsilon, S_{1,\varepsilon}\}$ are of the form \eqref{eq:innermapsymbolic} and \eqref{eq:scattmapsymbolic} and verify hypothesis \ref{it:itemB1}-\ref{it:itemB5} with $\gamma=\gamma(F_0)$.
        \item \label{it:Item3genericity} The scattering maps $\{S_{i,\varepsilon}\}_{i=1}^m$ satisfy assumption \ref{it:B6} with $n=1$ and some $\tilde r\geq 1$.
    \end{enumerate}
\end{prop}

The proof of Proposition \ref{prop:genericityofassumptions} is given in Section 8 of \cite{paradela2026robustlytransitivebehaviorsymplectic}. Here we just give a heuristic idea of why these assumptions are abundant around maps $F_0$ as that in Theorem \ref{thm:mainham}.

\begin{rem}
     The very same conclusion in Proposition \ref{prop:genericityofassumptions} holds if $F_0\in \mathrm{Symp}^\omega(M\times N)$ is a direct product of symplectomorphisms satisfying assumptions \ref{it:Hyp1}-\ref{it:Hyp2} and we consider an unfolding $F_{\varepsilon,\theta}=F_0\circ\Phi^{\varepsilon}_\theta\in \mathrm{Symp}^\omega(M\times N)$
   with $\Phi_\theta^{\varepsilon}$ as in \eqref{eq:unfoldingforsymplectic}. See Section 8 in \cite{paradela2026robustlytransitivebehaviorsymplectic} for more details.
\end{rem}
\medskip

First, Item \ref{it:Item1genericity} is straightforward since by assumption, for $F_0$ there exists a normally hyperbolic cylinder $\Lambda_0=\{\gamma_0\}\times N$ which is already $\kappa$-normally hyperbolic, $\kappa$-pinched and  exhibits one full  homoclinic channel. It is straightforward to observe that, in this setting this implies that there actually exists a countable collection of disjoint homoclinic channels. In particular, for any $m\in\mathbb N$ we can select a family of full and disjoint homoclinic channels $\{\Gamma_0^i\}_{i}\subset M\times N$. The corresponding results for $\varepsilon>0$ follow by regular perturbation theory (see \cite{MR501173}, \cite{DelaLLaveScattmap} or  Section 8 in \cite{paradela2026robustlytransitivebehaviorsymplectic}). Moreover, thanks to the definition of $\mathcal B^\omega$ (see \eqref{eq:Ballperturbations}) we can assume that the maps $T_\varepsilon,\{S_{i,\varepsilon}\}_{i=1}^m\subset \mathrm{Diff}^{r-1}(N)$ and that, moreover,  $T_\varepsilon$ is $C^{r-1}$ close to the map $\phi_N\in \mathrm{Ham}^\omega(N)$.

Concerning Item \ref{it:Item2genericity}, in view of assumption \ref{it:Hyp1} and the fact that $T_\varepsilon$ is $C^{r-1}$ close to $\phi_N\in \mathrm{Ham}^\omega(N)$, a standard argument shows that close to the non-degenerate elliptic equilibrium,  there exists an invariant circle which, for $\varepsilon>0$ small enough, admits a continuation as a KAM curve $\mathcal T_\varepsilon\subset\Lambda_\varepsilon$ (see \cite{PoschelKAM} or Section  8 in \cite{paradela2026robustlytransitivebehaviorsymplectic}). Around this curve the inner dynamics  can be put in normal form \eqref{eq:innermapsymbolic} provided $r$ is large enough. This gives us the assumptions \ref{it:itemB1}-\ref{it:itemB2} and \ref{it:itemB5}. On the other hand, for the assumptions involving the scattering map associated to the homoclinic channel $\Gamma_\varepsilon^1\subset M\times N$, one can make use of the perturbative tools developed in \cite{DelaLLaveScattmap} to show that the transversality-torsion assumption in \ref{it:itemB3} holds for perturbations $f$ from an open and dense subset of $\mathcal B^\omega(M\times N)$.

Finally, Item \ref{it:Item3genericity}, i.e. genericity of assumption \ref{it:B6}, can be easily deduced by means of the asymptotic expansions in \cite{DelaLLaveScattmap} and a suitable parametric transversality argument.

\begin{rem}
    To be precise, in \cite{paradela2026robustlytransitivebehaviorsymplectic} we showed typicality of assumption \ref{it:B6} with $n=0$. However, the reader can check that the very same proof gives typicality for any fixed $n\in \mathbb N$.
\end{rem}

\subsection{Proof of Theorem \ref{thm:mainham}}
Proposition \ref{prop:genericityofassumptions} shows that for $F_{0}$ as in Theorem \ref{thm:mainham}, assumptions \ref{it:H1}-\ref{it:H2} are verified and the corresponding inner map  $T_\varepsilon:N\to N$ and scattering maps $\{S_{i,\varepsilon}\}_{i=1}^m:N\to N$ satisfy \ref{it:itemB1}-\ref{it:B6}. In this setting, Proposition  \ref{prop:mainnhim} and Proposition \ref{prop:Reduction Skew-prod to IFS} imply that there exists a set $\mathcal X_\varepsilon\subset M\times N$, an open neighbourhood $\mathcal Q_\varepsilon\supset \mathcal X_\varepsilon$   and a homeomorphism $\Psi_\varepsilon:\Sigma\times N\to \mathcal X_\varepsilon$ such that:
 \begin{enumerate}
        \item $\mathcal X_\varepsilon$ is locally maximal and normally hyperbolic for the return map $\mathcal R_\varepsilon$ induced on $\mathcal Q$,
        \item Any natural measure $\mathfrak m_\varepsilon$ on $\mathcal X_\varepsilon$ is invariant for $\mathcal R_\varepsilon|_{\mathcal X_\varepsilon}:\mathcal X_\varepsilon\to\mathcal X_\varepsilon $ and,
        \item 
        Any natural measure  $\mathfrak m_\varepsilon$ on $\mathcal X_\varepsilon$ is  ergodic for Lebesgue almost every $\varepsilon\in (-\varepsilon_0,\varepsilon_0)$.
    \end{enumerate}
The proof of Theorem \ref{thm:mainham} is complete.

\medskip

\section{Outline for the proof of Theorem \ref{thm:main}}\label{sec:outline}

The purpose of this section is to materialize the approach introduced in Section \ref{sec:IFSintro} (recall the discussion after Theorem \ref{thm:main}) to  outline a proof of Theorem \ref{thm:main} in the particular case where the maps are affine.

\begin{thm}\label{thm:mainaffine}
   For every $\lambda\in(1/2,1)$ denote by 
   \begin{equation}\label{eq:affine group}
   \mathcal G_\lambda^{\mathrm{aff}}=\langle F_{\mathrm{aff},1,\lambda},F_{\mathrm{aff},2,\lambda}\rangle\subset \mathrm{Diff}^\infty(\mathbb R^2).
   \end{equation}
   the group of transformations generated by  the pair of affine maps in \eqref{eq:affinemaps}. Then,  for $\mathrm{Leb}$-a.e. $\lambda\in (1/2,1)$ the group action $(\mathcal G_\lambda^{\mathrm{aff}},\mathbb R^2,\mathrm{Leb})$ is locally ergodic on $D=[0,1]^2$. That is, for every Borel measurable set $E$ which is $\mathcal G_\lambda^{\mathrm{aff}}$-invariant we have 
    \[
    \mathrm{Leb}_{D}(E\cap D)\in\{0,1\}.
    \]
\end{thm}
Below we give a complete proof of Theorem \ref{thm:mainaffine}. To do so we show that for a typical $\lambda\in (1/2,1)$, the pushforward of the symmetric Bernoulli measure $\mu$  on $\Omega$  under the coding map \eqref{eq:affinetwosidedcodingmap} is equivalent to Lebesgue measure on $D$.  With a view towards the generalization to the nonlinear case, i.e. Theorem \ref{thm:main}, special emphasis will be given to the geometric aspects of the construction.  
\medskip

    Note that although our goal is simply to ensure the existence of a symbolic measure satisfying \eqref{it:1}-\eqref{it:2}, we have made a rather particular choice: $\mu$ being the symmetric Bernoulli measure. This choice is well   motivated by the fact that, as we will see below, its projection under the coding map \eqref{eq:affinetwosidedcodingmap} is well understood after the work of Solomyak \cite{MR1356783} on Bernoulli convolutions.

\subsection{Geometric preliminaries}
We recall that
\[
F_{\mathrm{aff},i,\lambda}(x,y)
=(f_{x,i,\lambda}(x),f_{y,i,\lambda}(y))
\]
with
\begin{equation}\label{eq:symmetry affine}
f_{x,1,\lambda}(s)=\lambda s=f_{y,1,\lambda}^{-1}(s),
\qquad
f_{x,2,\lambda}(s)=1-\lambda(1-s)=f_{y,2,\lambda}^{-1}(s).
\end{equation} 
For convenience, we also introduce  the one-dimensional intervals
\begin{equation}\label{eq: I_i intervals}
I_{i,\lambda}=f_{x,i,\lambda}([0,1])\subset[0,1].
\end{equation}

\subsection{The pushforward of  Bernoulli measure under the coding map}
Consider the pushforward of $\mu$ (the symmetric Bernoulli measure on $\Omega=\{1,2\}^\mathbb Z$) under the coding map \eqref{eq:affinetwosidedcodingmap} particularized for the affine system
\[
\nu_\lambda^{\mathrm{aff}}=(\Pi_\lambda^{\mathrm{aff}})_*\mu.
\]
The first observation is that in the affine case the horizontal and vertical dynamics are completely decoupled. Hence, the coding map can be written as
\[
\Pi_\lambda^{\mathrm{aff}}(\omega^-;\omega^+)
=\big(\Pi_\lambda^{\mathrm{aff},-}(\omega^-),
\Pi_\lambda^{\mathrm{aff},+}(\omega^+)\big),
\]
where
\[
\Pi_\lambda^{\mathrm{aff},-}(\omega^-)
=\lim_{m\to\infty}
f_{x,\omega_{-1},\lambda}\circ\cdots\circ
f_{x,\omega_{-m},\lambda}(0)
\]
and
\[
\Pi_\lambda^{\mathrm{aff},+}(\omega^+)
=\lim_{n\to\infty}
f_{y,\omega_0,\lambda}^{-1}\circ\cdots\circ
f_{y,\omega_{n-1},\lambda}^{-1}(0).
\]
Therefore, if we write (here $\mu^\pm$ is the pushforward of $\mu$ under the projection $\pi^\pm:\omega\mapsto \omega^\pm$)
\[
\mu=\mu^-\otimes \mu^+,
\]
the symmetry \eqref{eq:symmetry affine} implies that the respective  pushforwards of $\mu^-$ and $\mu^+$ under these maps coincide. Denoting their common value by $\eta_\lambda$ we obtain
\begin{equation}\label{eq:affine Bernoulli pushforward}
\nu_\lambda^{\mathrm{aff}}=\eta_\lambda\otimes\eta_\lambda.
\end{equation}
Thus, in order to establish property \eqref{it:1}, it is enough to prove that $\eta_\lambda$ is equivalent to Lebesgue measure on $[0,1]$.

\subsection{Typical absolute continuity and equivalence to Lebesgue} \label{sec:typical AC affine}

 Unfortunately, it is well known since the work of Erdös \cite{MR311} that the measure $\eta_\lambda$  is  not always absolutely continuous. However, one can show that typically in the parameter $\lambda$, these measures are equivalent to Lebesgue.

\begin{prop}[\cite{MR1356783} and \cite{MR1458276}]\label{prop:ACaffine}
For $\mathrm{Leb}$-a.e. $\lambda\in (1/2,1)$, the measure $\eta_\lambda$  has  a density $a_{\mathrm{aff},\lambda}\in L^2([0,1])$ and is strictly positive almost everywhere on $[0,1]$.
\end{prop}
\begin{proof}

Proposition \ref{prop:ACaffine} follows from well-known results in the theory of Bernoulli convolutions. 
On one hand, Solomyak showed in \cite{MR1356783} that the measure $\eta_\lambda\in \mathcal M([0,1])$ has a $L^2$ density for Lebesgue a.e. $\lambda\in(1/2,1)$ so the first part of  Proposition \ref{prop:ACaffine} holds. On the other hand, Mauldin and Simon showed in \cite{MR1458276} that either 
\[
\mathrm{Leb}_{[0,1]}\ll \eta_\lambda\qquad\qquad\text{or}\qquad\qquad \mathrm{Leb}_{[0,1]}\perp \eta_\lambda.
\]
Thus, we obtain the second part of Proposition \ref{prop:ACaffine}.
\end{proof}
We then  deduce the following making use of the expression \eqref{eq:affine Bernoulli pushforward}.
\begin{prop}\label{prop: properties of the affine measure}
    For $\mathrm{Leb}$-a.e. $\lambda\in (1/2,1)$ the measure symmetric Bernoulli measure $\mu$ on $\Omega$ satisfies properties \eqref{it:1}-\eqref{it:2}.
\end{prop}
From Proposition \ref{prop: properties of the affine measure} it is straightforward do complete the proof of Theorem \ref{thm:mainaffine}. This is done in Section \ref{sec:Hopfaffine}. Before continuing, having in mind the nonlinear problem in Theorem \ref{thm:main}, it is worth to analyze the proof of Proposition \ref{prop:ACaffine} from a more geometric viewpoint.
\medskip

\subsection*{The transversality argument}\label{sec:transversalityaffine}
With a view towards analyzing nonlinear one-dimensional IFS, a different proof of the main result in \cite{MR1356783} was given in \cite{MR1386842}. Roughly speaking, the idea goes as follows (the interested reader is referred to \cite{MR1386842} for the details).
\medskip

Given two different sequences $\omega,\eta\in \{1,2\}^\mathbb N$ we denote by 
\[
\omega\wedge \eta=\min\{i\in\mathbb N\cup\{0\}\colon \omega_i\neq \eta_i\}.
\]
A standard application of the Vitali covering theorem shows that the measure $\eta_\lambda\in \mathcal M([0,1])$  is absolutely continuous if and only if
\[
D(\eta_\lambda,x):=
\liminf_{r\to0}\frac{1}{2r}\eta_\lambda(B_r(x))
<\infty
\qquad\text{for }\eta_\lambda\text{-a.e. }x.
\]
If we moreover show that 
\begin{equation}\label{eq:finiteS}
\mathcal S=\int_{J}\int_{[0,1]} D(\eta_\lambda,x)\ \mathrm d\eta_\lambda(x) \mathrm d\lambda<\infty,
\end{equation}
then $\eta_\lambda$ has an $L^2$ density for a.e. $\lambda\in J\subset  (1/2,1)$. To establish \eqref{eq:finiteS} we use Fatou's lemma and the definition of $\eta_\lambda$ to show that 
\[
\mathcal S\leq \liminf_{r\to 0}\frac{1}{2r}\int_{\{1,2\}^\mathbb N}\int_{\{1,2\}^\mathbb N} \underbrace{\mathrm{Leb}\{\lambda\colon |\Pi_\lambda^{\mathrm{aff},-}(\tau)-\Pi_\lambda^{\mathrm{aff},-}(\omega)|\leq r\}}_{\star} \ d\mu^-(\omega)d\mu^-(\tau).
\]
In order to estimate $(\star)$, one first observes that for $\omega\wedge\tau=k$,
\[
\Pi_\lambda^{\mathrm{aff},-}(\tau)-\Pi_\lambda^{\mathrm{aff},-}(\omega)
=\pm(1-\lambda)\lambda^k g(\lambda),
\]
where $g$ is a power series of the form
\[
g(\lambda)=1+\sum_{n\geq1}b_n\lambda^n,
\qquad b_n\in\{-1,0,1\}.
\]
Suppose now that for the interval $J=[\lambda_0,\lambda_1]\subset(1/2,1)$ there exists $\delta>0$ such that for every $g$ as above
\begin{equation}\label{eq:Transversality}
|g(\lambda)|\leq\delta
\qquad\Longrightarrow\qquad
|g'(\lambda)|\geq\delta.
\end{equation}
Under Property \eqref{eq:Transversality} 
\[
\mathrm{Leb}\{\lambda\colon |\Pi^{\mathrm{aff},-}_\lambda(\tau)-\Pi^{\mathrm{aff},-}_\lambda(\omega)|\leq r\}\leq \delta^{-1} \lambda^{-k} r\leq \delta^{-1} \lambda_0^{-k} r
\]
and we deduce that 
\[
\mathcal S\leq \frac{1}{2\delta}\sum_{k\in \mathbb N} (2\lambda_0)^{-k}<\infty.
\]
To conclude the proof, the authors in \cite{MR1386842} first show that Property \eqref{eq:Transversality} holds on any interval of the form $J=[\lambda_0,2^{-2/3}]$ for any $\lambda_0>1/2$. This proves that $\nu_\lambda$ has a $L^2$ density for a.e. $\lambda \in J$. The extension of this conclusion to a.e. $\lambda\in(1/2,1)$ uses features inherent to the affine setting so we will not delve on the details and refer the interested reader to \cite{MR1386842}.

\subsection{Proof of Theorem \ref{thm:mainaffine}}
\label{sec:Hopfaffine}

We fix $\lambda\in(1/2,1)$ so that Proposition \ref{prop: properties of the affine measure} holds and omit the dependence on $\lambda$ from the notation. Let $E$ be a measurable subset invariant under the group action $\mathcal G^{\mathrm{aff}}$ in \eqref{eq:affine group} and define
\[
\widehat E=(\Pi^{\mathrm{aff}})^{-1}(E\cap D)\subset\Omega.
\]
Since
\[
\Pi^{\mathrm{aff}}(\sigma\omega)
=F_{\mathrm{aff},\omega_0}\Pi^{\mathrm{aff}}(\omega),
\]
we must have
\[
\omega\in\widehat E
\quad\Longrightarrow\quad
\sigma\omega\in\widehat E,
\]
so $\widehat E$ is $\sigma$-invariant. Ergodicity of $(\Omega,\sigma,\mu)$ then implies that
\begin{equation}\label{eq:0-1 law for E hat}
\mu(\widehat E)\in\{0,1\}.
\end{equation}
But then, since by Proposition \ref{prop: properties of the affine measure} we have $\Pi^{\mathrm{aff}}_*\mu\sim\mathrm{Leb}_D$, \eqref{eq:0-1 law for E hat} implies that
\[
\mathrm{Leb}_D(E)\in\{0,1\}.
\]
The proof of Theorem \ref{thm:mainaffine} is complete.

\section{Notation and organization of the remaining sections}
Here we describe the notation that will be used in the remaining part of the paper.

\subsection{Functional spaces}
Given two Banach  manifolds $X,Y$:

\subsection*{$C^k$ spaces}  Given $k\in\mathbb N$ we let 
$C^k(X,Y)$
be the space of $k$-times differentiable functions $f:X\to Y$ whose derivatives of order  $k$ are continuous. In the case where:
\begin{itemize}
    \item $Y=\mathbb R^d$ for $d\geq 2$ we denote by 
\[
\lVert f\rVert_{C^k(X)}=\sum_{0\leq |\alpha|\leq k} \lVert \partial^\alpha f\rVert_{C^0(X)} 
\]
which makes $C^k(X,\mathbb R^d)$ a Banach space.
\item $Y=\mathbb R$ we simply write 
\[
|f|_{C^k(X)}=\sum_{0\leq |\alpha|\leq k} |\partial^\alpha f|_{C^0(X)}
\]
which makes $C^k(X,\mathbb R)$ a Banach space. Moreover we will also use make use of the semi-norm 
\[
[f]_{C^k(X)}=\sum_{|\alpha|=k}|\partial^\alpha f|_{C^0(X)}
\]
\end{itemize}

\subsection{Word and sequence spaces}
Given some finite alphabet $A$ we let 
\[
\mathcal W_A= \bigcup_{n\geq 1} A^n
\]
be the space of finite words and denote by 
\[
\Omega_A= A^\mathbb Z
\]
be the space of bi-infinite sequences. We use the notation $\omega=(\omega_0,\dots,\omega_n)$ do denote an element of $\mathcal W_A$ (i.e. a finite word) and $\omega=(\dots,\omega_{-1};\omega_0,\dots)$ to denote an element of $\Omega_A$ (i.e. a sequence). We also make use of the notation 
\[
[\omega_{-m},\dots,\omega_{-1};\omega_0,\dots,\omega_n]=\{\tilde\omega\in \Omega_A\colon \tilde\omega_k=\omega_k,\ -m\leq k\leq n\}\subset \Omega_A
\]
for cylinders on $\Omega_A$. Sometimes we will just write $[\omega]$ when the specified symbols do not matter or are clear from the context. In addition, we will consider half-infinite sequences 
\[
\omega^-=(\dots,\omega_{-n},\dots,\omega_{-1})\in \Omega_A^-=A^\mathbb N\qquad\qquad\omega^+=(\omega_{0},\dots,\omega_{n},\dots)\in \Omega_A^+=A^\mathbb N
\]
Also, given $n\in\mathbb N$ and $\omega\in A^\mathbb N$ we will write 
\[
\omega|n=(\omega_0,\dots,\omega_{n-1})\in \mathcal W_A.
\]
\begin{rem}
    In most cases it will be clear from the context what alphabet  $A$ we are referring to so we will not include it in the notation.
\end{rem}

\subsection{Organization of the remaining sections}

In Section \ref{sec:geometric preli} we discuss relevant geometric features exhibited by $\epsilon$-good maps which are of importance to understand the coding map in \eqref{eq:affinetwosidedcodingmap}. Then, in Sections \ref{sec:typicalAC} and \ref{sec:equivalenceLebesgue} we show that the pushforward of the symmetric Bernoulli measure under the coding map is equivalent to Lebesgue, completing the proof of Theorem \ref{thm:main}. Finally, Sections \ref{sec:proof TT } and \ref{sec: from IFS to skewprod} are devoted to the proofs of Propositions \ref{prop:TT and HB} and \ref{prop:Reduction Skew-prod to IFS}.

\section{Geometric preliminaries}\label{sec:geometric preli}

In this section we present some geometric features exhibited by $\epsilon$-good maps. The main idea that we exploit is that $\epsilon$-good maps preserve  \textit{common}  (for all the maps in the collection)  stable and unstable cones. This will allow us to define certain spaces of un/stable curves which are invariant under the action of any element in the corresponding iterated function system. From this observation, we will ultimately be able to establish a symbolic coding which to any sequence $\omega^{\pm}\in \{1,2\}^\mathbb N$ assigns a un/stable curve.

\begin{rem}
    Most of  the results in this section follow from well known techniques so they are presented without proofs.
\end{rem}

\begin{rem}
    Throughout this section $C$ denotes a positive constant which is independent of $\epsilon$, $\lambda$ and (if applicable) of the sequence $\omega^{\pm}\in \{1,2\}^\mathbb N$ or word $\omega\in \mathcal W$. Moreover, the value of $C$ might change from line to line.
\end{rem}

\subsection{Uniformization of the geometry}
 Let $\{F_{1,\lambda},F_{2,\lambda}\}_{\lambda\in J}\in \mathrm{Diff}^3(N)$ be a family of $\epsilon$-good maps and let $\phi_\lambda:\mathcal D\to N$ and $\widetilde  F_{i,\lambda}:\mathcal D\to \mathbb R^2$, $i=1,2$ be as in Definition \ref{defn:epsilon good}. Recall also the definition of $D_\lambda\subset\mathcal D$ in Definition \ref{defn:hopfbox}. 
 \medskip

 We will find convenient to further introduce a $C^3$ coordinate transformation $\tilde\phi_\lambda:  D=[0,1]^2\to D_F$ 
 depending $C^3$ on $\lambda$ and satisfying (uniformly in $\lambda$ and $z\in D$),
\[
\lVert \tilde\phi_\lambda-\mathrm{id}\rVert _{C^3(D\times J)}\leq C\epsilon.
\]
We then define the maps 
\begin{equation}\label{eq:straightenedmaps}
 \mathtt F_{i,\lambda}:=\tilde\phi_\lambda^{-1}\circ \widetilde F_{i,\lambda}\circ\tilde\phi_\lambda :D\to \mathbb R^2
\end{equation}
and observe that the family $\{\mathtt F_{1,\lambda},\mathtt F_{2,\lambda}\}_\lambda:D\to \mathbb R^2$ satisfies 
\begin{equation}\label{eq:mathtt is epsilon good}
\max_{i=1,2}\  \lVert \mathtt F_{i,\lambda}-F_{\mathrm{aff},i,\lambda} \rVert_{C^3(D\times J)}\leq C\epsilon.
\end{equation}
We denote by 
\[
V_{i,\lambda}=\mathtt F_{i,\lambda}(D)\cap D\qquad\qquad H_{i,\lambda}=\mathtt F_{i,\lambda}^{-1}(D)\cap D
\]
so that $\mathtt F_{i,\lambda}:H_{i,\lambda}\to V_{i,\lambda}$  and $\mathtt F_{i,\lambda}^{-1}:V_{i,\lambda}\to H_{i,\lambda}$.
\medskip

\subsection{The graph transform operator}

Define the family of \textit{stable curves}
\begin{equation}\label{eq:setofcurvesvertical}
\Gamma^s=\{ (x,\gamma(x))\subset D\colon \gamma\in C^2,\ |\gamma'(x)|_{C^0}\leq 1,\ |\gamma''(x)|_{C^0}\leq 1,\ x\in[0,1] \}
\end{equation}
and \textit{unstable curves}
\begin{equation}\label{eq:setofcurveshorizontal}
\Gamma^u=\{ (\gamma(y),y)\subset D\colon \gamma\in C^2,\ |\gamma'(y)|_{C^0}\leq 1,\ |\gamma''(y)|_{C^0}\leq 1,\ y\in[0,1] \}.
\end{equation}
\begin{rem}
    In the following, we will repeteadly abuse notation and denote by $\gamma$ both an element of $\Gamma^{s,u}$ and the function appearing in its parametrization. In this sense, we can understand 
    \[
    \Gamma^{u,s}\subset C^2([0,1]).
    \]
\end{rem}
\medskip

Now, given $\gamma\in \Gamma^s$, for $i=1,2$ we denote by $x_{i,\lambda}:[0,1]\to[0,1]$ the unique function such that 
\[
x=\mathtt F_{i,\lambda,x}^{-1}(x_{i,\lambda}(x),\gamma(x_{i,\lambda}(x)))
\]
and, 
given $\gamma\in \Gamma^u$, for $i=1,2$ we denote by $y_{i,\lambda}:[0,1]\to[0,1]$ the unique function such that 
\begin{equation}\label{eq:implicit function theorem y coordinate}
y=\mathtt F_{i,\lambda,y}(\gamma(y_{i,\lambda}(y)),y_{i,\lambda}(y)).
\end{equation}
Then, we introduce the \textit{stable graph transform operator} associated to the map $\mathtt F_{i,\lambda}$
\[
\mathcal G_{i,\lambda}^{s}(\gamma)=\{(x, \gamma_{i,\lambda}(x))\subset D\colon \gamma_{i,\lambda}(x)=F_{i,\lambda,y}^{-1}(x_{i,\lambda}(x),\gamma(x_{i,\lambda}(x)),\ x\in [0,1]\}
\]
 and the \textit{unstable graph transform operator} associated to the map $\mathtt F_{i,\lambda}$
\[
\mathcal G_{i,\lambda}^{u}(\gamma)=\{(\gamma_{i,\lambda}(y),y)\subset D\colon \gamma_{i,\lambda}(y)=F_{i,\lambda,x}(\gamma(y_{i,\lambda}(y)),y_{i,\lambda}(y)),\  y\in [0,1]\}.
\]
\begin{lem}\label{lem:graphtransform}
   Fix any closed interval $J\subset (1/2,1)$ and let $\epsilon$ be small enough. Then,  for any $\lambda \in J$ and any  $\gamma \in \Gamma^{u,s}$ the following estimates hold
    \begin{equation}\label{eq: C1 estimate graph transform}
        [\mathcal G_{i,\lambda}^{u,s}(\gamma)]_{C^1}\leq (\lambda^2+C\epsilon)[\gamma]_{C^1}+C\epsilon
    \end{equation}
    and
     \begin{equation}\label{eq: C2 estimate graph transform}
        [\mathcal G_{i,\lambda}^{u,s}(\gamma)]_{C^2}\leq (\lambda^3+C\epsilon)[\gamma]_{C^1}+C\epsilon
    \end{equation}
    In addition, for any $\gamma,\gamma' \in \Gamma^{u,s}$ 
    \begin{equation}\label{eq:curves get closer}
        |\mathcal G^{u,s}_{i,\lambda}(\gamma)-\mathcal G_{i,\lambda}^{u,s}(\gamma')|_{C^1}\leq (\lambda+C\epsilon)|\gamma-\gamma'|_{C^1}.
    \end{equation}
\end{lem}
\medskip

We note that, in particular, Lemma \ref{lem:graphtransform} shows that 
\[
\mathcal G^{u,s}_{i,\lambda}(\Gamma^{u,s})\subset \Gamma^{u,s}
\]
and it is a contraction with respect to the $C^1$ norm.
\medskip

\subsubsection*{Parameter dependence:} The following estimate will be relevant to run the transversality argument in Section \ref{sec:typicalAC} (recall the discussion in Section \ref{sec:outline}). 

\begin{lem}
    Fix any closed interval $J\subset (1/2,1)$ and let $\epsilon$ be small enough. Then,  for any $\lambda \in \mathrm{int}(J)$ and  any $\gamma \in \Gamma^{u,s}$
    \begin{equation}\label{eq: parameter dependence graph transform}
    |\partial_\lambda\mathcal G^{u,s}_{i,\lambda}(\gamma)-\partial_\lambda\mathcal G^{u,s}_{\mathrm{aff},i,\lambda}(\gamma)|_{C^1}\leq C\epsilon.
    \end{equation}
    Moreover, for any $\gamma,\gamma' \in \Gamma^s$ \begin{equation}\label{eq:curves get closer 2}
        |\partial_\lambda(\mathcal G^{u,s}_{i,\lambda}(\gamma)-\mathcal G^{u,s}_{i,\lambda}(\gamma'))|_{C^0}\leq C|\gamma-\gamma'|_{C^1}.
    \end{equation}
\end{lem}

We will also make use of the following result.

\begin{lem}
   Fix any closed interval $J\subset (1/2,1)$ and let $\epsilon$ be small enough. Then,  for any $\lambda \in J$  the operator $\mathcal G^{u,s}_{i,\lambda}:\Gamma^{u,s}\to C^1$  is Fréchet differentiable and, at any $\gamma\in \Gamma^{u,s}$ its differential satisfies
    \begin{equation}\label{eq:frechet differential}
    \begin{split}
    |D\mathcal G_{i,\lambda}^{u,s}(\gamma)[v]|_{C^0}\leq &(\lambda+C\epsilon) |v|_{C^0}\\
    [D\mathcal G_{i,\lambda}^{u,s}(\gamma)[v]]_{C^1}\leq & C\epsilon|v|_{C^0}+(\lambda^2+C\epsilon) [v]_{C^1}
    \end{split}
    \end{equation}
    and
    \begin{equation}\label{eq:frechet differential 2}
     |D\mathcal G_{i,\lambda}^{u,s}(\gamma)[v]-D\mathcal G^{u,s}_{\mathrm{aff},i,\lambda}(\gamma)[v]|_{C^0}\leq C\epsilon |v|_{C^1}.
    \end{equation}
    Moreover, for any $\gamma,\gamma'\in \Gamma^{u,s}$
\begin{equation}\label{eq:frechet differential 3}
     |D\mathcal G_{i,\lambda}^{u,s}(\gamma)[v]-D\mathcal G_{i,\lambda}^{u,s}(\gamma')[v]|_{C^0}\leq C\epsilon|\gamma-\gamma'|_{C^1} |v|_{C^1}.
    \end{equation}
\end{lem}
\begin{proof}
 We omit the dependence on $\lambda$ as it plays no role in our discussion. We only analyze the unstable case, the stable one being analogous.  It follows from the estimates \eqref{eq:mathtt is epsilon good} and the definition of $\mathcal G^{u}_i$ and the implicit function theorem that the map $\mathcal G^{u}_i:\Gamma^{u}\to C^1$ is Fréchet differentiable. A straightforward computation shows that 
    \begin{equation}\label{eq: frechet derivative formula}
D\mathcal G_i^u(\gamma)[v](y)
=
\underbrace{\frac{\det \mathrm  D\mathtt F_i(\gamma(y_i),y_i)}
{\partial_y \mathtt F_{i,y}(\gamma(y_i),y_i)
+\partial_x \mathtt F_{i,y}(\gamma(y_i),y_i)\gamma'(y_i)}}_{A_i(y)}
\ v(y_i)
\end{equation}
where $y_i:[0,1]\to [0,1]$ is definied implicitly by \eqref{eq:implicit function theorem y coordinate}. After differentiation
\begin{equation}\label{eq: frechet derivative formula derivative}
\left(D\mathcal G_i^u(\gamma)[v]\right)'(y)=A_i'(y)v(y_i)+A_i(y)v'(y_i(y))y_i'(y)
\end{equation}
Note that for the affine case the above expression reduces to  
\[
D\mathcal G_{\mathrm{aff},i}^u(\gamma)[v](y)
=\lambda\ v(y_{i,\mathrm{aff}}).
\]
For the non-linear case the estimates \eqref{eq:frechet differential} follow from \eqref{eq: frechet derivative formula} after  a routine computation which exploits that 
\[
|\det \mathrm D\mathtt F_i-1|_{C^1}\leq C\epsilon,\qquad\qquad |\partial_y \mathtt F_{i,y}-\lambda^{-1}|_{C^1}\leq C\epsilon\qquad\qquad |\partial_x \mathtt F_{i,y}|_{C^1}\leq C\epsilon.\qedhere
\]
\end{proof}

\medskip

\subsection{Un/stable manifolds and the coding map}

Having studied the action of the maps $\mathtt F_{i,\lambda}$  on un/stable curves, we now introduce what we call stable and unstable rectangles. Namely, a subset $H\subset D$ is a \textit{stable rectangle} if there exist $\gamma_d,\gamma_u\in \Gamma^s$ such that 
\[
H=\{(x,y)\in D\colon \gamma_d(x)\leq y\leq \gamma_u(x),\ x\in[0,1]\}.
\]
We define the width of a stable rectangle $H$ as 
\[
\ell_H=\min_{x\in [0,1]}\underbrace{\gamma_u(x)-\gamma_d(x)}_{\ell_H(x)}>0.
\]
Analogously, a subset $V\subset D$ is a \textit{unstable rectangle} if there exist $\gamma_l,\gamma_r\in \Gamma^u$ such that 
\[
V=\{(x,y)\in D\colon \gamma_l(y)\leq x\leq \gamma_r(y),\ y\in[0,1]\}.
\]
We define the width of a unstable rectangle $V$ as 
\[
\ell_V=\min_{y\in [0,1]}\underbrace{\gamma_r(y)-\gamma_l(y)}_{\ell_V(y)}>0.
\]
Making use of the third item in Lemma \ref{lem:graphtransform}, it is not difficult to check that (the intersection of $D$ with) the image of any unstable (resp. stable) rectangle under the map $\mathtt F_{i,\lambda}$ (resp. $\mathtt F_{i,\lambda}^{-1}$) is again an unstable (resp. stable) rectangle with strictly smaller width. Hence, we can iterate this construction. That is, given any $\tau=(\tau_0,\dots,\tau_{n-1})\in \mathcal W$ the sets
    \[
    V_{\tau,\lambda}=D\cap \mathtt F_{\tau,\lambda}^{n}(D)\qquad\qquad H_{\tau,\lambda}=D\cap \mathtt F_{\tau,\lambda}^{-n}(D)
    \]
are, respectively, vertical and horizontal rectangles. 
By construction $V_\tau$ (resp. $H_\tau$) consists of points  $z\in D$ for which the finite orbit segment 
\[
\{\mathtt F_{\tau_{n-k},\lambda}^{-1}
\circ\cdots\circ
\mathtt F_{\tau_{n-1},\lambda}^{-1}(z)
\}_{k=1}^{n}\subset D\qquad\qquad (\text{resp. }\{\mathtt F_{\tau_k,\lambda}\circ \cdots\circ \mathtt F_{\tau_{0},\lambda}(z)\}_{k=0}^{n-1}\subset D).
\]
The following technical lemma will prove useful later.
\begin{lem}\label{lem:widht of rectangles}
    Given any $\tau=(\tau_0,\dots, \tau_{n-1})\in \mathcal W$ 
    \[
  (\lambda-C\epsilon)^n\leq  \ell_{V_{\tau,\lambda}}(y)\leq (\lambda+C\epsilon)^n\qquad\qquad  (\lambda-C\epsilon)^n\leq  \ell_{H_{\tau,\lambda}}(x)\leq (\lambda+C\epsilon)^n\
    \]
    for all $y\in[0,1]$ and all $x\in [0,1]$.
\end{lem}

The proof of Lemma \ref{lem:widht of rectangles} follows from iteration of the estimate \eqref{eq:curves get closer} and it is left to the reader. We now replace  words by sequences $\omega\in\{1,2\}^{\mathbb N}$. 

\begin{lem}\label{lem: construction un/stable manifolds}
   For any $\omega^-=(\dots, \omega_{-n},\dots,\omega_{-1})\in\{1,2\}^{\mathbb N}$ and $\omega^+=(\omega_0,\dots,\omega_n,\dots)\in\{1,2\}^{\mathbb N}$ the sets
    \[
    \gamma^u_{\omega^-,\lambda}=\bigcap_{n\geq 1} V_{(\omega_{-n},\dots,\omega_{-1})}\qquad\qquad \gamma^s_{\omega^+,\lambda}=\bigcap_{n\geq 1}H_{(\omega_{0},\dots,\omega_{n-1})}
    \]
    are, respectively, unstable and stable curves. That is, $\gamma^u_{\omega,\lambda}\in \Gamma^u$ and $\gamma^s_{\omega,\lambda}\in \Gamma^s$. Moreover, 
    \begin{equation}\label{eq: unstable manifolds are almost vertical}
    [\gamma^u_{\omega,\lambda}]_{C^1},[\gamma^s_{\omega,\lambda}]_{C^1}\leq C\epsilon.
    \end{equation}
\end{lem}
The first part of the statement is an easy  consequence of  Lemma \ref{lem:widht of rectangles}. For the second one it is enough to iterate the inequality \eqref{eq: C1 estimate graph transform}.
\medskip

Lemma \ref{lem: construction un/stable manifolds} allows us to assign to each infinite sequence $\omega\in \{1,2\}^\mathbb N$ a curve $\gamma^{u,s}_{\omega,\lambda}$. For convenience, since, in contrast with the affine case, the dependence of unstable curves on $y$ (resp. of stable curves on $x$) is non-trivial we will find convenient to introduce the family of coding maps
\begin{equation}\label{eq:codingnonlinearproof}
\Pi_{y,\lambda}^-:\omega^-=(\dots,\omega_{-n},\dots,\omega_{-1})\in \{1,2\}^\mathbb N\to (\gamma^u_{\omega^-,\lambda}(y),y)\in D
\end{equation}
for $y\in [0,1]$ and 
\begin{equation}\label{eq:coding map x}
\Pi_{x,\lambda}^+:\omega^+=(\omega_0,\dots,\omega_n,\dots)\in \{1,2\}^\mathbb N\to (x,\gamma^s_{\omega^+,\lambda}(x))\in D
\end{equation}
for $x\in [0,1]$.
\subsection{Regularity of the coding with respect to $\lambda$}\label{sec:regularitycoding}
The following technical result will prove crucial to run the transversality argument in Section \ref{sec:typicalAC}.

\begin{lem}
Fix any closed interval $J\subset(1/2,1)$. Then, for any $\epsilon$ small enough, any $\lambda\in J$, any $y\in [0,1]$ and any $\omega\in \{1,2\}^\mathbb N$
\begin{equation}\label{eq:regularity coding}
\lVert \Pi_{y,\lambda} (\omega)-\Pi_{y,\lambda}^{\mathrm{aff}}(\omega)\rVert \leq C\epsilon,\qquad\qquad \lVert \partial_\lambda(\Pi_{y,\lambda} (\omega)-\Pi_{y,\lambda}^{\mathrm{aff}}(\omega))\rVert \leq C\epsilon.
\end{equation}
Moreover, 
\begin{equation}\label{eq: regularity vertical coding}
\lVert \partial_y(\Pi_{y,\lambda} (\omega)-\Pi_{y,\lambda}^{\mathrm{aff}}(\omega))\rVert \leq C\epsilon.
\end{equation}
Similar results hold for the coding map $\Pi_{x,\lambda}$.
\end{lem}

\begin{proof}
 We omit the dependence in $\lambda$ from the notation and write $\gamma_\omega$ instead of $\gamma_\omega^u$. In the following we will use repeatedly that for any $\omega=(\dots,\omega_{-n},\dots,\omega_{-1})\in\{1,2\}^\mathbb N$
 \begin{equation}\label{eq:recursion coding}
 \gamma_\omega=\mathcal G_{\omega_0}(\gamma_{\sigma\omega})
 \end{equation}
  as elements of $\Gamma^u$.
\medskip

\noindent\textit{Estimate for $\Pi_{y} (\omega)-\Pi_y^{\mathrm{aff}}(\omega)$:} 
Define 
\[
E_0=\sup_{\omega,y}\lVert \Pi_{y} (\omega)-\Pi_y^{\mathrm{aff}}(\omega)\rVert=\sup_{\omega}|\gamma_\omega-\gamma_\omega^{\mathrm{aff}}|_{C^0}
\]
We write
\[
\gamma_\omega-\gamma_\omega^{\mathrm{aff}}=\mathcal G_{\omega_0}(\gamma_{\sigma\omega})-\mathcal G_{\omega_0}^{\mathrm{aff}}(\gamma_{\sigma\omega}^{\mathrm{aff}})= \mathcal G_{\omega_0}(\gamma_{\sigma\omega})-\mathcal G_{\omega_0}(\gamma_{\sigma\omega}^{\mathrm{aff}})+\mathcal G_{\omega_0}(\gamma_{\sigma\omega}^{\mathrm{aff}})-\mathcal G_{\omega_0}^{\mathrm{aff}}(\gamma_{\sigma\omega}^{\mathrm{aff}})
\]
Then, from the estimate \eqref{eq:curves get closer} and the fact that for $i=1,2$ the operators $\mathcal G_{i},\mathcal G_{i}^{\mathrm{aff}}$ are $O(\epsilon)$-close we deduce that 
\[
|\gamma_\omega-\gamma_\omega^{\mathrm{aff}}|_{C^0}\leq (\lambda+C\epsilon)|\gamma_{\sigma\omega}-\gamma_{\sigma\omega}^{\mathrm{aff}}|_{C^0}+C\epsilon.
\]
Taking the supremum over $\omega\in\{1,2\}^\mathbb N$ we obtain
\[
E_0\leq (\lambda+C\epsilon)E_0+C\epsilon\qquad\Longrightarrow\qquad E_0\leq \frac{C\epsilon}{1-(\lambda+C\epsilon)}.
\]
\medskip

\noindent\textit{Estimate for $\partial_\lambda(\Pi_{y} (\omega)-\Pi_y^{\mathrm{aff}}(\omega))$:} Define 
\[
W_\omega(y)=\partial_\lambda \gamma_\omega(y)\qquad \qquad W_\omega^{\mathrm{aff}}(y)=\partial_\lambda \gamma_\omega^{\mathrm{aff}}(y)
\]
and let 
\[
E_1=\sup_{\omega,y}\lVert \partial_\lambda(\Pi_{y} (\omega)-\Pi_y^{\mathrm{aff}}(\omega))\rVert=\sup_{\omega}|W_\omega-W_\omega^{\mathrm{aff}}|_{C^0}.
\]
Differentiating \eqref{eq:recursion coding}
\[
V_{\omega,\lambda}=\partial_\lambda \mathcal G_{\omega_0}(\gamma_{\sigma\omega})+D\mathcal G_{\omega_0}(\gamma_{\sigma\omega})[W_{\sigma\omega}].
\]
where $D\mathcal G_{i}(\gamma)[\cdot]$ stands for the Fréchet differential of the operator $\mathcal G_i$ evaluated at the curve $\gamma\in \Gamma^u$. We write 
\[
W_\omega-W_\omega^{\mathrm{aff}}=\partial_\lambda \mathcal G_{\omega_0}(\gamma_{\sigma\omega})-\partial_\lambda \mathcal G_{\omega_0}^{\mathrm{aff}}(\gamma_{\sigma\omega}^{\mathrm{aff}})+D\mathcal G_{\omega_0}(\gamma_{\sigma\omega})[W_{\sigma\omega}]-D\mathcal G_{\omega_0}^{\mathrm{aff}}(\gamma_{\sigma\omega}^{\mathrm{aff}})[W_{\sigma\omega}^{\mathrm{aff}}]
\]
On one hand, using \eqref{eq:curves get closer 2}
\[
\sup_{\omega}|\partial_\lambda \mathcal G_{\omega_0}(\gamma_{\sigma\omega})-\partial_\lambda \mathcal G_{\omega_0}^{\mathrm{aff}}(\gamma_{\sigma\omega}^{\mathrm{aff}})|_{C^0}\leq C E_0+C\epsilon.
\]
On the other hand, making use of \eqref{eq:frechet differential}-\eqref{eq:frechet differential 3}
\[
\sup_\omega|D\mathcal G_{\omega_0}(\gamma_{\sigma\omega})[W_{\sigma\omega}]-D\mathcal G_{\omega_0}^{\mathrm{aff}}(\gamma_{\sigma\omega}^{\mathrm{aff}})[W_{\sigma\omega}^{\mathrm{aff}}]|_{C^0}\leq (\lambda+C\epsilon) E_1+ C E_0+C\epsilon.
\]
Combining both 
\[
E_1\leq (\lambda+C\epsilon)E_1+CE_0+C\epsilon\leq (\lambda+C\epsilon)E_1+C\epsilon\qquad\Longrightarrow \qquad E_1\leq \frac{C\epsilon}{1-(\lambda+C\epsilon)}
\]
\medskip

\noindent\textit{Estimate for $\partial_y\Pi_{y,\lambda} (\omega)$}: Define
\[
U_\omega(y)=\partial_y \gamma_\omega(y)
\]
and write 
\[
E_2=\sup_\omega|U_\omega|_{C^0}.
\]
From the recursion \eqref{eq:recursion coding} and the estimate \eqref{eq: C1 estimate graph transform} 
\[
[U_\omega]_{C^0}=[\gamma_\omega]_{C^1}\leq (\lambda^2+C\epsilon)[\gamma_{\sigma\omega}]_{C^1}+C\epsilon
\]
and, arguing as above, 
\[
E_2\leq C\epsilon.\qedhere
\]
\end{proof}

\section{Typical absolute continuity}\label{sec:typicalAC}
We now deploy the powerful parameter-transversality theory  developed in the framework of Bernoulli convolutions (of particular relevance for our analysis are \cite{MR1749437,MR4385129}) to establish typical absolute continuity of the measure 
\begin{equation}\label{eq:pushforward definitive}
\nu_{\lambda}=(\Pi_\lambda)_* \mu
\end{equation}
with $\mu$ being the symmetric Bernoulli measure on $\Omega=\{1,2\}^\mathbb Z$ and $\Pi_\lambda$ being the coding map \eqref{eq:affinetwosidedcodingmap}.

\begin{thm}\label{thm:absolutecontinuity} Fix any closed interval $J\subset (1/2,2^{-2/3})$. Then, for $\epsilon$ small enough and almost every $\lambda\in J$ the measure $\nu_\lambda\in \mathcal M(D)$ defined in \eqref{eq:pushforward definitive}  is absolutely continuous with respect to Lebesgue. 
\end{thm}

Before proceeding with the proof the statement of Theorem \ref{thm:absolutecontinuity} calls for a comment. As discussed already in Section \ref{sec:outline}, our approach to establish Theorem \ref{thm:absolutecontinuity} is via the transversality method developed within the framework of Bernoulli convolutions. That is, we will rely on the transversality estimate \eqref{eq:Transversality} which is known to hold only on closed intervals $J\subset(1/2,2^{-2/3})$, thus the restriction in the statement of Theorem \ref{thm:absolutecontinuity}.

\subsection{Outline and proof of Theorem \ref{thm:absolutecontinuity}}\label{sec:outline proof AC}
For the proof of Theorem \ref{thm:absolutecontinuity} we would like to follow the same scheme we presented in Section \ref{sec:transversalityaffine} for the affine case. 
As we will see below, one has to face several annoying technicalities when adapting these ideas. Hence, we find it convenient to first reduce Theorem \ref{thm:absolutecontinuity} to a more technical statement, Proposition \ref{prop:ACtechnicalnonlinear}, which will be proved in the rest of this section.
\medskip

Fix any $y\in [0,1]$ (resp.  $x\in [0,1]$)  and let $\eta_{y,\lambda}\in \mathcal M(\gamma_y)$ (resp. $\eta_{x,\lambda}\in \mathcal M(\gamma_x)$)  be the  measure 
\begin{equation}\label{eq: final projection measure}
\eta^-_{y,\lambda}=(\Pi_{y,\lambda}^-)_*\   \mu^-\qquad\qquad (\text{resp. } \eta^+_{x,\lambda}=(\Pi_{x,\lambda}^+)_*\   \mu^+
\end{equation}
with 
\[
\Pi_{y,\lambda}^-:\{1,2\}^\mathbb N\to \gamma_y\subset D\qquad\qquad (\text{resp. }\Pi_{y,\lambda}^+:\{1,2\}^\mathbb N\to \gamma_x\subset D)
\]
 being the coding map in \eqref{eq:codingnonlinearproof} (resp. \eqref{eq:coding map x}). Proceeding as in Section \ref{sec:transversalityaffine}, for  $\lambda\in J$ we define 
\[
G^-(y,\lambda):=\int_{\gamma_y} D(\eta^-_{y,\lambda},x)\ \mathrm d \eta^- _{y,\lambda}(x),\qquad\qquad G^+(x,\lambda):=\int_{\gamma_x} D(\eta^+_{x,\lambda},y)\ \mathrm d \eta^+ _{x,\lambda}(y)
\]
where 
\[
D(\eta^-_{y,\lambda},x)=\liminf_{r\to 0}\frac{1}{2r} \eta^-_{y,\lambda}(B_{r}(x))\qquad\qquad  D(\eta^+_{x,\lambda},y)=\liminf_{r\to 0}\frac{1}{2r} \eta^+_{x,\lambda}(B_{r}(y)).
\]

\begin{prop}\label{prop:ACtechnicalnonlinear} 
Fix any closed interval $J\subset (1/2,2^{-2/3})$. Then, the function $(y,\lambda)\mapsto G^-(y,\lambda)$ (resp.  $(x,\lambda)\mapsto G^+(x,\lambda)$ ) is measurable and
\begin{equation}\label{eq:Upperbound Lowerderivative Final}
  \mathcal S^-=  \sup_{y\in[0,1]} \int_J G^-(y,\lambda)\ \mathrm d\lambda <\infty\qquad\qquad (\text{resp.  }   \mathcal S^+=  \sup_{x\in[0,1]} \int_J G^+(x,\lambda)\ \mathrm d\lambda <\infty).
\end{equation}
\end{prop}
\medskip

The finite bound \eqref{eq:Upperbound Lowerderivative Final} (compare to estimate \eqref{eq:finiteS} in Section \ref{sec:transversalityaffine})   is the heart of the proof of Theorem \ref{thm:absolutecontinuity}. Making use of the definition \eqref{eq: final projection measure}  one can easily see that the first inequality in  \eqref{eq:Upperbound Lowerderivative Final} holds if we obtain a (finite) upper bound, uniform in $y\in[0,1]$ for 
\begin{equation}\label{eq:lowerderivativeheuristic}
\liminf_{r\to 0} \frac{1}{2r} \int_J \int_{\{1,2\}^\mathbb N}\int_{\{1,2\}^\mathbb N} \mathbbm 1_{\{|\Pi^-_{y,\lambda}(\omega)-\Pi^-_{y,\lambda}(\eta)|\leq r\}} \ d\mu^-(\omega)d \mu^-(\eta) \mathrm d\lambda.
\end{equation}
The proof of Proposition \ref{prop:ACtechnicalnonlinear} occupies the rest of Section \ref{sec:typicalAC} (Sections \ref{sec: bdd distortion}-\ref{sec:proof AC}). Before giving the proof of Proposition \ref{prop:ACtechnicalnonlinear}, we first show how to complete the proof of Theorem \ref{thm:absolutecontinuity}.

\medskip

\begin{proof}[Proof of Theorem \ref{thm:absolutecontinuity}]
The proof is divided into two parts. We start by observing that \eqref{eq: unstable manifolds are almost vertical} implies that for all $x,y\in D$ and all $r>0$
\[
\Pi^{-1}_\lambda(B^2_r(x,y))\subset \{\omega^-\in \Omega^{-}\colon |\gamma^u_{\omega^-}-x|\leq (1+C\epsilon)r\} \times  \{\omega^+\in \Omega^{+}\colon |\gamma^s_{\omega^+}-y|\leq (1+C\epsilon)r\}. 
\]
Hence
\begin{equation}\label{eq:uppercontrol by product measure}
\nu_\lambda (B^2_r(x,y))\leq \eta^-_{y,\lambda} (B_{(1+C\epsilon)r})\ \eta^+_{x,\lambda} (B_{(1+C\epsilon)r}).
\end{equation}
We now make use of Proposition \ref{prop:ACtechnicalnonlinear} to  control the right hand side in \eqref{eq:uppercontrol by product measure}.   Fix $y\in [0,1]$, by \eqref{eq:Upperbound Lowerderivative Final} we have 
    \[
    D(\eta^-_{y,\lambda},x)<\infty\qquad\qquad \text{for } \mathrm{Leb}_J\text{-a.e. }\lambda\quad \text{and }\quad \nu_{y,\lambda}\text{-a.e. } x.
    \]
    Hence, using Vitali's covering theorem and Lebesgue's differentiation theorem we deduce the existence of a set $\mathcal U\subset J\times[0,1]$ of full measure such that for any pair $(\lambda,y)\in \mathcal U$ there exists a function $a_{\lambda}(\cdot,y)\in L^2([0,1])$  such that 
    \begin{equation}\label{eq:density final final}
    \frac{\mathrm d \eta^-_{y,\lambda}}{\mathrm d x}=a_{\lambda}(\cdot,y).
    \end{equation}
It then follows that for $(\lambda,y)\in \mathcal U$ we must have $G(y,\lambda)=\lVert \tilde a_{\lambda}(\cdot,y)\rVert^2_{L^2}$. In particular, $y\mapsto \lVert a_{\lambda}(\cdot,y)\rVert^2_{L^2}$ is measurable and, using \eqref{eq:Upperbound Lowerderivative Final} plus Tonelli's theorem, we conclude that 
\begin{equation}\label{eq:L2estimatepre}
\int_0^1 \lVert a_{\lambda}(\cdot,y)\rVert^2_{L^2} \mathrm d y<\infty \qquad\qquad \text{for a.e. }\lambda \in J.
\end{equation}
That is, $a_\lambda\in L^2(D)$. Reasoning analogously we deduce the existence of $b_\lambda\in L^2(D)$ such 
\[
\frac{\mathrm d \eta^+_{x,\lambda}}{\mathrm d y}=b_{\lambda}(x,\cdot).
\]
Fix now any $\lambda\in J$ for which both $a_\lambda,b_\lambda\in L^2(D)$. We can then define the maximal functions
\[
Ma_\lambda(x,y)=\sup_{r>0}\frac 1{2r}\int_{x-r}^{x+r} a_\lambda(s,y)\mathrm ds \qquad\qquad Mb_\lambda(x,y)=\sup_{r>0}\frac 1{2r}\int_{y-r}^{y+r} b_\lambda(x,s)\mathrm ds.
\]
Hardy-Littlewood's maximal inequality shows that $M{a_\lambda},M{b_\lambda}\in L^2(D)$ so, by Cauchy-Schwartz
\[
g_\lambda=M{a_\lambda}M{b_\lambda}\in L^1(D).
\]
On the other hand, inequality \eqref{eq:uppercontrol by product measure} implies that for a.e. $z\in D$ and any ball $B_r(z)\subset D$ we have 
\[
\nu_\lambda (B_r(z))\leq C r^2  g_\lambda(z)
\]
what implies that  $\nu_\lambda \ll \mathrm{Leb}_D$.
\end{proof}

\medskip

\medskip

\subsection{Organization of the proof of Proposition \ref{prop:ACtechnicalnonlinear}}

Concerning the proof of Proposition \ref{prop:ACtechnicalnonlinear}, on a first reading, we recommend the reader to skip directly to Section \ref{sec:proof AC}, where Proposition \ref{prop:ACtechnicalnonlinear}  is proved making use of some auxiliary lemmas, to get an overall feeling of the main ideas and only afterwards come back to the more technical Sections \ref{sec: bdd distortion}-\ref{sec:finiteenergy}, containing all the estimates used in Section \ref{sec:proof AC}.

\subsection{Bounded distortion}\label{sec: bdd distortion}

We first need to gather a couple of technical results. 

\begin{lem}[Bounded distortion]\label{lem:bounded distortion}    Given $\gamma\in\Gamma$ and $\tau\in\mathcal W$, write
$n=|\tau|$, 
for $z\in\mathtt F^{-n}_{\tau,\lambda}(\gamma\cap V_{\tau,\lambda})$ let
\[
z_{0,\lambda}=z,\qquad
\gamma_{0,\lambda}=\mathtt F^{-n}_{\tau,\lambda}(\gamma\cap V_{\tau,\lambda}),
\]
and, for $j=0,\dots,n-1$, set
\[
z_{j+1,\lambda}=\mathtt F_{\tau_j,\lambda}(z_{j,\lambda}),
\qquad
\gamma_{j+1,\lambda}
=\mathtt F_{\tau_j,\lambda}(\gamma_{j,\lambda}).
\]
Define
\begin{equation}\label{eq: unstable jacobian orbit segment}
\mathrm{Jac}_{\gamma,\tau,\lambda}(z)
=
\prod_{j=0}^{n-1}
\mathrm{Jac}_{\gamma_{j+1,\lambda}}
\mathtt F_{\tau_j,\lambda}^{-1}(z_{j+1,\lambda}),
\end{equation}
where $\mathrm{Jac}_{\gamma}G(z)$ denotes the arclength
Jacobian of a diffeomorphism $G$ restricted to the curve $\gamma$ at the point $z\in \gamma$.

    Then, there exists a  $C^1$ function $\psi_{\tau,\gamma}:\mathtt F^{-n}_{\tau,\lambda}(\gamma\cap V_{\tau,\lambda})\to [1/2,2]$   such that for any $z\in \mathtt F^{-n}_{\tau,\lambda}(\gamma\cap V_{\tau,\lambda})$
    \begin{equation}\label{eq:bdd distortion}
    \mathrm{Jac}_{\gamma,\tau,\lambda}(z)=\frac{1}{\ell_\gamma (\gamma\cap V_{\tau,\lambda})} \ \psi_{\tau,\gamma,\lambda}(z)
    \end{equation}
where $\ell_\gamma$ stands for the arclength measure on $\gamma$.  Moreover, there exists $C>0$ independent of $\tau,\gamma$ and $\lambda $ such that
    \begin{equation}\label{eq:uniform C1 bound psi Jacobian}
    |\psi_{\tau,\gamma,\lambda}|_{C^1}\leq C.
    \end{equation}
\end{lem}

\begin{proof}
We estimate separately the $C^0$ and $C^1$ norm.

\noindent\textit{$C^0$ estimate:} We start by observing that \eqref{eq:mathtt is epsilon good} implies that for  for any $\gamma\in \Gamma$ and any $z\in \gamma$
\begin{equation}\label{eq:unstable Jacobian}
\mathrm{Jac}_{\gamma}\mathtt F_i^{-1}(z)=\lambda+O_{C^1}(\epsilon)
\end{equation}
Since, by construction,  for all $j=0,\dots,n$ the curve $\gamma_j$ corresponds to  a segment of a curve in $\Gamma$, the estimate  \eqref{eq:unstable Jacobian} implies  that  for all $j=0,\dots,n$
\begin{equation}\label{eq:intermediate distance}
\operatorname{dist}(z_j,z_j')
\leq C(\lambda+C\epsilon)^j.
\end{equation}
Hence it follows from \eqref{eq:intermediate distance}  that for all $j=0,\dots,n-1$ 
\begin{equation}\label{eq:intermediate unstable jac}
\left|\log \mathrm{Jac}_{\gamma_{j+1}}\mathtt F^{-1}_{\tau_j}(z_{j+1})-\log \mathrm{Jac}_{\gamma_{j+1}}\mathtt F^{-1}_{\tau_j}(z_{j+1}')
\right|\leq C\epsilon(\lambda+C\epsilon)^{j+1}.
\end{equation}
Finally, to estimate the variation of \eqref{eq: unstable jacobian orbit segment} we write 
  \[
\log\frac{\mathrm{Jac}_{\gamma,\tau}(z)}
     {\mathrm{Jac}_{\gamma,\tau}(z')}=
\sum_{j=0}^{n-1}\left(\log\mathrm{Jac}_{\gamma_{j+1}}\mathtt F^{-1}_{\tau_j}(z_{j+1})-\log\mathrm{Jac}_{\gamma_{j+1}}\mathtt F^{-1}_{\tau_j}(z_{j+1}')\right),
\]
   and deduce from \eqref{eq:intermediate unstable jac} that 
   \[
   \left| 
   \log \frac{\mathrm{Jac}_{\gamma,\tau}(z)}{\mathrm{Jac}_{\gamma,\tau}(z')}\right|\leq C\epsilon\sum_{j=0}^{n-1}(\lambda+C\epsilon)^{j+1}\leq C\epsilon.
   \]
   This,  together with the fact that $\gamma\cap V_\tau$ gets map to a curve of length of order one,  implies the existence of a  $C^0$ function $\psi_{\tau,\gamma}:(\mathtt F_\tau^{|\tau|})^{-1}(\gamma\cap V_\tau)\to [1/2,2]$ with $|\psi_{\tau,\gamma}|_{C^0}\leq C$   such that for any $z\in (\mathtt F_\tau^{|\tau|})^{-1}(\gamma\cap V_\tau)$
    \begin{equation}\label{eq:bdd distortion proof}
    \mathrm{Jac}_{\gamma,\tau}(z)=\frac{1}{\ell_\gamma (\gamma\cap V_{\tau})} \  \psi_{\tau,\gamma}(z).
    \end{equation}
\medskip

\noindent\textit{$C^1$ estimate:} 
Parametrize $\gamma_0$ by its
horizontal coordinate $x$, and write
\[
z_j(x)=(x_j(x),\gamma_j(x_j(x))),\qquad x_0(x)=x.
\]
Since
\[
x_{j+1}
=
\mathtt F_{\tau_j,x}(x_j,\gamma_j(x_j)),
\]
we have
\begin{equation}\label{eq:contraction x coordinate}
\left|\frac{\mathrm d x_j}{\mathrm d x}\right|
\leq C(\lambda+C\epsilon)^j,
\qquad j=0,\dots,n.
\end{equation}
Define
\[
b_j(s)=
\log\mathrm{Jac}_{\gamma_{j+1}}\mathtt F^{-1}_{\tau_j}
(s,\gamma_{j+1}(s)).
\]
The uniform $C^2$ bounds on the admissible curves imply
\[
|b_j'|_{C^0}\leq C.
\]
Therefore,
\[
\left|
\frac{\mathrm d}{\mathrm dx}
\log\mathrm{Jac}_{\gamma,\tau}(x,\gamma_0(x))
\right|=\left|\sum_{j=0}^{n-1}
b_j'(x_{j+1}(x))\frac{\mathrm d x_{j+1}}{\mathrm dx}\right|\leq C\sum_{j=0}^{n-1}
(\lambda+C\epsilon)^{j+1}
\leq C.
\]
Together with \eqref{eq:bdd distortion proof}, this gives
the claimed uniform $C^1$ bound on $\psi_{\tau,\gamma}$.
\end{proof}
\medskip

\subsubsection*{Regularity of the Jacobian w.r.t $\lambda$}
\begin{lem}\label{lem:estimate for kappa dependence on lambda}
Consider the setting of Lemma \ref{lem:bounded distortion}. Then,
\begin{equation}\label{eq: difference of kappas}
 |\partial_\lambda \log (\mathrm{Jac}_{\gamma,\tau,\lambda})|_{C^0}\leq C|\tau|.
\end{equation}
\end{lem}
\begin{proof}
    Let $\{\gamma_{j,\lambda}\}\subset \Gamma$ be the collection of curves in Lemma \ref{lem:bounded distortion}. In the following we will exploit that
    \begin{equation}\label{eq: estimate curves gamma lambda}
   [\gamma_{j,\lambda}']_{C^1}, [\gamma_{j,\lambda}']_{C^2}\leq C\epsilon\qquad\qquad \max_{0\leq j\leq |\tau|} |\partial_\lambda \gamma_{j,\lambda}|_{C^1}\leq C.
    \end{equation}
    The latter inequality it is deduced by writing $\partial_\lambda\gamma_{j,\lambda}=
\partial_\lambda\mathcal G_{\tau_j,\lambda}(\gamma_{j+1,\lambda})+D\mathcal G_{\tau_j,\lambda}(\gamma_{j+1,\lambda})[\partial_\lambda\gamma_{j+1,\lambda}]$ and making use of inequalities \eqref{eq: parameter dependence graph transform} and \eqref{eq:frechet differential}.
\medskip

We now can understand the dependence on $\lambda$ of the orbit segment $\{z_{j,\lambda}\}_j\subset\{\gamma_{j,\lambda}\}_j$. To do so we write
\[
z_{j,\lambda}(x)
=
(x_{j,\lambda}(x),
 \gamma_{j,\lambda}(x_{j,\lambda}(x))),
\qquad x_{0,\lambda}(x)=x.
\]
Then
\[
x_{j+1,\lambda}
=
\mathtt F_{\tau_j,\lambda,x}
(x_{j,\lambda},\gamma_{j,\lambda}(x_{j,\lambda})).
\]
Differentiating this relation with respect to $\lambda$ and making use of \eqref{eq:mathtt is epsilon good} and \eqref{eq: estimate curves gamma lambda} 
\begin{equation}\label{eq: estimate x coord gamma lambda}
|\partial_\lambda x_{j+1,\lambda}|\leq C+C\epsilon+(\lambda+C\epsilon)|\partial_\lambda x_{j,\lambda}|\qquad \Longrightarrow\qquad \max_{0\leq j\leq n} |\partial_\lambda x_{j,\lambda}|\leq C.
\end{equation}
Together, \eqref{eq: estimate curves gamma lambda} and \eqref{eq: estimate x coord gamma lambda} imply that $z_{j,\lambda}=(x_{j,\lambda},\gamma_{j,\lambda}(x_{j,\lambda}))$ satisfies
\begin{equation}\label{eq: estimate points lambda}
\max_{0\leq j\leq n} |\partial_\lambda z_{j,\lambda}|\leq C.
\end{equation}

Finally, having understood the orbit dependence on $\lambda$ we can now study that of the weighted Jacobian \eqref{eq: unstable jacobian orbit segment}. To do so, we define 
\[
b_{j,\lambda}(x)
=
\log\mathrm{Jac}_{\gamma_{j+1,\lambda}}
\mathtt F_{\tau_j,\lambda}^{-1}
(x,\gamma_{j+1,\lambda}(x)).
\]
and note that (use \eqref{eq:mathtt is epsilon good}, \eqref{eq: estimate curves gamma lambda} and \eqref{eq: estimate points lambda}) 
\begin{equation}\label{eq:bdd log jacobian lambda}
\left|
\partial_\lambda
\bigl(b_{j,\lambda}(x_{j+1,\lambda}(x))\bigr)
\right|\leq C.
\end{equation}
Therefore,
\[
\left|
\partial_\lambda\log\mathrm{Jac}_{\gamma,\tau,\lambda}(x,\gamma_{0,\lambda}(x))\right|\leq\sum_{j=0}^{n-1}\left|
\partial_\lambda\bigl(b_{j,\lambda}(x_{j+1,\lambda}(x))\bigr)
\right|\leq Cn
\]
as we wanted to show.
\end{proof}

\medskip

\subsection{The transversality estimate}
We now want to reproduce the transversality argument (see Section \ref{sec:outline}) in the nonlinear setting. To that end, given $y\in [0,1]$ and $\omega,\eta\in \{1,2\}^\mathbb N$ we let 
\begin{equation}\label{eq:delta omega eta}
   \Delta_{y,\lambda}(\omega,\eta)=\Pi_{y,\lambda}(\omega)-\Pi_{y,\lambda}(\eta).
\end{equation}
where $\Pi_{y,\lambda}$ is the coding map in \eqref{eq:codingnonlinearproof}. 
\subsubsection*{Transversality of degree $\beta$ in the affine case} We first note that in the affine setting the function 
\[
\Delta_{\lambda}^{\mathrm{aff}}(\omega,\eta)=\Pi^{\mathrm{aff}}_{\lambda}(\omega)-\Pi^{\mathrm{aff}}_{\lambda}(\eta).
\]
does not depend on $y\in (0,1)$. Moreover, as we discussed in Section \ref{sec:outline}, it was shown in  \cite{MR1386842} that there exists $\delta_{\mathrm{aff}}>0$ such that for any $\lambda\in J_0\subset (1/2,2^{-2/3})$ and any $\omega,\eta\in \{1,2\}^\mathbb N$ with $\omega_0\neq \eta_0$  
\begin{equation}\label{eq:affinetransv pf}\tag*{(T0)}
|\Delta^{\mathrm{aff}}_{\omega,\eta}(\lambda)|\leq \delta_{\mathrm{aff}}\qquad\qquad \Longrightarrow\qquad\qquad |\partial_\lambda\Delta^{\mathrm{aff}}_{\omega,\eta}(\lambda)|\geq \delta_{\mathrm{aff}}.
\end{equation}
On the other hand, for $\omega=\tau\omega'$ and $\eta=\tau\eta'$ with $\omega_0'\neq\eta_0'$ it is not difficult to observe that 
\begin{equation}\label{eq:heuristics transv degree beta}
\Delta_{\lambda}^{\mathrm{aff}}(\omega,\eta)=\lambda^{|\tau|}\Delta_{\lambda}^{\mathrm{aff}}(\omega',\eta')
\end{equation}
Fix now any $\beta>0$ and let $c_\beta>0$ a sufficiently small constant to be chosen.  Then, if 
\[
|\Delta_{\lambda}^{\mathrm{aff}}(\omega,\eta)|\leq c_\beta \lambda^{(1+\beta)|\tau|}
\]
\eqref{eq:heuristics transv degree beta} implies that
\[
|\partial_\lambda \Delta_{\lambda}^{\mathrm{aff}}(\omega,\eta)|=\lambda^{|\tau|}\left||\tau|\lambda^{-1}\Delta_{\lambda}^{\mathrm{aff}}(\omega',\eta')+\partial_\lambda \Delta_{\lambda}^{\mathrm{aff}}(\omega',\eta')\right|\geq \lambda^{|\tau|}(\delta_{\mathrm{aff}}-|\tau|\lambda^{-1}c_\beta\lambda^{\beta|\tau|}).
\]
Hence, if we choose $c_\beta>0$ sufficiently small so that for any $n\in \mathbb N$
\[
\delta_{\mathrm{aff}}-n\lambda^{-1}c_\beta \lambda^{\beta n}\geq \frac{\delta_{\mathrm{aff}}}{2}
\]
we have proved that 
\begin{equation}\label{eq:affinetransv pf beta}
|\Delta_{\lambda}^{\mathrm{aff}}(\omega,\eta)|\leq c_\beta \lambda^{(1+\beta)|\tau|}\qquad \Longrightarrow\qquad|\partial_\lambda \Delta_{\lambda}^{\mathrm{aff}}(\omega,\eta)|\geq c_\beta \lambda^{(1+\beta)|\tau|}.
\end{equation}
The property \eqref{eq:affinetransv pf beta} is called transversality of degree $\beta$ and was introduced by Peres and Schlag in \cite{MR1749437}. Remarkably, this property is flexible enough so that we can establish an analogue in our non-linear two dimensional setting. 
\medskip

For later convenience, given $\omega=\tau\omega'$ and $\eta=\tau\eta'$ with $\omega_0'\neq \eta_0'$ we define the metric in $\{1,2\}^\mathbb N$ given by
\begin{equation}\label{eq: affine metric}
d_{\mathrm{aff},\lambda}(\omega,\eta)=\lambda^{|\tau|}.
\end{equation}
Making use of this notation, \eqref{eq:affinetransv pf beta} reads 
\begin{equation}\label{eq:affinetransv pf beta good}\tag*{(T0$\beta$)}
|\Delta_{\lambda}^{\mathrm{aff}}(\omega,\eta)|\leq c_\beta \ d_{\mathrm{aff},\lambda}(\omega,\eta)^{(1+\beta)}\qquad \Longrightarrow\qquad|\partial_\lambda \Delta_{\lambda}^{\mathrm{aff}}(\omega,\eta)|\geq c_\beta\  d_{\mathrm{aff},\lambda}(\omega,\eta)^{(1+\beta)}.
\end{equation}
\medskip

\subsubsection*{First-symbol transversality} The first step in order to run the  transversality argument in the nonlinear setting is to check that first symbol transversality is an open property.
\begin{lem}\label{lem:first symbol transv}
    Fix any closed interval $J\subset(1/2,2^{-2/3})$. There exists $\delta>0$ such that for any $\epsilon>0$ small enough the following holds. Let $\omega,\eta\in\{1,2\}^\mathbb N$ be such that $\omega_0\neq \eta_0$. Then, for any $y\in[0,1]$ and any $\lambda\in J$
    \begin{equation}\label{eq:1 symbol transversality}\tag*{(T)}
    |\Delta_{y,\lambda}(\omega,\eta)|\leq \delta\qquad\qquad \Longrightarrow \qquad\qquad |\partial_\lambda \Delta_{y,\lambda}(\omega,\eta)|\geq \delta.
    \end{equation}
\end{lem}
\begin{proof}
    Let $\delta_{\mathrm{aff}}>0$ be as in \ref{eq:affinetransv pf}. We omit the dependence on $y,\lambda$. Then, for $\epsilon$ small enough,  by \eqref{eq:regularity coding} and the first symbol transversality estimate \ref{eq:affinetransv pf} 
     \[
    |\Delta(\omega,\eta)|\leq \frac {\delta_{\mathrm{aff}}}{2}\  \Longrightarrow \  | \Delta^{\mathrm{aff}}(\omega,\eta)|\leq \frac {\delta_{\mathrm{aff}}}{2}+C\epsilon\leq \delta_{\mathrm{aff}}
    \]
    implies
    \[
    |\partial_\lambda \Delta^{\mathrm{aff}}(\omega,\eta)|\geq \delta_{\mathrm{aff}}\ \Longrightarrow\ |\partial_\lambda \Delta(\omega,\eta)|\geq \delta_{\mathrm{aff}}-C\epsilon \geq \frac {\delta_{\mathrm{aff}}}{2}.\qedhere
    \]
\end{proof}
\medskip

\subsubsection*{Geometric and weighted distance}
The second step in the nonlinear extension is to define a suitable metric in the space of sequences $\{1,2\}^\mathbb N$. Motivated by Lemma \ref{lem:bounded distortion} given $\omega=\tau\omega'$ and $\eta=\tau\eta'$ with $\omega_0'\neq \eta_0'$ we  introduce the geometric distance
\begin{equation}\label{eq:geometric distance}
 d_{y,\lambda}(\omega,\eta)=\ell_{\gamma_y}(\gamma_y \cap V_{\tau,\lambda})
\end{equation}
The following lemma relates this metric to the affine metric \eqref{eq: affine metric}.
\begin{lem}\label{lem:comparingdistances}
    Fix any $\beta>0$ and a closed interval $J\subset(1/2,1)$. Then, for any  $\epsilon>0$ small enough (depending only on $\beta$ and $J$),  any $y\in[0,1]$, any  $\omega,\eta\in \{1,2\}^\mathbb N$ and any $\lambda\in J$ we have
    \begin{equation}\label{eq:comparingdistances}
    \frac 12  d_{y,\lambda}(\omega,\eta)^{1+\beta}\leq d_{\mathrm{aff},\lambda}(\omega,\eta)\leq  2  d_{y,\lambda}(\omega,\eta)^{\frac{1}{1+\beta}}.
    \end{equation}
    
\end{lem}
\begin{proof}
    From the proof of Lemma \ref{lem:bounded distortion} it is not difficult to check that if $\omega=\tau\omega'$ and $\eta=\tau\eta'$ with $\omega'_0\neq \eta'_0$ 
    \[
    \frac 12(\lambda-C\epsilon)^{|\tau|}\leq \ell_{\gamma_y}(\gamma_y \cap V_{\tau,\lambda})\leq 2(\lambda+C\epsilon)^{|\tau|}
    \]
    what automatically gives \eqref{eq:comparingdistances}. 
\end{proof}
\medskip

\subsubsection*{Technical lemmas}

As third step of our nonlinear extension we now collect a couple of technical lemmas which are needed to estimate the derivatives of \eqref{eq:delta omega eta} with respect to $\lambda$.

\begin{lem}\label{lem: first technical lemma transversality}
   Fix any closed interval $J\subset(1/2,1)$ and suppose that $\epsilon>0$ is small enough. Let $\omega,\eta\in\{1,2\}^\mathbb N$ of the form $\omega=\tau\omega'$ and $\eta=\tau\eta'$ with $\omega'_0\neq\eta_0'$. Then, for every $y\in (0,1)$ there exists $Y_\lambda(y)\in (0,1)$ such that 
    \[
   C^{-1}  d_{y,\lambda}(\omega,\eta)^{-1}    |\Delta_{y,\lambda}(\omega,\eta)|\leq  |\Delta_{Y_\lambda(y),\lambda}(\omega',\eta')|\leq C   d_{y,\lambda}(\omega,\eta)^{-1}    |\Delta_{y,\lambda}(\omega,\eta)|.
    \]
\end{lem}

\begin{proof}
We omit the dependence on $\lambda$ and denote by 
\[
\ell_{\tau,y}=\ell_{\gamma_y}(\gamma_y\cap V_\tau).
\]
Let $(x_\omega(y),y)\in \gamma_\omega^u$ and $(x_\eta(y),y)\in \gamma_\eta^u$, and define
\[
z_\omega(y)=(x_{\omega'}(Y_\omega(y)),Y_\omega(y))=\mathtt F_{\tau}^{-|\tau|}(x_\omega(y),y)\in \gamma^u_{\omega'}\qquad z_\eta(y)=(x_{\eta'}(Y_\eta(y)),Y_\eta(y))=\mathtt F_{\tau}^{-|\tau|}(x_\eta(y),y)\in \gamma^u_{\eta'}.
\]
Then, in view of Lemma \ref{lem:bounded distortion} (see in particular \eqref{eq:bdd distortion}), it follows from the fundamental theorem of calculus that 
    \[
    \lVert z_\omega(y)-z_\eta(y)\rVert\leq \left|\int_{x_\omega(y)}^{x_\eta(y)}\mathrm{Jac}_{\gamma_y,\tau} \mathtt F_\tau^{-|\tau|}(s,y) \mathrm ds\right|\leq \frac{2}{\ell_{\tau,y}}|\Delta_{y}(\omega,\eta)|.
    \]
    Therefore,
    \[
    |x_{\omega'}(Y_\omega(y))-x_{\eta'}(Y_\omega(y))|\leq  |x_{\omega'}(Y_\omega(y))-x_{\eta'}(Y_\eta(y))|+ |x_{\eta'}(Y_\eta(y))-x_{\eta'}(Y_\omega(y))|
    \]
    For the first term we clearly have 
    \[
    |x_{\omega'}(Y_\omega(y))-x_{\eta'}(Y_\eta(y))|\leq \lVert z_\omega(y)-z_\eta(y)\rVert\leq \frac{C}{\ell_{\tau,y}}\ |\Delta_{y}(\omega,\eta)|.
    \]
    On the other hand, for the second term we use \eqref{eq: unstable manifolds are almost vertical} to deduce that
    \[
    |x_{\eta'}(Y_\eta(y))-x_{\eta'}(Y_\omega(y))|\leq C\epsilon |Y_\eta(y)-Y_\omega(y)|\leq C\epsilon \lVert z_\omega(y)-z_\eta(y)\rVert\leq \frac{C\epsilon}{\ell_{\tau,y}}\ |\Delta_{y}(\omega,\eta)|. 
    \]
    Combining both estimates 
    \[
     |x_{\omega'}(Y_\omega(y))-x_{\eta'}(Y_\omega(y))|\leq \frac{C}{\ell_{\tau,y}}\ |\Delta_{y}(\omega,\eta)|
    \]
    as we wanted to show. The lower bound is obtained from an analogous computation.
\end{proof}
\medskip

\begin{lem}
Consider the setting of Lemma \ref{lem: first technical lemma transversality}. Then, there exists $C>0$ such that the following holds 
\begin{equation}\label{eq: second technical lemma transversality}
\left|\partial_\lambda\Delta_{y,\lambda}(\omega,\eta)\right|
\ge
C^{-1} d_{y,\lambda}(\omega,\eta)
\left|\partial_\lambda
\Delta_{Y(y),\lambda}(\omega',\eta')\right|
-
C  d_{y,\lambda}(\omega,\eta)
\left(
|\tau|\left|\Delta_{Y(y),\lambda}(\omega',\eta')\right|
+\epsilon
\right).
\end{equation}
\end{lem}

\begin{proof}
Write $n=|\tau|$ and
$\ell_{\tau,y,\lambda}=\ell_{\gamma_y}(\gamma_y\cap V_{\tau,\lambda})$.
Parametrize the pulled-back curve by
\[
\mathtt F_{\tau,\lambda}^{-n}(\gamma_y\cap V_{\tau,\lambda})
=
\{(u,g_{\tau,y,\lambda}(u)):u\in[0,1]\},
\]
and define
\[
K_{\tau,y,\lambda}(x)
=
\pi_x\mathtt F_{\tau,\lambda}^{n}
(x,g_{\tau,y,\lambda}(x)).
\]
We start by recalling the estimates
\begin{equation}\label{eq:consequences bdd distortion}
C^{-1}\ell_{\tau,y,\lambda}
\leq K_{\tau,y,\lambda}'(u)
\leq C\ell_{\tau,y,\lambda},
\qquad
|K_{\tau,y,\lambda}''|_{C^0}
\leq C\ell_{\tau,y,\lambda},
\qquad
|\partial_\lambda K_{\tau,y,\lambda}'|_{C^0}
\leq Cn\ell_{\tau,y,\lambda},
\end{equation}
which follow from
Lemmas \ref{lem:bounded distortion} and
\ref{lem:estimate for kappa dependence on lambda},
together with the uniform $C^2$ bound on
$g_{\tau,y,\lambda}$ and the uniform $C^1$ bound on
$\partial_\lambda g_{\tau,y,\lambda}$. We now write
\begin{equation}\label{eq:decomposition technical}
\Delta_{Y_{\omega,\lambda}(y),\lambda}(\omega',\eta')
=
\underbrace{
x_{\omega',\lambda}(Y_{\omega,\lambda}(y))
-
x_{\eta',\lambda}(Y_{\eta,\lambda}(y))
}_{A_{\omega,\eta}(\lambda,y)}
+
\underbrace{
x_{\eta',\lambda}(Y_{\eta,\lambda}(y))
-
x_{\eta',\lambda}(Y_{\omega,\lambda}(y))
}_{E_{\omega,\eta}(\lambda,y)}    
\end{equation}
and estimate the derivatives of the two terms with respect
to $\lambda$. To do so we first introduce some notation. Using the notation of the preceding lemma, set
\[
u_\lambda=x_{\omega',\lambda}(Y_{\omega,\lambda}(y)),
\qquad
v_\lambda=x_{\eta',\lambda}(Y_{\eta,\lambda}(y))
\]
so
\begin{equation}\label{eq: Delta in terms of K}
\Delta_{y,\lambda}(\omega,\eta)=
K_{\tau,y,\lambda}(u_\lambda)-K_{\tau,y,\lambda}(v_\lambda).
\end{equation}
Using \eqref{eq:regularity coding},
$|g_{\tau,y,\lambda}'|\leq1$ and the uniform bound on
$\partial_\lambda g_{\tau,y,\lambda}$, one may check that
\[
|\partial_\lambda u_\lambda|+
|\partial_\lambda v_\lambda|+
|\partial_\lambda Y_{\omega,\lambda}(y)|+
|\partial_\lambda Y_{\eta,\lambda}(y)|
\leq C.
\]

\medskip

\noindent\textit{Estimate for $A_{\omega,\eta}$:}
Since $A_{\omega,\eta}=u_\lambda-v_\lambda$, differentiating
\eqref{eq: Delta in terms of K} gives
\begin{align*}
\partial_\lambda\Delta_{y,\lambda}(\omega,\eta)
={}&
K_{\tau,y,\lambda}'(u_\lambda)
\partial_\lambda A_{\omega,\eta}+
\left(
K_{\tau,y,\lambda}'(u_\lambda)
-
K_{\tau,y,\lambda}'(v_\lambda)
\right)\partial_\lambda v_\lambda+
(\partial_\lambda K_{\tau,y,\lambda})(u_\lambda)
-
(\partial_\lambda K_{\tau,y,\lambda})(v_\lambda).
\end{align*}
The first estimate in \eqref{eq:consequences bdd distortion}
implies
\[
|K_{\tau,y,\lambda}'(u_\lambda)
\partial_\lambda A_{\omega,\eta}|
\geq
C^{-1}\ell_{\tau,y,\lambda}
|\partial_\lambda A_{\omega,\eta}|.
\]
Using the second estimate and
$|\partial_\lambda v_\lambda|\leq C$, we obtain
\[
\left|
\left(
K_{\tau,y,\lambda}'(u_\lambda)
-
K_{\tau,y,\lambda}'(v_\lambda)
\right)\partial_\lambda v_\lambda
\right|
\leq C\ell_{\tau,y,\lambda}|A_{\omega,\eta}|.
\]
Finally, the third estimate gives
\[
|(\partial_\lambda K_{\tau,y,\lambda})(u_\lambda)
-
(\partial_\lambda K_{\tau,y,\lambda})(v_\lambda)|
\leq Cn\ell_{\tau,y,\lambda}|A_{\omega,\eta}|.
\]
Combining these estimates, for $n\geq1$ we conclude that
\[
|\partial_\lambda\Delta_{y,\lambda}(\omega,\eta)|
\geq
C^{-1}\ell_{\tau,y,\lambda}
|\partial_\lambda A_{\omega,\eta}|
-
Cn\ell_{\tau,y,\lambda}|A_{\omega,\eta}|.
\]

\medskip

\noindent\textit{Estimate for $E_{\omega,\eta}$:}
Since the two pulled-back points lie on the graph of
$g_{\tau,y,\lambda}$ and $|g_{\tau,y,\lambda}'|\leq1$,
\[
|Y_{\eta,\lambda}(y)-Y_{\omega,\lambda}(y)|
\leq |u_\lambda-v_\lambda|
=|A_{\omega,\eta}|.
\]
Together with $|\partial_y x_{\eta',\lambda}|\leq C\epsilon$,
this gives
\[
|E_{\omega,\eta}|
\leq C\epsilon|A_{\omega,\eta}|,
\qquad
|A_{\omega,\eta}|
\leq
C|\Delta_{Y_{\omega,\lambda}(y),\lambda}(\omega',\eta')|.
\]
Differentiating the definition of $E_{\omega,\eta}$,
we obtain
\begin{align*}
\partial_\lambda E_{\omega,\eta}
=&
(\partial_\lambda x_{\eta',\lambda})(Y_{\eta,\lambda}(y))
-
(\partial_\lambda x_{\eta',\lambda})(Y_{\omega,\lambda}(y))+
(\partial_y x_{\eta',\lambda})(Y_{\eta,\lambda}(y))
\partial_\lambda Y_{\eta,\lambda}(y)\\
&-
(\partial_y x_{\eta',\lambda})(Y_{\omega,\lambda}(y))
\partial_\lambda Y_{\omega,\lambda}(y).
\end{align*}
The last two terms are $O(\epsilon)$ by the small-slope
estimate and the uniform parameter bounds established above.
The first difference is also $O(\epsilon)$ by
\eqref{eq:regularity coding}, since
$\partial_\lambda x_{\eta',\lambda}$ is uniformly
$O(\epsilon)$-close to its affine counterpart, which is
independent of the slice coordinate. Therefore,
\[
|\partial_\lambda E_{\omega,\eta}|\leq C\epsilon.
\]
\textit{Conclusion:}
To conclude, differentiating the decomposition \eqref{eq:decomposition technical} gives
\[
\partial_\lambda A_{\omega,\eta}
+
\partial_\lambda E_{\omega,\eta}
=
\left.
\partial_\lambda\Delta_{Y,\lambda}(\omega',\eta')
\right|_{Y=Y_{\omega,\lambda}(y)}
+
\partial_Y\Delta_{Y_{\omega,\lambda}(y),\lambda}(\omega',\eta')
\partial_\lambda Y_{\omega,\lambda}(y).
\]
The last term is $O(\epsilon)$, and hence
\[
|\partial_\lambda A_{\omega,\eta}|
\geq
\left|
\left.
\partial_\lambda\Delta_{Y,\lambda}(\omega',\eta')
\right|_{Y=Y_{\omega,\lambda}(y)}
\right|-C\epsilon.
\]
Substituting this and the bound on $|A_{\omega,\eta}|$
into the estimate obtained above yields
\[
\begin{aligned}
|\partial_\lambda\Delta_{y,\lambda}(\omega,\eta)|
\geq{}&
C^{-1}\ell_{\tau,y,\lambda}
\left|
\left.
\partial_\lambda\Delta_{Y,\lambda}(\omega',\eta')
\right|_{Y=Y_{\omega,\lambda}(y)}
\right|-
C\ell_{\tau,y,\lambda}
\left(
n|\Delta_{Y_{\omega,\lambda}(y),\lambda}(\omega',\eta')|
+\epsilon
\right).
\end{aligned}
\]
Since $\ell_{\tau,y,\lambda}=d_{y,\lambda}(\omega,\eta)$,
this proves the claim with
$Y_\lambda(y)=Y_{\omega,\lambda}(y)$.
The empty-word case follows immediately by taking
$Y_\lambda(y)=y$.
\end{proof}
\medskip

\subsubsection*{Transversality of degree $\beta$}

Finally, we have all the ingredients needed to show that the coding map is ``transverse of degree $\beta$'' with respect to $\lambda$. 

\begin{prop}[Transversality of degree $\beta$]\label{prop:transv degree beta}
    Fix any $\beta>0$ and any closed interval $J\subset(1/2,2^{-2/3})$. Then, for any $c_\beta>0$ sufficiently small and any  $\epsilon\geq 0$ small enough  the following holds.  For any $\lambda\in J$, any $y\in [0,1]$ and any  pair $\omega,\eta\in\{1,2\}^\mathbb N$ 
   \begin{equation}\tag*{(T$\beta$)}\label{eq:transversality degree beta}
    |\Delta_{y,\lambda}(\omega,\eta)|\leq c_\beta\  d_{\mathrm{aff},\lambda}(\omega,\eta)^{1+\beta} \qquad \Longrightarrow\qquad  |\partial_\lambda\Delta_{y,\lambda}(\omega,\eta)|\geq c_\beta\  d_{\mathrm{aff},\lambda}(\omega,\eta)^{1+\beta}.
    \end{equation}
\end{prop}

\begin{proof}
Fix any $\beta>0$ and any $y\in (0,1)$ and suppose that 
\[
|\Delta_{y,\lambda}(\omega,\eta)|\leq c_\beta\  d_{\mathrm{aff},\lambda}(\omega,\eta)^{1+\beta}.
\]
Using \eqref{eq:comparingdistances} we may assume that 
\begin{equation}\label{eq: decreasing beta proof}
|\Delta_{y,\lambda}(\omega,\eta)|\leq c_\beta\  d_{y,\lambda}(\omega,\eta)^{1+\beta/2}.
\end{equation}
Then,  Lemma \ref{lem: first technical lemma transversality} implies that 
\[
|\Delta_{Y_{\lambda}(y),\lambda}(\omega',\eta')|\leq C  d_{y,\lambda}(\omega,\eta)^{-1}    |\Delta_{y,\lambda}(\omega,\eta)|\leq C c_\beta d_{y,\lambda}(\omega,\eta)^{\beta/2}.
\]
Hence, for $c_\beta>0$ small enough (uniformly for any $\omega,\eta$ and $y\in [0,1]$) and the  transversality estimate \ref{eq:1 symbol transversality} we must have 
\[
|\partial_\lambda \Delta_{Y_\lambda(y),\lambda}(\omega',\eta')|\geq \delta>0.
\]
Therefore, by inequality \eqref{eq: second technical lemma transversality} after decreasing if necessary the value of $c_\beta$,
\[
|\partial_\lambda \Delta_{y,\lambda}(\omega,\eta)|\geq   d_{y,\lambda}(\omega,\eta)\left( C^{-1}\delta-C c_\beta |\tau|   d_{y,\lambda}(\omega,\eta)^{\beta/2}-C\epsilon \right)\geq C  d_{y,\lambda}(\omega,\eta)\geq d_{\mathrm{aff},\lambda}(\omega,\eta)^{1+\beta},
\]
where in passing to the last inequality we have used  Lemma \ref{lem:comparingdistances}.
\end{proof}
\medskip

\subsection{Finite energy}\label{sec:finiteenergy}

Fix a positive constant $\sigma$. Then, given a metric $d$ and a measure $\mu$ on $\{1,2\}^\mathbb N$ we define their $\sigma$-\textit{correlation energy} as
\begin{equation}\label{eq:Energy}
        \mathcal E_{1+\sigma}(d,\mu)=\int_{\{1,2\}^\mathbb N} \int_{\{1,2\}^\mathbb N} d(\omega,\eta)^{-1-\sigma} \mathrm d\mu(\omega)\mathrm d\mu(\eta).
    \end{equation}
\begin{prop}\label{prop:finiteenergy}
 For  $\lambda\in  (1/2,1)$ let $\mu^-$ be the symmetric Bernoulli measure on $\Omega^-=\{1,2\}^\mathbb N$   and let  $d_{\mathrm{aff},\lambda}$ be the metric introduced in \eqref{eq: affine metric}.  For any closed interval $J\subset(1/2,1)$ there exists $\sigma >0$ such that 
    \begin{equation}\label{eq:finiteEnergy}
\sup_{\lambda\in J}\ \mathcal E_{1+\sigma}(d_{\mathrm{aff},\lambda},\mu^-)<\infty.
    \end{equation}
\end{prop}

\begin{proof}
We omit the dependence on $\lambda$, which plays no role in the discussion.  Clearly, 
\[
d_{\mathrm{aff}}(\omega,\eta)^{-1-\sigma}\leq \sum_{k=0}^\infty \sum_{|\tau|=k} \lambda^{-(1+\sigma)|\tau|}\  \mathbbm 1_{[\tau]}(\omega)\mathbbm 1_{[\tau]}(\eta).
\]
It then follows from Tonelli's theorem and the definition of the correlation energy that 
\begin{align*}
\mathcal E_{1+\sigma}(d_{\mathrm{aff}},\mu^-)\leq & \sum_{k=0}^\infty \sum_{|\tau|=k} \lambda^{-(1+\sigma)|\tau|}\ (\mu^-([\tau])^2= \sum_{k=0}^\infty \sum_{|\tau|=k}  \left(4\lambda^{1+\sigma}\right)^{-|\tau|}\leq \sum_{k=0}^\infty (2\lambda^{1+\sigma})^{-k}
\end{align*}
and the conclusion follows by taking $\sigma$ so that $2\lambda^{1+\sigma}>1$ for all $\lambda \in J\subset (1/2,1)$.
\end{proof}

\medskip

\subsection{Proof of Proposition \ref{prop:ACtechnicalnonlinear} }\label{sec:proof AC}
The joint measurability of $G^-(y,\lambda)$ follows exactly as in Section 4 of \cite{MR1386842}, using the joint continuity of the coding map
$(y,\lambda,\omega)\mapsto\Pi^-_{y,\lambda}(\omega)$. 
\medskip

It only remains to establish finiteness of \eqref{eq:Upperbound Lowerderivative Final}. In the following we omit the dependence on $y$ as it plays no role in the discussion. We want  to show that the for the one-parametric family of  measures $\{\eta^-_\lambda\}_\lambda$ defined in \eqref{eq: final projection measure}
\begin{equation}\label{eq:estimate S final0}
\mathcal S=\int_{J}\int_{\gamma_y} D(\eta^-_\lambda,x)\ \mathrm d\eta^-_\lambda(x) \mathrm d\lambda<\infty,
\end{equation}
where 
\[
D(\eta^-_\lambda, x):=\liminf_{r\to 0}\frac{1}{2r}\eta^-_\lambda (B_r(x)).
\]
Again, using Fatou's lemma and the definition \eqref{eq: final projection measure},  in order to establish \eqref{eq:estimate S final0} it is enough to show that 
\begin{equation}\label{eq:estimate S final}
\mathcal S\leq \liminf_{r\to 0}\frac{1}{2r} \int_J \int_{\{1,2\}^\mathbb N}\int_{\{1,2\}^\mathbb N} \mathbbm 1_{\{|\Pi^-_\lambda(\omega)-\Pi^-_\lambda(\eta)|\leq r\}} \ \mathrm d\mu^-(\omega)\mathrm d\mu^-(\eta) \mathrm d\lambda.
\end{equation}
It follows from  the transversality estimate \ref{eq:transversality degree beta} that for any $\beta>0$ and any $c_\beta>0$ sufficiently small this integral can be estimated by (here $\lambda_-=\min\{\lambda\in J\}$)
\[
\mathcal S\leq c_\beta^{-1}\int_{\{1,2\}^\mathbb N} \int_{\{1,2\}^\mathbb N} d_{\mathrm{aff},\lambda_-}^{-1-\beta}(\omega,\eta)\ \mathrm d\mu^-(\omega)\mathrm d\mu^-(\eta)=c_\beta^{-1} \mathcal E_{1+\beta}(d_{\mathrm{aff},\lambda_-},\mu^-)
\]
with $\mathcal E_{1+\beta}$ being the correlation energy introduced in \eqref{eq:Energy}. Hence, using Lemma \ref{prop:finiteenergy} it follows that $
\mathcal S<\infty$ (by taking  $\beta$ sufficiently small). The proof of Theorem \ref{thm:absolutecontinuity} is complete.
\medskip

\section{Equivalence to Lebesgue and proof of Theorem \ref{thm:main}}\label{sec:equivalenceLebesgue}

We have shown in Theorem \ref{thm:absolutecontinuity} that for a typical $\lambda\in J\subset (1/2,2^{-2/3})$ the measure $\nu_\lambda\in \mathcal M(D)$ in  \eqref{eq:pushforward definitive} is absolutely continuous with respect to Lebesgue. However, in order to conclude the proof of Theorem \ref{thm:main} we actually need to show that $\nu_\lambda\sim \mathrm{Leb}$. This is the main result of the present section.

\begin{thm}[Equivalence to Lebesgue]\label{thm:equivalence to Lebesgue}
  There exists a non-empty closed interval $J_0\subset (1/2,2^{-2/3})$ and  $\epsilon_0>0$ such that the following holds. For any $0<\epsilon\leq \epsilon_0$ and for a.e. $\lambda\in J_0$ the measure $\nu_\lambda\in \mathcal M(D)$ in \eqref{eq:pushforward definitive} is equivalent to Lebesgue measure on $D$.
\end{thm}

From this result it is straightforward to complete the proof of Theorem \ref{thm:main}.

\subsection{Proof of Theorem \ref{thm:main}}

Fix any $\lambda$ so that Theorem \ref{thm:equivalence to Lebesgue} holds.  Then,  $\nu_\lambda \simeq  \mathrm{Leb}_D$ and  the  proof of Theorem \ref{thm:main} is obtained by the very same argument used in Section \ref{sec:Hopfaffine} to conclude the proof of Theorem \ref{thm:mainaffine}.

\subsection{ Heuristics for the proof of Theorem \ref{thm:equivalence to Lebesgue}}

In the affine setting, thanks to the decoupling of vertical and horizontal dynamics we were able to use a well-known result by Mauldin-Simon \cite{MR1458276} to conclude equivalence. Unfortunately, we have not been able to generalize their argument to our two-dimensional hyperbolic picture. Our proof of equivalence follows a different route, detailed below. 
\medskip

More precisely, our argument goes as follows. By Theorem \ref{thm:absolutecontinuity} we know that for  typical $\lambda\in (1/2,2^{-2/3})$ the slice measures $\mu^-_{y,\lambda}$ and $\mu^+_{x,\lambda}$ in \eqref{eq: final projection measure} have densities
\begin{equation}\label{eq:slice densities again}
a_\lambda(\cdot,y)=\frac{\mathrm d\mu^-_{y,\lambda}}{\mathrm dx}\qquad\qquad b_\lambda(x,\cdot)=\frac{\mathrm d\mu^+_{x,\lambda}}{\mathrm dy}.
\end{equation}
\begin{rem}
    Recall that these slice measures can be used to control the measure $\nu_\lambda\in \mathcal M(D)$ in \eqref{eq:pushforward definitive} (See Section \ref{sec:outline proof AC}).
\end{rem}

We will show that there exists a inteval of values of $\lambda$ on which, typically, the sets $\{a_\lambda=0\}$ and $\{b_\lambda=0\}$ have zero Lebesgue measure. We describe the argument for $a_\lambda$, the argument for
$b_\lambda$ uses the analogous counting function below  for horizontal
cylinders. To do so, at any $z\in D$ we introduce the set
\[
\mathcal A_{\lambda,n} (z)=\{\tau\in \{1,2\}^n\colon z\in V_{\tau,\lambda}\}
\]
and the counting function
\begin{equation}\label{eq:countingfunction}
N_{\lambda,n}(z)=\# \mathcal A_{\lambda,n}(z).
\end{equation}
On one hand it is not difficult to see that this the average growth of this quantity can be controlled by\footnote{For the affine case we have the exact relation $\int_D N_{\lambda,n}^{\mathrm{aff}}(z)\mathrm dz=(2\lambda)^n$ for all $z\in D$.}
\begin{equation}\label{eq:heuristic upper control counting}
\int_D N_{\lambda,n}(z)\ \mathrm dz\leq (2\lambda+C\epsilon)^n.
\end{equation}
On the other hand, we will show that for a typical $\lambda$ there exists $s(\lambda)>0$ such that at Lebesgue almost any $z\in\{a_\lambda=0\}$ and sufficiently large $n$
\begin{equation}\label{eq:heuristic upper control counting for zero point}
N_{\lambda,n}(z)\leq (2\lambda^{1+s}+C\epsilon)^n.
\end{equation}
This (a priori) extra gain is associated to the fact that indeed, for a typical $\lambda$, the densities \eqref{eq:slice densities again} belong to the Sobolev space $H^s([0,1])$ (see Theorem \ref{thm:sobolev} for a precise statement).
\medskip

Suppose now that we could find some $\rho$ satisfying
\begin{equation}\label{eq:heuristic rho}
2\lambda^{1+s(\lambda)}\rho<1
\end{equation}
and an interval of values of $\lambda$ for which (note that in view of \eqref{eq:heuristic upper control counting} this $\rho$ must satisfy $(2\lambda+C\epsilon)\rho \geq 1$) 
\begin{equation}\label{eq:heuristic upper control counting special lambda}
\int_D N_{\lambda, n}(z)^{-1}\mathrm dz\leq C \rho^n.
\end{equation}
Then, a simple Borel-Cantelli argument  would show that for sufficienly small $\epsilon$ we must have $
\mathrm{Leb}\{a_\lambda=0\}=0$. The analogous statement would give that $\mathrm{Leb} \{b_\lambda=0\}=0$ for almost every $\lambda$ in that interval as well, establishing equivalence of the measure $\nu_\lambda\in \mathcal M(D)$ in \eqref{eq:pushforward definitive} to Lebesgue measure on $D$ (a complete argument is given in Section \ref{sec:proof equivalence}). 

\medskip

The proof of Theorem \ref{thm:equivalence to Lebesgue} boils down to make rigourous this heuristic discussion. It occupies the rest of this section and is divided in a number of steps. 
\medskip

 In Section \ref{sec: zero set sobolev density} we show that, if the density $a_\lambda(\cdot,y)$ is Sobolev with exponent $s$, the number of cylinder overlaps on the set $\{a_\lambda=0\}$ grows like \eqref{eq:heuristic upper control counting for zero point}. Then, in Section \ref{sec:Sobolev} we prove that for a typical $\lambda$ the density $a_\lambda(\cdot,y)$ belongs to a suitable Sobolev space.

 In Section \ref{sec:average growth cylinder overlaps} we find an interval of values of $\lambda$ for which  \eqref{eq:heuristic upper control counting special lambda} holds with a value of $\rho$ satisfying \eqref{eq:heuristic rho}.
 
 Combining the results of Sections \ref{sec: zero set sobolev density}-\ref{sec:Sobolev} with those in Section \ref{sec:average growth cylinder overlaps}, in Section \ref{sec:zero set} we show that typically $\mathrm{Leb}\{a_\lambda=0\}=0$. An analogous result is deduced for the density $b_\lambda$ associated to the vertical slice measures $\eta^+_{\lambda,x}$. Finally, in Section \ref{sec:proof equivalence} we complete the proof of Theorem \ref{thm:equivalence to Lebesgue}.

\subsection{Growth rate at the zero set of a Sobolev density}\label{sec: zero set sobolev density}

The main idea behind the proof of Theorem \ref{thm:equivalence to Lebesgue} is that the density $a_{\lambda}(\cdot,y)$ of $\eta^-_{\lambda,y}$ being regular enough forces the counting function $N_{\lambda,n}(z)$ to grow much more slowly than expected at points where the density vanishes.
Trivially we have the following.
\begin{lem}\label{lem:trivial bound counting function}
    There exists $C>0$ such that for every $y\in [0,1]$, every $\lambda\in (1/2,1)$, any $n\in\mathbb N$ and any $z=(x,y)\in D$ we have 
    \[
    \eta_{y,\lambda}^-(B_{(\lambda+C\epsilon)^n}(x))\geq 2^{-n} N_{\lambda,n}(z).
    \]
\end{lem}
We now analyze what happens to densities enjoying suitable regularity assumptions.

\begin{defn}
    Fix any $s\in (0,1)$. We say that a measurable function $h:[0,1]\to \mathbb R$ is in $H^s([0,1])$ if 
    \[
    \int_0^1 \int_0^1 \frac{|h(v)-h(u)|^2}{|v-u|^{1+2s}}\ \mathrm du \mathrm dv<\infty
    \]
\end{defn}

For densities in $H^s$ we deduce the following.

\begin{lem}\label{lem: growth on zero set}
    There exists $C>0$ such that the following holds. Suppose that for Lebesgue almost every $y\in [0,1]$ the density $a_\lambda(\cdot,y)=\mathrm d\eta^-_{y,\lambda}/\mathrm dx$   is in $H^s([0,1])$ for some $s\in(0,1)$. Then, for Lebesgue a.e. $z\in \{a_\lambda=0\}$ there exists $n_0(z)<\infty$ such that for all $n\geq n_0(z)$
    \begin{equation}\label{eq:small overlap zero sobolev}
    N_{\lambda,n}(z)\leq C\ (2\lambda^{1+s}+C\epsilon)^n.
    \end{equation}
\end{lem}
\begin{proof}
    Fix $y\in[0,1]$ for which $a_\lambda(\cdot,y)\in H^s([0,1])$. Let $r>0$ be small enough. We claim that for Lebesgue a.e. $z=(x,y)\in \{a_\lambda=0\}$ there exists $r_0(z)>0$ such that for all $0<r\leq r_0(z)$
    \begin{equation}\label{eq: sobolev small mass at vanish}
    \int_{x-r}^{x+r} a_\lambda(u,y)\mathrm du \leq r^{1+s}.
    \end{equation}
    Hence, \eqref{eq:small overlap zero sobolev} follows immediately from Lemma \ref{lem:trivial bound counting function} after choosing  $r=(\lambda +C\epsilon)^n$. To verify the claim \eqref{eq: sobolev small mass at vanish} it is enough to observe that the Sobolev norm being bounded implies that for a.e. $z=(x,y)\in \{a_\lambda=0\}$ we have 
    \[
   \int_0^1\frac{|a_\lambda(x,y)-a_\lambda(u,y)|^2}{|x-u|^{1+2s}}\mathrm du= \int_0^1\frac{a_\lambda^2(u,y)}{|x-u|^{1+2s}}\mathrm du<\infty.
    \]
    Hence, 
    \[
    \int_{x-r}^{x+r} a^2_\lambda(u,y)\mathrm du\leq r^{1+2s}\int_{x-r}^{x+r} \frac{a_\lambda^2(u,y)}{|x-u|^{1+2s}}\mathrm du=o(r^{1+2s})
    \]
    where in the last inequality we used  the absolute continuity of the Lebesgue integral. In particular, for $r$ small enough (depending a priori on $x$)
    \[
    \int_{x-r}^{x+r} a^2_\lambda(u,y)\mathrm du\leq r^{1+2s},
    \]
    so \eqref{eq: sobolev small mass at vanish} follows from Cauchy-Schwartz.
\end{proof}

\subsection{Sobolev regularity of the slice densities}\label{sec:Sobolev}

The following result guarantees that we can benefit from the gain in the previous section. 

\begin{thm}[Sobolev regularity]\label{thm:sobolev} Fix any closed interval $J\subset (1/2,2^{-2/3})$ and any $\delta>0$. For $\epsilon$ small enough, any $y\in [0,1]$, almost every $\lambda\in J$  and any
\[ 
0<s<s_*(\lambda)-\delta,\qquad\qquad s_*(\lambda)=\frac 12\left(\frac{\log 2}{\log(1/\lambda)}-1\right),
\]
the density 
\begin{equation}\label{eq:Sobolev estimate}
a_\lambda(\cdot,y)=\frac{\mathrm d\eta_{y,\lambda}^-}{\mathrm d x}\in H^s([0,1]).
\end{equation}  
\end{thm}

The proof of Theorem \ref{thm:sobolev} builds again on the transversality   techniques  developed in \cite{MR1749437,MR4385129} in the framework of Bernoulli convolutions.  We defer this proof to Appendix \ref{sec:sobolev appendix} as the main ideas are very similar to those in Section \ref{sec:typicalAC}. Note that from the definition of $s_*(\lambda)$ we have 
\[
\lambda^{1+s_*(\lambda)}=\sqrt{\lambda/2}.
\]
so \eqref{eq:heuristic rho} reads
\begin{equation}\label{eq:heuristic rho with slambda}
\sqrt{2\lambda} \ \rho<1.
\end{equation}

\subsection{Fast average growth of the counting function for suitable values of $\lambda$}\label{sec:average growth cylinder overlaps}

We now look for some $\rho$ satisfying \eqref{eq:heuristic rho} for which \eqref{eq:heuristic upper control counting special lambda} holds on a non-empty interval of values of $\lambda$.  At first sight \eqref{eq:heuristic upper control counting special lambda} and \eqref{eq:heuristic rho} only imply
\[
(2\lambda)^{-1}< \rho <(2\lambda)^{-1/2},
\]
so one might hope that the desired conclusion holds for all values of $\lambda\in (1/2,1)$. We show in Appendix \ref{sec: average growth proof appendix} that this is not the case, highlighting the subtlety of the construction.

\begin{prop}\label{prop:averagegrowth}
    There exists $C$, a non-empty interval interval $J_0\subset(1/2,2^{-2/3})$ on which the following holds. There exists $C>1$ and $\epsilon_0>0$ such that for any $0\leq \epsilon\leq \epsilon_0$ and any $\lambda\in J_0$ 
    \[
    \int_D N_{\lambda,n}(z)^{-1}\ \mathrm dz\leq C \rho_\lambda^n
    \]
    with $\rho_\lambda (2\lambda)^{1/2}<1$.
\end{prop}

The proof of Proposition \ref{prop:averagegrowth} is somewhat long so we present it in Appendix \ref{sec: average growth proof appendix} in order to do not interrupt our ongoing discussion.

\subsection{The zero set of the density a Sobolev density}\label{sec:zero set}

We now compare the average growth of the counting function (see Proposition \ref{prop:averagegrowth}) with its growth on the zero set (see Lemma \ref{lem: growth on zero set}) to conclude that the zero set must have zero Lebesgue measure.

\begin{prop}
    Let $J_0\subset (1/2,2^{-2/3})$ be the interval constructed in Proposition \ref{prop:averagegrowth}. For any $\epsilon$ small enough and almost every $\lambda\in J_0$ we have 
    \[
    \mathrm{Leb}_D\{a_\lambda=0\}=0.
    \]
\end{prop}

\begin{proof}
 For almost every $\lambda\in (1/2,2^{-2/3})$ Lemma \ref{lem: growth on zero set} and  Theorem \ref{thm:sobolev}  imply that for almost every $z\in E=\{a_\lambda=0\}$ there exists $n_0(z)<\infty$ such that 
 \[
 N_{\lambda,n}(z) \leq (\sqrt{2\lambda}+C\epsilon)^n
 \]
for any $n\geq n_0(z)$. For each $k\in \mathbb N$ let $E_k\subset E$ consist of the points satisfying
\[
N_{\lambda,n}(z)\leq (\sqrt{2\lambda}+C\epsilon)^n
\quad\text{for every }n\geq k.
\]
Clearly $\operatorname{Leb}(E\setminus\bigcup_kE_k)=0$. To complete the proof we now show that $\mathrm{Leb}(E_k)=0$ for any $k\in\mathbb N$.
\medskip

On one hand, it follows from the definition of $E_k$ that for $n\geq k$,
\[
\operatorname{Leb}(E_k)(\sqrt{2\lambda}+C\epsilon)^{-n}
\leq\int_DN_{\lambda,n}^{-1}\ \mathrm dz.
\]
On the other hand, it follows from Proposition \ref{prop:averagegrowth} that if we also restrict to $\lambda\in J_0$ we have 
\[
\int_DN_{\lambda,n}^{-1}\ \mathrm dz\leq C \rho^n.
\]
for some $\rho$ satisfying $\rho\sqrt{2\lambda}<1$.
Consequently,
\[
\operatorname{Leb}(E_k)
\leq C(\rho (\sqrt{2\lambda}+C\epsilon))^n
\longrightarrow0.
\]
as we wanted to show.
\end{proof}

\subsection{Proof of Theorem \ref{thm:equivalence to Lebesgue}}\label{sec:proof equivalence}

On one hand, an analogous argument for the vertical slice measures $\eta^+_{x,\lambda}$ in \eqref{eq:pushforward definitive} shows that for $J_0\subset(1/2,2^{-2/3})$, any $\epsilon$ small enough and almost every $\lambda\in J_0$ we have at the same time
\[
\operatorname{Leb}\{a_\lambda=0\}=0=\operatorname{Leb}\{b_\lambda=0\}
\]
On the other hand, proceeding as in Section \ref{sec:outline proof AC} one can show that
\[
\nu_\lambda (B^2_r(x,y))\geq \eta^-_{y,\lambda} (B_{(1+C\epsilon)r})\ \eta^+_{x,\lambda} (B_{(1+C\epsilon)r})
\]
so we are done.

\section{Transversality-torsion and Hopf boxes: proof of Proposition \ref{prop:TT and HB}}\label{sec:proof TT }

We want to show that it is possible to find a  $C^3$ function 
\[
\lambda\in J\mapsto s(\lambda)\in (0,1)
\]
such that for any $\epsilon>0$ small enough,  the family of maps $\{T_{\varepsilon(\epsilon, \lambda)},S_{\varepsilon(\epsilon, \lambda)}\}_{\lambda \in J}$ with 
\[
\varepsilon(\epsilon,\lambda)=\epsilon^r s(\lambda)
\]
(for some $r>0$ independent of $\epsilon$ and $\lambda$) is $\epsilon$-good. We do so in a number of steps. Although this precise statement was not proved in \cite{paradela2026robustlytransitivebehaviorsymplectic} (because it was not needed)  it follows plainly from the arguments in Section 3 of that paper, as we now show. 

\subsection{Local affine approximation}

We begin by obtaining a suitable normal form for iterates reproducing the so-called transversality-torsion mechanism (recall the discussion below Proposition \ref{prop:TT and HB} or see Section 2.2. in \cite{paradela2026robustlytransitivebehaviorsymplectic}).

\begin{lem}\label{lem:C1control}
    Fix any $C>0$. Then, there exists $\ell>0$ independent of $\varepsilon$ such that for any $0\leq n\leq C/\varepsilon$ the map 
    \[
    F_{n,\varepsilon}:=T^n_\varepsilon\circ S_\varepsilon: \mathbb A\to \mathbb A
    \]
    can be expressed as
    \begin{equation}\label{eq:localaffineapprox}
    F_{n,\varepsilon}:\binom{\varphi}{J}=\underbrace{\binom{[n\beta]}{0}}_{\bs b_n}+\underbrace{\begin{pmatrix}
    1+n\varepsilon ab&na\\
    \varepsilon b&1
\end{pmatrix}}_{A^{(n)}_\varepsilon}\binom{\varphi}{J}+\mathcal E_\varepsilon(\varphi,J)
    \end{equation}
   with $\mathcal E_\varepsilon=(\mathcal E_{\varepsilon,\varphi},\mathcal E_{\varepsilon,J})$ such that at any $(\varphi,J)\in \{|\varphi|\leq 1/2,\ |J|\leq \ell \varepsilon\}$
    \begin{equation}\label{eq:errorestimatesIFSapprox}
    \mathcal E_{\varepsilon,\varphi}=O(n\varepsilon |\varphi|^2,\varepsilon)\qquad\qquad \mathcal E_{\varepsilon,J}=O(\varepsilon|\varphi|^2,\varepsilon^2).
    \end{equation}
    Moreover, for these values of $(\varphi,J)$ we have $C^1$ estimates
    \begin{equation}\label{eq:errorestimatesIFSapproxC1}
    D\mathcal E_{\varepsilon}(\varphi,J)=\begin{pmatrix}O(n\varepsilon|\varphi|,\varepsilon)&O(1)\\ O(\varepsilon|\varphi|,\varepsilon^2)&O(\varepsilon) \end{pmatrix}
    \end{equation}
    $C^2$ estimates \begin{equation}\label{eq:errorestimatesIFSapprox C^2}
    D^2\mathcal E_{\varepsilon,\varphi}(\varphi,J)=\begin{pmatrix}O(n\varepsilon)&O(n\varepsilon)\\ O(n\varepsilon)&O(n) \end{pmatrix}\qquad\qquad D^2\mathcal E_{\varepsilon,J}(\varphi,J)=\begin{pmatrix}O(\varepsilon)&O(\varepsilon)\\ O(\varepsilon)&O(\varepsilon) \end{pmatrix},
    \end{equation}
    and $C^3$ estimates
\begin{equation}\label{eq:errorestimatesIFSapprox C^3}
    \begin{split}
|\partial_{\varphi^3}^3\mathcal E_{\varepsilon,\varphi}(\varphi,J)|,|\partial_{\varphi^2,J}^3\mathcal E_{\varepsilon,\varphi}(\varphi,J)|\leq &n\varepsilon,\quad |\partial_{\varphi, J^2}^3\mathcal E_{\varepsilon,\varphi}(\varphi,J)|\leq n,\quad |\partial_{J^3}^3\mathcal E_{\varepsilon,\varphi}(\varphi,J)|\leq n/\varepsilon\\
|\partial_{\varphi^3}^3\mathcal E_{\varepsilon,J}(\varphi,J)|,|\partial_{\varphi^2,J}^3\mathcal E_{\varepsilon,J}(\varphi,J)|\leq& \varepsilon,\quad |\partial_{\varphi, J^2}^3\mathcal E_{\varepsilon,J}(\varphi,J)|\leq n\varepsilon,\quad |\partial_{J^3}^3\mathcal E_{\varepsilon,J}(\varphi,J)|\leq n.
\end{split}
    \end{equation}
\end{lem}

\begin{proof}
    We proceed exactly as in the proof of Lemma 3.1. in \cite{guardia2025partiallyhyperbolicdynamics3body} to obtain that for any $n\in\mathbb N$ as in the statement and $(\varphi,J)\in \{|\varphi|\leq 1/2,\ |J|\leq \ell \varepsilon\}$ we have
    \begin{equation}\label{eq:recursion C1 control}
        \begin{split}
    \varphi_n(\varphi,J)=&\varphi+[n\beta]+na J+\varepsilon P_{\varepsilon,\varphi}(\varphi,J)+na\varepsilon P_{\varepsilon,J}(\varphi,J)+\rho_{\varepsilon,\varphi}^{(n)}(\varphi,J)\\
        J_n(\varphi,J)=&J+\varepsilon P_{\varepsilon,J}(\varphi,J)+\rho_{\varepsilon,J}^{(n)}(\varphi,J)     
        \end{split}
    \end{equation}
    with $\rho_\varepsilon^{(n)}=(\rho_{\varepsilon,\varphi}^{(n)},\rho_{\varepsilon,J}^{(n)})$ satisfying
    \[
    |\rho_{\varepsilon}^{(n)}|_{C^0}=\binom{O(n\varepsilon^2)}{O(n\varepsilon^3)}\qquad\qquad |D\rho_{\varepsilon}^{(n)}|_{C^0}=\begin{pmatrix}
        O(n\varepsilon^2)&O(n\varepsilon)\\
        O(n\varepsilon^3)&O(n\varepsilon^2)
    \end{pmatrix}.
    \]
    To obtain the estimates \eqref{eq:errorestimatesIFSapprox}-\eqref{eq:errorestimatesIFSapproxC1} we write 
    \begin{align*}
    P_{\varepsilon,J}(\varphi,J)=&P_{0,J}(\varphi,J)+O(\varepsilon)=\overbrace{\partial_\varphi P_{0,J}(0,0)}^b\  \varphi+O(|J|,(|\varphi|+|J|)^2)+O(\varepsilon)=b\ \varphi+O(|\varphi|^2,\varepsilon)\\
    \partial_\varphi P_{\varepsilon,J}(\varphi,J)=&\partial_\varphi P_{0,J}(\varphi,J)+O(\varepsilon)=\partial_\varphi P_{0,J}(0,0)\  +O(|\varphi|,|J|)+O(\varepsilon)=b+O(|\varphi|,\varepsilon)\\
   \partial_J P_{\varepsilon,J}(\varphi,J)=&O(1)
    \end{align*}
    and use that $|P_{\varepsilon,\varphi}|_{C^1}=O(1)$.
    \medskip
    
    In order to deduce the $C^2$ estimates, a simple computation shows that 
    \[
    D^2 T^n_{\varepsilon,\varphi}(\varphi,J)=\begin{pmatrix}O(n\varepsilon^2)&O(n\varepsilon)\\ O(n\varepsilon)&O(n) \end{pmatrix}\qquad\qquad D^2T^n_{\varepsilon,J}(\varphi,J)=\begin{pmatrix}0&0\\0&0\end{pmatrix}.
    \]
    Since $|S_{\varepsilon}-\mathrm{id}|_{C^2}\leq \varepsilon$ the estimates \eqref{eq:errorestimatesIFSapprox C^2} follow from the chain rule. The computation for the $C^3$ estimates is similar and it is left to the reader.
\end{proof}
We observe that for any $n\in \mathbb N$ the matrix $A^{(n)}_\varepsilon$  in \eqref{eq:localaffineapprox} is hyperbolic (it has  determinant one and its trace is larger than two). A simple computation shows that its eigenvalues $\{\alpha,1/\alpha\}\subset(0,\infty)$ are given by 
\begin{equation}\label{eq:eigenvalue affine thmA}
\alpha(n,\varepsilon)=\tilde\alpha(n\varepsilon ab)\qquad\qquad\text{where}\qquad\qquad  \tilde\alpha(s)=\frac{s+2-\sqrt{s^2+4s}}{2}
\end{equation}
For later convenience we observe that
\begin{equation}
   \tilde\alpha^{-1}(1/2,1)=(0,1/2)
\end{equation}
and define the inverse function $\tilde s:(1/2,1)\to (0,1/2)$  satisfying
\begin{equation}\label{eq:reparametrization}
\tilde \alpha(\tilde s(\lambda))=\lambda
\end{equation}
for all $\lambda\in (1/2,1)$.

\subsubsection*{Parameter dependence}
We will also need to understand the dependence of these maps on $\varepsilon$ (for a fixed $n\in\mathbb N$).

\begin{lem}
    Consider the setting of Lemma \ref{lem:C1control}. Then, at any $(\varphi,J)\in \{|\varphi|\leq 1/2,\ |J|\leq \ell \varepsilon\}$
    \begin{equation}\label{eq:error derivative varepsilon}
    \partial_\varepsilon \mathcal E_{\varepsilon}(\varphi,J)=\binom{O(n|\varphi|)}{O(1)}\qquad\qquad \partial_\varepsilon D\mathcal E_{\varepsilon}(\varphi,J)=\begin{pmatrix}
        O(1,n|\varphi|)&O(n)\\O(|\varphi|,\varepsilon)&O(1)
    \end{pmatrix}
    \end{equation}
    and
\begin{equation}\label{eq:error derivative varepsilon higher}
\partial_\varepsilon^2\mathcal E_\varepsilon
=
\begin{pmatrix}
O(n)\\
O(1)
\end{pmatrix} 
\ 
\partial_\varepsilon^2D\mathcal E_\varepsilon
=
\begin{pmatrix}
O(n)&O(n/\varepsilon)\\
O(1)&O(n)
\end{pmatrix}\ 
\partial_\varepsilon D^2\mathcal E_{\varepsilon,\varphi}
=
\begin{pmatrix}
O(n)&O(n)\\
O(n)&O(n/\varepsilon)
\end{pmatrix}
\ 
\partial_\varepsilon D^2\mathcal E_{\varepsilon,J}
=
\begin{pmatrix}
O(1)&O(1)\\
O(1)&O(n)
\end{pmatrix}
    \end{equation}
\end{lem}
\begin{proof}
    After some tedious but straightforward work one finds that $\rho_\varepsilon^{(n)}$ in \eqref{eq:recursion C1 control} also satisfies 
    \[
    |\partial_\varepsilon \rho_{\varepsilon}^{(n)}|_{C^0}=\binom{O(n\varepsilon)}{O(n\varepsilon^2)}\qquad\qquad |\partial_\varepsilon D\rho_{\varepsilon}^{(n)}|_{C^0}=\begin{pmatrix}
        O(n\varepsilon)&O(n)\\
        O(n\varepsilon^2)&O(n\varepsilon)
    \end{pmatrix}
    \]
    Differentiating \eqref{eq:recursion C1 control} and plugging these estimates yields \eqref{eq:error derivative varepsilon}. The similar argument gives the higher derivatives \eqref{eq:error derivative varepsilon higher}.
\end{proof}

\medskip

\subsection{A good range of iterates}
In what follows we fix any  closed interval 
\[
I\subset(0,1)
\]
take  $\epsilon>0$  and let $s\in I$. Then, motivated by the scaling in \eqref{eq:eigenvalue affine thmA}, we write (later we will substitute $\delta$ by a power of $\epsilon$ to agree with the notation in Proposition \ref{prop:TT and HB})
\[
\varepsilon_\delta(s)=\delta s
\]
and let
\[
N_\delta=\left[\frac{1}{\delta ab}\right], \qquad\qquad N_\kappa= \left[\frac{1}{\kappa \gamma}\right].
\]
where $\gamma>0$ is the Diophantine constant of $\beta\in \mathrm{DC}_\gamma$ (recall Assumption \ref{it:itemB2}) and $\kappa$ is a small parameter to be chosen and  which will determine the size of the domain in which we will make use of the linear approximation in Lemma \ref{lem:C1control}.
\begin{lem}\label{lem:diophantine approx}
Let  $
\mathcal N_{\delta,\kappa}=\{N_\delta,\dots,N_\delta+N_\kappa\}\subset\mathbb N$. Then, there exists $C_0>0$ and $\delta_0>0$ such that for any $0<\delta<\delta _0\kappa \ll 1$ the following holds:
\begin{enumerate}\label{lem:good iterates range}
    \item for any $n\in\mathcal N_{\delta,\kappa}$ and any $s\in I$ we have that  
    \[  \alpha(n,\varepsilon_\delta(s))=\tilde\alpha(s)+ \hat\alpha_{n,\delta,\kappa}(s)
    \]
    for some $C^3$ function $\hat \alpha_{n,\delta,\kappa}$ satisfying $|\partial_s \hat \alpha_{n,\delta,\kappa}|_{C^0(I)},|\partial^2_{s^2} \hat \alpha_{n,\delta,\kappa}|_{C^0(I)}\leq C \delta/\kappa$.
    \item The set $[n\beta]_{n\in\mathcal N_{\delta,\kappa}}\subset \mathbb T$ is $C_0\kappa$-dense.
\end{enumerate}
\end{lem}
The first item follows after a trivial computation. For the second one see Theorem VI in Chapter 5 of \cite{Cassels72}.
\medskip

\subsection{Uniform diagonalization}

Key to our construction is that the matrices $A^{(n)}_\varepsilon$ can be diagonalized simultaneously (up to a small error) for all $n\in \mathcal N_{\delta,\kappa}$. 

\begin{lem}\label{lem: uniform diagonalization}
  For any $\delta>0$ let   $Q_\delta\in GL(2,\mathbb R)$ be the matrix such that $Q_\delta^{-1} A^{(N_\delta)} Q_\delta=\mathrm{diag}(\alpha,1/\alpha)$. Then, for any $n\in \mathcal N_{\delta,\kappa}$ and any $s\in I$ we have  
    \begin{equation}\label{eq:uniform diagonal}
    Q^{-1}_\delta A^{(n)}_{\varepsilon_{\delta}(s)}Q_\delta(s)=\underbrace{\mathrm{diag}(\tilde\alpha(s),1/\alpha(s))}_{A(s)}+\widetilde A_{n,\delta,\kappa}(s)
    \end{equation}
with 
\begin{equation}\label{eq:error simultaneous diagonal}
|\widetilde A_{n,\delta,\kappa}|=O(\delta/\kappa)\qquad\qquad |\partial_s \widetilde A_{n,\delta,\kappa}|=O(\delta/\kappa) \qquad\qquad |\partial^2_{s^2} \widetilde A_{n,\delta,\kappa}|=O(\delta/\kappa).
\end{equation}
\end{lem}
One can easily understand this result by observing that for $n\in\mathcal N_{\delta,\kappa}$ we have $n/N_\delta-1=O(\delta/\kappa)$ and that for $n=N_\delta$ \eqref{eq:uniform diagonal} hols without error term. For details on the proof we refer the reader to the proof of Lemma 4.3. in \cite{paradela2026robustlytransitivebehaviorsymplectic} (see also Lemma  3.4. in \cite{guardia2025partiallyhyperbolicdynamics3body}).
\medskip

Lemma \ref{lem: uniform diagonalization} provides us with a uniform coordinate system
\[
\text{(Diagonalization)}\qquad\qquad\qquad \psi_\delta:\binom{\tilde x}{\tilde y}\mapsto \binom{\varphi}{J}=Q_\delta\binom{\tilde x}{\tilde y}
\]
to study (simultaneously) the family of maps $\{F_n\}_{n\in\mathcal N_{\delta,\kappa}}$. To make the non-linear errors negligible we will also find convenient to introduce a scaling
\[
\text{(Scaling)}\qquad\qquad\qquad \tilde\psi_\kappa:\binom{x}{ y}\mapsto \binom{\tilde x}{\tilde y}=R_\kappa\binom{x}{y}.
\]
where $R_\kappa=\mathrm{diag}(\sqrt\kappa,\sqrt\kappa)$ and $C_0$ is the constant in the second item of Lemma \ref{lem:diophantine approx}. For any $\delta>0$ and $\kappa>0$ denote their composition by 
\begin{equation}\label{eq: pre-good coordinate system}
\psi_{\delta,\kappa}=\psi_\delta\circ\tilde\psi_\kappa:\mathbb R^2\to \mathbb R^2.
\end{equation}
\begin{lem}\label{lem: uniform diag plus scaling}
Let $\psi_{\delta,\kappa}$ be defined as in \eqref{eq: pre-good coordinate system} and define
    \[
    \widehat F_{n,\delta,\kappa,s}=\psi_{\delta,\kappa}^{-1}\circ F_{n,\varepsilon_\delta(s)}\circ\psi_{\delta,\kappa}.
    \]
   Then, there exists $\delta_0>0$ such that for any $0<\delta\leq \delta_0  \kappa$ and any $n\in\mathcal N_{\delta,\kappa}$ we have 
    \begin{equation}
    \widetilde F_{n,\delta,\kappa,s}:\binom{x}{y}\mapsto \tilde{\bs b}_n+A_s\binom{x}{y}+\widehat{\mathcal{E}}_{n,\delta,\kappa,s}(x,y),
    \end{equation}
    with $\tilde{\bs b}_n=P_{\delta,\kappa}^{-1}\binom{[n\beta]}{0}$ and $\widehat{\mathcal E}_{n,\delta,\kappa,s}$ verifying 
    \begin{equation}\label{eq:pre C2 estimates error TT}
   |\widehat {\mathcal E}_{n,\delta,\kappa,s}|_{C^3(\mathcal D)},|  \partial_s \widetilde {\mathcal E}_{n,\delta,\kappa,s}|_{C^2(\mathcal D)},|  \partial_{s^2}^2 \widetilde {\mathcal E}_{n,\delta,\kappa,s}|_{C^1(\mathcal D)}=O(\sqrt\kappa, \delta/\kappa).
\end{equation}
\end{lem}

\begin{proof}
    In the following we will use withouth mentioning it that for $(x,y)\in \mathcal D$ we have $|\varphi(x,y)|=O(\sqrt\kappa)$. We now establish \eqref{eq:pre C2 estimates error TT} for 
    \[
    \widehat{\mathcal E}_{n,\delta,\kappa,s}(x,y)=\left( Q^{-1}_\delta A^{(n)}_{\varepsilon_{\delta}(s)}Q_\delta-A_s\right)\binom{x}{y}+\underbrace{\psi_{\delta,\kappa}^{-1}\circ \mathcal E_{n,\varepsilon_{\delta}(s)}\circ\psi_{\delta,\kappa}}_{\breve{\mathcal E}_{n,\delta,\kappa,s}}(x,y).
    \]
    \medskip
    
    \noindent\textit{$C^0$ estimates:}  We have shown in Lemma \ref{lem:C1control} that 
\[
|\mathcal E\circ\psi_{\delta,\kappa}|_{C^0}=\binom{O(\delta,\kappa)}{O(\delta\kappa,\delta^2)}
\]
Hence, 
\[
|\breve{\mathcal E}_{n,\delta,\kappa,s}|_{C^0}=O(\delta/\sqrt\kappa,\sqrt\kappa)
\]
which, combined with \eqref{eq:uniform diagonal}, implies that 
\[
|\widehat{\mathcal E}_{n,\delta,\kappa,s}|_{C^0}=O(\delta/\kappa,\sqrt\kappa).
\]
    \medskip

    \noindent\textit{$C^1$ estimates:} We have shown in Lemma \ref{lem:C1control} that 
\[
|(D\mathcal E_{\varepsilon_\delta(s)})\circ \psi_{\delta,\kappa}|_{C^0}=\begin{pmatrix}
    O(\sqrt\kappa,\delta)&O(1)\\
    O(\delta\sqrt\kappa,\delta^2)&O(\delta)
\end{pmatrix}
\]
Hence, 
\[
|D\breve{\mathcal E}_{n,\delta,\kappa,s}|_{C^0}=O(\sqrt\kappa,\delta)=O(\sqrt\kappa)
\]
which, combined with \eqref{eq:uniform diagonal}, implies that 
\[
|D\widehat {\mathcal E}_{n,\delta,\kappa,s}|_{C^0}=O(\delta/\kappa,\sqrt\kappa).
\]

    \medskip
    
    \noindent\textit{$C^2$ estimates:} Let $P_{\delta,\kappa}=Q_\delta R_\kappa$. From the chain rule it is not difficult to deduce that, for $i=x,y$ (we omit the subscripts $n,\delta,\kappa,s$ to alleviate the notation)
    \[
    D^2 \breve{\mathcal E}_i=\sum_{j=\varphi,J} P^{-1}_{j,i} P^T D^2\mathcal E_j P.
    \]
    Using \eqref{eq:errorestimatesIFSapprox C^2} we deduce that
    \[
    |P^{T} D^2\mathcal E_\varphi P|_{C^0}=O(N_\delta\delta \kappa)=O(\kappa)\qquad\qquad |P^{T} D^2\mathcal E_J P|_{C^0}=O(\delta\kappa).
    \]
    Hence, for every combination $i=x,y$ and $j=\varphi,J$
    \[
    |P^{-1}_{j,i} P^T D^2\mathcal E_j P|_{C^0}=O(\sqrt\kappa)
    \]
    and we deduce that $|D^2\widehat{\mathcal E}|_{C^0}=O(\sqrt\kappa)$. 
\medskip

 \noindent\textit{$C^3$ estimates:} These are obtained from \eqref{eq:errorestimatesIFSapprox C^3} in the very same way as we obtained the $C^2$ estimates. We leave the computation to the reader.

 \medskip

    \noindent\textit{Derivative with respect to $s$:} The estimates for $\partial_s \widehat{\mathcal E}$ are obtained from \eqref{eq:error derivative varepsilon},\eqref{eq:error simultaneous diagonal} and the chain rule in a similar fashion.
\end{proof}
\medskip

\subsection{Proof of Proposition \ref{prop:TT and HB}}
We now conclude the proof of Proposition \ref{prop:TT and HB} by showing the following.

\begin{prop}\label{prop: exists an epsilon good family}
    Let $I\subset (0,1)$ be a closed interval and let $\tilde s:(1/2,1)\to (0,1)$ be as in \eqref{eq:reparametrization}. There exists $C>0$ such that the following holds. For any $\epsilon>0$ there exist $n_1,n_2\in \mathbb N$  such that, if we write $\varepsilon_\epsilon(\lambda)=\epsilon^3 \tilde s(\lambda)$ the family of maps $\{F_{n_1,\varepsilon_\epsilon(\lambda)},F_{n_2,\varepsilon_\epsilon(\lambda)}\}_{\lambda\in J}$ is $C\epsilon$-good.
\end{prop}

\begin{proof}
  Let $\epsilon_0>0$ be given by Lemma \ref{lem: uniform diag plus scaling} and choose any $0<\epsilon\leq \epsilon_0\kappa$. Then, for any $n\in\mathcal N_{\epsilon^3,\kappa}$ any map $\widehat F_{n,\epsilon^3,\kappa,s}$ is $O_{C^3}(\sqrt\kappa, \epsilon/\kappa)$-close to its affine part, given by 
  \[
  \widetilde F_{\mathrm{aff},n,s}:\binom{x}{y}=\bs{\tilde b}_n+ \begin{pmatrix}
      \lambda(s)&0\\
      0&\lambda^{-1}(s)
  \end{pmatrix}\binom{x}{y}.
  \]
   We now choose suitable $n_1, n_2$ so the affine part coincides with those in the definition of $\epsilon$-good map. To do so we choose any  $n_1,n_2\in\mathcal N_{\epsilon^3,\kappa}$ for which 
  \[
  |[n_1\beta]|\leq C_0\kappa \qquad\qquad |[n_2\beta]-\sqrt\kappa|\leq C_0\kappa.
  \]
 On one hand, it follows from the definition of the matrix $P_{\epsilon^3,\kappa}$ that 
 \[
 \tilde{\bs b}_{n_1}=O(\sqrt\kappa)
 \]
 On the other hand, there exist $c_1,c_2\neq 0$ (independent of $s\in I$, $\kappa$ and $\epsilon$) such that
 \[
 \tilde {\bs b}_{n_2}=P_{\epsilon^3}^{-1}\binom{[n_2\beta]
 }{0}= \binom{c_1}{c_2}+O(\sqrt\kappa).
 \]
 We introduce one final change of coordinates $\psi_s$ associated to the diagonal matrix 
 \[
 M_s=\mathrm{diag}\bigg(\frac{c_1^{-1}}{1-\lambda(s)},\frac{c_2^{-1}}{1-\lambda^{-1}(s)}\bigg),
 \]
 and observe that, in the new coordinate system, the maps $\breve F_{n_i,\epsilon^3,\kappa,s}=\psi^{-1}_s\circ  \widehat F_{n_i,\epsilon^3,\kappa,s}\circ\psi$ satisfy
 \[
 |\breve F_{n_i,\epsilon^3,\kappa,s}-F_{\mathrm{aff},n_i,s}|_{C^3(\mathcal D)}=O(\epsilon^3/\kappa, \sqrt\kappa)\qquad\qquad |\partial_s(\breve F_{n_i,\epsilon^3,\kappa,s}-F_{\mathrm{aff},n_i,s})|_{C^2(\mathcal D)}=O(\epsilon^3/\kappa, \sqrt\kappa)
 \]
 Hence, if we choose $\kappa=\epsilon^{2}$ the conclusion follows with 
 \[
 \widetilde F_{i,\lambda}=\breve F_{n_i,\epsilon^3,\epsilon^2, s(\lambda)}.\qedhere
 \]
\end{proof}

\medskip

\section{From the IFS to the skew-product: proof of Proposition \ref{prop:Reduction Skew-prod to IFS}}\label{sec: from IFS to skewprod}

We now give the proof of Proposition \ref{prop:Reduction Skew-prod to IFS}. We first recall the standing assumptions. 
Consider the setting in Proposition \ref{prop:mainnhim} and let $\mathcal X_\varepsilon$ be the normally hyperbolic lamination constructed therein for the return map $\mathcal R_\varepsilon$. 
Recall that $\mathcal R_\varepsilon|_{\mathcal X_\varepsilon}:\mathcal X_\varepsilon\to \mathcal X_\varepsilon$ preserves any natural measure $\mathfrak m_\varepsilon$ on $\mathcal X_\varepsilon$ and that there exists a homeomorphism $\Psi_\varepsilon$ which conjugates this map to the skew-product
\[
\mathcal F_\varepsilon:(\bs\omega,z)\in \Omega\times N\to (\sigma \bs\omega, F_{\bs\omega,\varepsilon}(z))\in \Omega\times N
\]
in \eqref{eq:skewprodnhim}. An expression for $F_{\bs\omega,\varepsilon}$ in terms of the maps   $T_\varepsilon:N\to N$ and  $\{S_{i,\varepsilon}\}_{i=1}^m:N\to N$ defined in \eqref{eq:innerdynamicsdefn} and \eqref{eq:scattmapdynamicsdefn} is given in \eqref{eq:skewprodformula}. Recall that $C^2$ estimates for the remainder $R_{\omega,\varepsilon}$ in \eqref{eq:skewprodformula} are given in terms of (see Remark \ref{rem: barlambda})
\[
\bar\lambda=\alpha^{\kappa_0}\lambda<1
\]
for some $0\leq \kappa_0\leq \kappa$. 
\medskip

We want to prove that, provided the collection of maps $\{T_\varepsilon, \{S_{i,\varepsilon}\}_{i=1}^m\}$ satisfies \ref{it:itemB1}-\ref{it:B6}, there exists $\varepsilon_0>0$ such that, for a.e. $\varepsilon\in (0,\varepsilon_0)$, any natural invariant measure $(\Psi_\varepsilon)_*\mathfrak m_\varepsilon$ on $\Sigma\times N$ is ergodic for  $\mathcal F_\varepsilon$. 
\begin{rem}
    In order to avoid introducing extra notation we will simply denote by $\mathfrak m_\varepsilon$ the measure $(\Psi_\varepsilon)_*\mathfrak m_\varepsilon$ on $\Sigma\times N$. We will write 
    \[
    \mathrm{d} \mathfrak m_\varepsilon(\omega,z)=\mu_{\mathrm{Ber}}(\omega) \mathrm{vol}_{\omega,\varepsilon}(z).
    \]
\end{rem}

The idea of the proof is to reduce ergodicity of the skew-product to that of an IFS and then use the ideas in Section \ref{sec:IFSintro} to deduce ergodicity of the corresponding IFS. The proof is divided in a number of steps:
\begin{enumerate}
    \item We first define the holonomy maps between different fibers in the lamination and give certain quantitative estimates which will prove useful. 
    \item Then, given any measurable $\mathcal F_\varepsilon$-invariant set $B\subset \Omega \times N$  we produce an essentially bi-saturated representative $B_{\mathrm{good}}$.
    \item Making use of Fubini, and provided  $\mathfrak m_\varepsilon(B)>0$,  we choose a reference leaf $N_{\bs\omega^*}$ for which the set 
    \[
    E_*=B_{\mathrm{good}}\cap N_{\bs\omega^*}
    \]
    satisfies $\mathrm{vol}_{\bs\omega^*,\varepsilon}(E_*)> 0$. Moreover, exploiting the bi-saturation property of $B_{\mathrm{good}}$, we will be able to make use of a suitable collection of holonomy maps to construct a tuple of diffeomorphisms to $N_{\bs\omega^*}$ which are $\epsilon$-good (with some $\epsilon(\varepsilon)\to 0$ as $\varepsilon\to 0$). This implies local ergodicity of the corresponding iterated function system.

    \item Finally, we deduce global ergodicity by making use of  Proposition \ref{prop:IFS property R} and a shadowing argument.
\end{enumerate}

\subsection{The holonomy maps}
\begin{figure}
    \centering
\includegraphics[scale=0.45]
{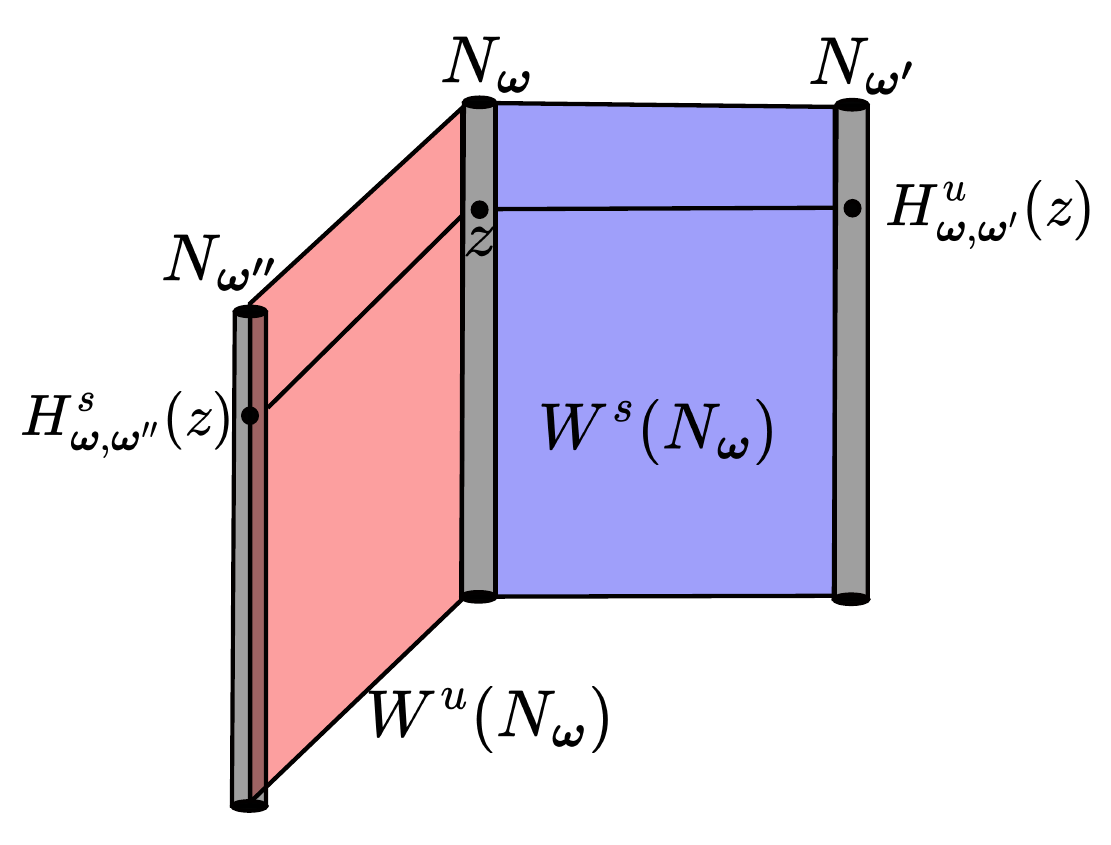}
    \caption{The holonomy maps between different leaves.}
    \label{fig:fig6}
\end{figure}

Given symbols $\bs\omega,\bs\omega'\in \Omega$ with $\bs\omega'\in W^u(\bs\omega)$ (resp. $\bs\omega'\in W^s(\bs\omega)$) we observe that at any $z\in N_{\bs\omega'}$ 
\[
 T_z N_{\bs\omega'} \oplus T_z W^{u}_{z^{u}}=T_z W^{u}(N_{\bs\omega})\qquad\qquad (\text{resp. }T_z N_{\bs\omega'} \oplus T_z W^{s}_{z^{s}}=T_z W^{s}(N_{\bs\omega}))
\]
where $z^{u}\in N_{\bs\omega}$ (resp. $z^{s}\in N_{\bs\omega}$)  is the unique point for which $z\in W^{u}_{z^u}$ (resp. $z\in W^{s}_{z^s}$) (see Figure \ref{fig:fig6}).
\medskip

It then follows from the discussion in Section \ref{sec:nihimintro}  that for $\bs\omega'\in W^u(\bs\omega)$ the  unstable wave map
\[
\Omega^u_{\bs\omega,\bs\omega',\varepsilon}:z\in N_{\bs\omega'}\mapsto z^u\in N_{\bs\omega}
\]
and for 
$\bs\omega'\in W^s(\bs\omega)$ the  stable wave map
\[
\Omega^s_{\bs\omega,\bs\omega',\varepsilon}:z\in N_{\bs\omega'}\mapsto z^s\in N_{\bs\omega}
\]
are $C^{\min\{r,\kappa\}}$ diffeomorphisms (recall that $\kappa$ is the bunching exponent in \ref{it:H1}). We then define the 
\textit{unstable holonomy map}
\begin{equation}\label{eq:unstableholonomy}
H_{\bs\omega,\bs\omega',\varepsilon}^u=(\Omega^u_{\bs\omega,\bs\omega',\varepsilon})^{-1}\in \mathrm{Diff}^{\min\{r,\kappa\}}(N_{\bs\omega},N_{\bs\omega'})
\end{equation}
and the 
\textit{stable holonomy map}
\begin{equation}\label{eq:stableholonomy}
H_{\bs\omega,\bs\omega',\varepsilon}^s=(\Omega^s_{\bs\omega,\bs\omega',\varepsilon})^{-1}\in \mathrm{Diff}^{\min\{r,\kappa\}}(N_{\bs\omega},N_{\bs\omega'}).
\end{equation}

\begin{lem}\label{lem: holonomies are AC}
    The holonomy maps \eqref{eq:unstableholonomy}, \eqref{eq:stableholonomy} are nonsingular. That is, for any ${\bs\omega}\in \Omega$ and ${\bs\omega}'\in W^{u,s}(\bs\omega)$ we have 
    \[
    (H^{u,s}_{{\bs\omega},{\bs\omega}',\varepsilon})_*\ \mathrm{vol}_{{\bs\omega},\varepsilon}\ll \mathrm{vol}_{{\bs\omega}',\varepsilon}.
    \]
\end{lem}

The following formula will prove useful to study the holonomy maps from a quantitative viewpoint.

\begin{lem}\label{lem:formulas holonomies}
  Let  $\kappa$ be the exponent in \ref{it:H1}-\ref{it:H1+1/2}. Then, there exists $\kappa_0'>0$ and  $\varepsilon_0>0$ for which, provided $\kappa\geq \kappa_0'$, then for any $0\leq \varepsilon\leq \varepsilon_0$ the  holonomy maps defined in \eqref{eq:unstableholonomy}, \eqref{eq:stableholonomy} satisfy
    \[
    H^u_{{\bs\omega},{\bs\omega}',\varepsilon}=\lim_{n\to \infty} (F_{{\bs\omega}',\varepsilon}^{-n})^{-1}\circ F_{{\bs\omega},\varepsilon}^{-n}\qquad\qquad H^s_{{\bs\omega},{\bs\omega}',\varepsilon}=\lim_{n\to \infty} (F_{{\bs\omega}',\varepsilon}^{n})^{-1}\circ F_{{\bs\omega},\varepsilon}^{n}
    \]
    in the $C^3$ topology.
\end{lem}

\begin{proof}
Throughout the proof $C$ denotes a positive constant which does not depend on $\varepsilon$, ${\bs\omega},{\bs\omega}'$ and $n$. We omit the dependence on $\varepsilon$ as it plays no role in the discussion. We first give the proof for the stable case and then indicate at the end the changes needed to treat the unstable case.  We start by checking the convergence of the above limits in the $C^2$ topology, and then show that they agree with the definitions \eqref{eq:unstableholonomy}, \eqref{eq:stableholonomy}. 
\medskip

\noindent\textit{Stable holonomy}: For the first part we proceed as follows. For each $n\in\mathbb N$ let 
\[
H^{s,n}_{{\bs\omega},{\bs\omega}'}=(F_{{\bs\omega}'}^{n})^{-1}\circ F_{\bs\omega}^{n}
    \]
and write \begin{equation}\label{eq:holonomy factoring} 
   H^{s,n+1}_{{\bs\omega},{\bs\omega}'}=(F_{{\bs\omega}'}^{n+1})^{-1}\circ F_{\bs\omega}^{n+1}=(F_{{\bs\omega}'}^{-n})^{-1}\circ \underbrace{F_{\sigma^n({\bs\omega}')}^{-1}\circ F_{\sigma^n({\bs\omega})}}_{d_n}\circ  \ F_{\bs\omega}^{n}=\underbrace{(F_{{\bs\omega}'}^{n})^{-1}\circ d_n\circ F_{{\bs\omega}'}^{n}}_{G_n}\circ\  H^{s,n}_{{\bs\omega},{\bs\omega}'}.
    \end{equation}
     Thus, after iteration we deduce that for any $n\in\mathbb N$
   \begin{equation}\label{eq:holonomy iteration} 
   H^{s,n}_{{\bs\omega},{\bs\omega}'}=G_{n-1}\circ\cdots \circ G_0.
    \end{equation}
    We now deduce suitable  estimates for $\lVert G_{j}-\mathrm{id}\rVert_{C^3}$ which guarantee the convergence of $H^{s,n}_{{\bs\omega},{\bs\omega}'}$ as $n\to \infty$. To do so we observe that since ${\bs\omega}'\in W^s({\bs\omega})$, for any $j\geq 0$ we have $\sigma^j({\bs\omega})_0=\sigma^j({\bs\omega}')_0$. Thus, it follows from the formula  \eqref{eq:skewprodformula} that 
    \[
    d_j({\bs\omega},{\bs\omega}')=F_{\sigma^j({\bs\omega}')}^{-1}\circ F_{\sigma^j({\bs\omega})}=(S_{\iota_j}+R_{\sigma^j({\bs\omega}')})^{-1}\circ (S_{\iota_j}+R_{\sigma^j({\bs\omega})}).
    \]
        It then follows from  \eqref{eq: comparing Romega ad inv manifolds}, that 
        
        \[
    |R_{{\bs\omega}'}-R_{{\bs\omega}}|_{C^3}\leq C  \bar\lambda^{\min\{{\omega}_{-1}',{\omega}_{-1}\}}
        \]
       and for $j\geq 1$
    \[
    |R_{\sigma^j({\bs\omega}')}-R_{\sigma^j({\bs\omega})}|_{C^3}\leq C  \bar\lambda^{\tau_{j}}. \qquad\qquad \tau_j=\min\{{\omega}_{-1}',{\omega}_{-1}\}+{\omega}_0+\dots+{\omega}_{j-1}.
    \]
    Hence, if we define $\tau_0=\min\{{\omega}_{-1}',{\omega}_{-1}\}$ it follows that for all $j\geq 0$
    \[
  |d_j-\mathrm{id}|_{C^3}\leq C \bar\lambda^{\tau_{j}}.
    \]
    To conclude we observe that (recall that $G_j$ was introduced in \eqref{eq:holonomy factoring}) 
    \begin{align*}
    |G_j-\mathrm{id}|_{C^0}\leq& C|F_{\bs\omega}'^{(n)}|_{C^1}|d_j-\mathrm{id}|_{C^0}\leq C\alpha^{\tau_j}||d_j-\mathrm{id}|_{C^0} \leq C(\alpha\bar\lambda)^{\tau_j} \\
    |G_j-\mathrm{id}|_{C^1}\leq& C |F_{{\bs\omega}'}^{(n)}|_{C^2}|(F_{{\bs\omega}'}^{(n)})^{-1}|_{C^2} |d_j-\mathrm{id}|_{C^1}\leq C \alpha^{2\tau_j} |d_j-\mathrm{id}|_{C^1}\leq C(\alpha^2\bar\lambda)^{\tau_j}
    \end{align*}
   and, similarly
   \[
   |G_j-\mathrm{id}|_{C^3}\leq  C (\alpha^4\bar\lambda)^{\tau_j}.
   \]
   Since $\bar\lambda=\alpha^{\kappa_0+1}\lambda$, by enlarging the exponent $\kappa$ in \ref{it:H1}-\ref{it:H1+1/2},  we can assume that 
   \[
   \alpha^4\bar\lambda=\alpha^{\kappa_0+5}\lambda<1.
   \]
   Then, the above estimates imply that \eqref{eq:holonomy iteration} converge in the $C^2$ topology to a diffeomorphism $K_{{\bs\omega},{\bs\omega}'}^s:N_{\bs\omega}\to N_{{\bs\omega}'}$. 
\medskip

It now remains to show that $K^s_{{\bs\omega},{\bs\omega}'}$ coincides with \eqref{eq:stableholonomy}. To do so observe that 
\[
K_{\sigma{\bs\omega},\sigma{\bs\omega}'}^s\circ F_{\bs\omega}
=
F_{{\bs\omega}'}\circ K_{{\bs\omega},{\bs\omega}'}^s.
\]
After iteration
\[
F_{{\bs\omega}'}^{
n}\circ K^s_{{\bs\omega},{\bs\omega}'}= K^s_{\sigma^n{\bs\omega},\sigma^n{\bs\omega}'}\circ F_{\bs\omega}^{n}.
\]
Now, because \(\sigma^n{\bs\omega},\sigma^n{\bs\omega}'\) have an increasingly long common past, \eqref{eq: comparing Romega ad inv manifolds} and the preceding estimates yield the improved bound
\begin{equation}\label{eq:K^s identity} 
\|K_{\sigma^n{\bs\omega},\sigma^n{\bs\omega}'}^s-\mathrm{id}\|_{C^0}
\leq C\bar\lambda^{\tau_n}.
\end{equation}
Indeed, the difference of the fiber maps at the shifted sequences starts at order $\bar\lambda^{\tau_n}$, and the successive conjugations contribute a summable geometric series with ratio controlled by \(\alpha\bar\lambda<1\). Using \eqref{eq:K^s identity} we obtain
\[
\lVert 
F_{{\bs\omega}'}^{n}(K_{{\bs\omega},{\bs\omega}'}^s(z))
-
F_{\bs\omega}^{n}(z)
\rVert_{C^0}
\leq C \bar\lambda^{\tau_n}.
\]
This says that the orbits of \(z\) and \(K_{{\bs\omega},{\bs\omega}'}^s(z)\) approach one another  at exponential rate. Hence, it is a classical fact (see \cite{MR501173}) that we must have 

\[
K^s_{{\bs\omega},{\bs\omega}'}(z)= W^s(z) \cap N_{{\bs\omega}'},
\]
that is, $K^s_{{\bs\omega},{\bs\omega}'}(z)=H^s_{{\bs\omega},{\bs\omega}'}(z)$.
\medskip

\noindent\textit{Unstable holonomy}: Finally, for the unstable case, we write
\begin{align*}
H^{u,n+1}_{{\bs\omega},{\bs\omega}'}=&(F_{{\bs\omega}'}^{-n-1})^{-1}\circ F_{\bs\omega}^{-n-1}=(F_{{\bs\omega}'}^{(-n)})^{-1}\circ \underbrace{F_{\sigma^{-n-1}({\bs\omega}')}\circ F_{\sigma^{-n-1}({\bs\omega})}^{-1}}_{\tilde d_n}\circ  \ F_{\bs\omega}^{(-n)}\\
=&\underbrace{(F_{{\bs\omega}'}^{(n)})^{-1}\circ \tilde d_n\circ F_{{\bs\omega}'}^{(-n)}}_{\tilde G_n}\circ\  H^{u,n}_{{\bs\omega},{\bs\omega}'}.
\end{align*}
Since for any $j\geq 0$ we have 
\[
\tilde d_j=F_{\sigma^{-j-1}({\bs\omega}')}\circ F_{\sigma^{-j-1}({\bs\omega})}^{-1}=T^{{\omega}_{-j-1}}\circ (S_{\iota_{j-1}}+R_{\sigma^{-j-1}({\bs\omega}')})\circ (S_{\iota_{j-1}}+R_{\sigma^{-j-1}({\bs\omega})})^{-1} \circ T^{-{\omega}_{-j-1}}
\]
and (also for any $j\geq 0$)
\[
|R_{\sigma^{-j-1}({\bs\omega}')}-R_{\sigma^{-j-1}({\bs\omega})}|_{C^3}\leq C \varepsilon  \bar\lambda^{\tilde\tau_j},\qquad\qquad \tilde \tau_j=\min\{{\omega}_0',{\omega}_0\}+{\omega}_{-1}+\dots+{\omega}_{-j-1}
\]
it is not difficult to check that 
\[
|\tilde d_j-\mathrm{id}|_{C^3}\leq C  \alpha^{4{\omega}_{-j-1}}\bar\lambda^{\tilde\tau_j}
\]
so 
\[
|\tilde G_j-\mathrm{id}|_{C^3}\leq C(\alpha^4\bar\lambda)^{\tilde \tau_j}
\]
and we are done. Indeed, $C^3$ convergence of $H^{u,n}_{{\bs\omega},{\bs\omega}'}$ follow directly and the equivalence of the limit definition and \eqref{eq:unstableholonomy} can be deduced by the very same argument used in the stable case.
\end{proof}

\subsubsection*{Quantitative estimates}
We now gather a few quantitative estimates obtained in the proof of Lemma \ref{lem:formulas holonomies} which will be of importance later. 

\begin{prop}\label{prop: close to id holonomies}
   Consider the setting of Lemma \ref{lem:formulas holonomies}. Then, there exists $\bar\lambda\in (0,1)$ (possibly larger than that in Proposition \ref{prop:mainnhim}), $\kappa_0''>0$ and $\varepsilon_0>0$ such that for any $\kappa\geq \kappa_0''$ and $0\leq \varepsilon \leq \varepsilon_0$ the following holds. For any ${\bs\omega}\in \Omega $ and any ${\bs\omega}'\in W^u({\bs\omega})$ (resp. ${\bs\omega}'\in W^s({\bs\omega})$) we have 
    \begin{equation}\label{eq:almostidentitys Holonomy}
        H^u_{{\bs\omega},{\bs\omega}',\varepsilon}(z)=z+O_{C^3}(\bar\lambda^{{\omega}_{-1}+\min\{{\omega}_0',{\omega}_0\}})\qquad\qquad (\text{resp. } H^s_{{\bs\omega},{\bs\omega}',\varepsilon}(z)=z+O_{C^3}(\bar\lambda^{\min\{{\omega}_{-1}',{\omega}_{-1}\}})).
    \end{equation}
\end{prop}
\begin{proof}
We give the proof of the stable and unstable case separately.
\medskip

\noindent\textit{Stable holonomy:} Let $\bar\lambda\in (0,1)$ be as in Proposition \ref{prop:mainnhim}. It follows from the estimates in the proof of Lemma \ref{lem:formulas holonomies} that (recall the expansion \eqref{eq:holonomy iteration})
    \[
    |H^{s,n}_{{\bs\omega},{\bs\omega}'}-\mathrm{id}|_{C^3}\leq C \sum_{j=0}^{n} (\alpha^4\bar\lambda)^{\tau_j}\leq  C(\alpha^4\bar\lambda)^{\mathrm{min}\{{\omega}_{-1},{\omega}_{-1}'\}}.
    \]
\noindent\textit{Unstable holonomy:}  Similarly, the estimates in the proof of Lemma \ref{lem:formulas holonomies} imply that  
  \[
    |H^{u,n}_{{\bs\omega},{\bs\omega}'}-\mathrm{id}|_{C^3}\leq C \sum_{j=0}^{n} (\alpha^4\bar\lambda)^{\tilde\tau_j}\leq  C (\alpha^4\bar\lambda)^{{\omega}_{-1}+\min\{\omega_0',\omega_0\}}.
    \]
    The result follows after relabeling $\bar\lambda$.
\end{proof}

The following result, which allows us to compare different holonomy maps, will also be useful later.

\begin{lem}
Consider the setting of Proposition \ref{prop: close to id holonomies}. For any ${\bs\omega}\in {\bs\omega} $ and any ${\bs\omega}',{\bs\omega}''\in W^u({\bs\omega})$ (resp. ${\bs\omega}',{\bs\omega}''\in W^s({\bs\omega})$) we have 
\begin{equation}\label{eq:comparing holonomies}
|H^u_{{\bs\omega},{\bs\omega}',\varepsilon}- H^u_{{\bs\omega},{\bs\omega}'',\varepsilon}|_{C^3}\leq C \bar\lambda^{{\omega}_{-1}+\min\{{\omega}_0',{\omega}_0''\}}   \qquad\qquad (\text{resp. } |H^s_{{\bs\omega},{\bs\omega}',\varepsilon}-H^s_{{\bs\omega},{\bs\omega}'',\varepsilon}|_{C^3} \leq C \bar\lambda^{\min\{{\omega}_{-1}',{\omega}_{-1}''\}}) 
\end{equation}
\end{lem}

\begin{proof}
     We borrow the terminology from the proof of Lemma \ref{lem:formulas holonomies} and  Proposition \ref{prop: close to id holonomies}. From their definition
\[
H_{{\bs\omega},{\bs\omega}'}^{u,s}
=
H_{{\bs\omega}'',{\bs\omega}'}^{u,s}
\circ H_{{\bs\omega},{\bs\omega}''}^{u,s}.
\]
Hence, 
\[
H^{u,s}_{{\bs\omega},{\bs\omega}'}-H^{u,s}_{{\bs\omega},{\bs\omega}''}=H_{{\bs\omega}'',{\bs\omega}'}^{u,s}
\circ H_{{\bs\omega},{\bs\omega}''}^{u,s}-H^{u,s}_{{\bs\omega},{\bs\omega}''}=(H^{u,s}_{{\bs\omega}'',{\bs\omega}'}-\mathrm{id})\circ H^{u,s}_{{\bs\omega},{\bs\omega}''}.
\]
But now observe that since ${\bs\omega}',{\bs\omega}''\in W^{u,s}({\bs\omega})$ the estimate \eqref{eq:almostidentitys Holonomy} implies the desired conclusion.
\end{proof}

\medskip

\subsection{A good representative}

Now, given a $\mathcal F_\varepsilon$-invariant measurable set $B\subset\Omega\times N$ we construct a representative with suitable properties. 
\begin{defn}
    Given a measurable set $B\subset\Omega\times N$ we say that it is \textit{s-saturated} (resp. \textit{u-saturated)} if for every $x\in B$ then $W^s(x)\subset B$ (resp. $W^u(x)\subset B$). We say it is \textit{essentially s-saturated} (resp. \textit{essentially u-saturated}) if there exists a measurable set  $B^s$ (resp. $B^u$) with $\frak m_\varepsilon(B\triangle B^s)=0$ (resp. $\frak m_\varepsilon(B\triangle B^u)=0$) which is s-saturated (resp. u-saturated). We say that it is \textit{essentially bi-saturated} if it is both $s$-saturated and $u$-saturated.
\end{defn}

\begin{lem}\label{lem:essential saturation} 
    Consider the setting of Proposition \ref{prop:mainnhim} and let $B\subset \Omega\times N$ be a $\mathcal F_\varepsilon$-invariant measurable set. Then, there exists a measurable set $A_0$ which is essentially bi-saturated and satisfies $\mathfrak m_\varepsilon(A_0\triangle A)=0$.
\end{lem}
\begin{proof}
The proof of this result follows plainly from the very same saturation step in the classical version of Hopf's  argument (see for instance 
Section 6.3. in  \cite{MR1963683}). Here we  only indicate the main idea for convenience of the reader. 

Let $B\subset \Omega\times N$ be a $\mathcal F_\varepsilon$-invariant measurable set. Choose a sequence of bounded Lipschitz functions $\varphi_j: \Omega\times N\to[0,1]$ such that
\[
\|\varphi_j-\mathbbm 1_G\|_{L^1(\mathfrak m)}\le 2^{-2j}.
\]
Fix some $j\in\mathbb N$. For every $x\in \Omega\times N$, define
\[
\varphi_j^\pm(x)
=
\limsup_{n\to\infty}
\frac1n\sum_{k=0}^{n-1}
\varphi_j(\mathcal F_\varepsilon^{\pm k}(x)).
\]
We claim that $\varphi_j^+$ is constant along strong-stable leaves. Indeed, if $y\in W^s(x)$, then
\[
d(\mathcal F^n_\varepsilon(x),\mathcal F^n_\varepsilon(y))\to 0.
\]
Since $\varphi_j$ is uniformly continuous,
\[
\frac1n\sum_{k=0}^{n-1}
\left|\varphi_j(\mathcal F^k_\varepsilon(x))
-\varphi_j(\mathcal F^k_\varepsilon(y))\right|
\to 0.
\]
Consequently,
\[
\varphi_j^+(x)=\varphi_j^+(y).
\]
Analogously, \(\varphi_j^-\) is constant along strong-unstable leaves.
\medskip

On the other hand, Birkhoff's ergodic theorem gives
\[
\varphi_j^\pm
=
\mathbb E(\varphi_j\mid\mathcal I)
\qquad \mathfrak m_\varepsilon\text{-a.e.},
\]
where $\mathcal I$ is the invariant sigma-algebra of $\mathcal F_\varepsilon$. Since $B$ is invariant, we have
\[
\mathbb E(\mathbbm 1_B|\mathcal I)=\mathbbm 1_B.
\]
Thus, by the $L^1$-contraction property of conditional expectation,
\[
\|\varphi_j^\pm-\mathbbm 1_B\|_{L^1(\mathfrak m)}
\leq \|\varphi_j-\mathbbm 1_B\|_{L^1(\mathfrak m)}
\leq 2^{-2j}.
\]
It follows that
\[
\varphi_j^\pm(x)\longrightarrow \mathbbm 1_B(x)
\qquad\text{for }\mathfrak m_\varepsilon\text{-a.e. }x.
\]
Define
\[
A^s=\left\{x:\limsup_{j\to\infty}\varphi_j^+(x)>\frac12\right\},
\qquad
A^u=\left\{x:\limsup_{j\to\infty}\varphi_j^-(x)>\frac12\right\}.
\]
By construction, \(A^s\) is  \(s\)-saturated and \(A^u\) is  \(u\)-saturated. Moreover,
\[
\mathfrak m_\varepsilon (A^s\triangle A)=
\mathfrak m_\varepsilon(A^u\triangle A)=0.
\]
Taking $A_0=A^u$, we obtain an essentially bi-saturated measurable representative of $A$, as required. \end{proof}
    
    Now, given any cylinder $a=[a_0,\dots,a_n]\in \Omega$ we define the \textit{forward insertion map} 
\[
    {\bs\omega}=(\dots {\bs\omega}_{-1};{\bs\omega}_0,\dots)\mapsto \iota_a^+({\bs\omega})=(\dots {\bs\omega}_{-1};a,{\bs\omega}_0,\dots),
\]
the \textit{backward insertion map}
\[
{\bs\omega}=(\dots {\bs\omega}_{-1};{\bs\omega}_0,\dots)\mapsto \iota^-_a({\bs\omega})=(\dots {\bs\omega}_{-1},a;{\bs\omega}_0,\dots),
\]
and denote by 
\begin{equation}\label{eq:insertion map Skewprod}
\mathcal I_{a,\varepsilon}^+({\bs\omega},z)=(\iota_a^+({\bs\omega}), H^u_{{\bs\omega},\iota^+_a({\bs\omega}),\varepsilon}(z))\qquad\qquad \mathcal I_{a,\varepsilon}^-({\bs\omega},z)=(\iota_a^-({\bs\omega}), H^s_{{\bs\omega},\iota_a^-({\bs\omega}),\varepsilon}(z)).
\end{equation}

\begin{lem}\label{lem:good representative}
     Consider the setting of Proposition \ref{prop:mainnhim} and let $B\subset \Omega\times N$ be a $\mathcal F_\varepsilon$-invariant measurable set. Let  
     \[
     \mathcal W=\{({\bs\omega}_0,\dots,{\bs\omega}_n),\colon n\in \mathbb N\}.
     \] 
     be the set of all finite words. There exists a measurable set $B_{\mathrm{good}}$ which satisfies the following:
    \begin{enumerate}
        \item \textit{(Equivalence to $B$)}:   $\mathfrak m_\varepsilon(B_{\mathrm{good}}\triangle B)=0$,
        \item \textit{(Invariance under the insertion map)}: For any $a\in \mathcal W$ and any $\star=\pm$
        \[
        ({\bs\omega},z)\in B_{\mathrm{good}}\Longrightarrow \mathcal I_{a,\varepsilon}^\star({\bs\omega},z)\in B_{\mathrm{good}}.
        \]
        \item \textit{($\mathcal F$-invariance)}: $({\bs\omega},z)\in B_{\mathrm{good}}\Longleftrightarrow (\sigma({\bs\omega}), F_{{\bs\omega},\varepsilon}(z))\in B_{\mathrm{good}}$
    \end{enumerate}
\end{lem}
\begin{proof}
We only give the proof in the case where $a=[a_0]$ (i.e. its length is one) as the general case follows from the very same argument and the fact that $\mathcal W$ is countable. We omit the dependence on $\varepsilon$. The proof is divided in a number of steps.
\medskip

\noindent\textit{Non-singularity of the insertion map:}   We start by observing that the transformation \eqref{eq:insertion map Skewprod} is non-singular, that is, for any measurable set $B$ with $\mathfrak m(B)=0$ we have $\mathfrak m(\mathcal I_a^{-1}(B))=0$.  To see this we proceed as follows. Given any ${\bs\omega}_*\in \Omega$ we introduce the slice  $B_{\bs\omega}=B\cap\{{\bs\omega}={\bs\omega}_*\}$. Clearly, we must have $\mathrm{vol}_{\bs\omega}(B_{\bs\omega})=0$ for $\nu$-a.e. ${\bs\omega}$.

First, we note that the symbolic transformation $\iota_a:\Omega\to\Omega$ is non-singular. Indeed, for any measurable symbolic set $C\subset \mathcal B(\Omega)$ (one may check this using cylinders)
   \[
   \mu_{\mathrm{Ber}}(\iota_a^{-1}(C))=\frac{1}{p_a}\mu_{\mathrm{Ber}}(C\cap \{{\bs\omega}_0=a\}).
   \]
In particular, for $\mu_{\mathrm{Ber}}$-a.e. ${\bs\omega}$
\[
\mathrm{vol}_{\iota_a({\bs\omega})}(B_{\iota_a}({\bs\omega}))=0.
\]
Second, we have shown in Lemma \ref{lem: holonomies are AC} that  $(H^u_{{\bs\omega},\iota_a({\bs\omega})})_*\mathrm{vol}_{\bs\omega}\ll \mathrm{vol}_{\iota_a({\bs\omega})}$ for any ${\bs\omega}$. Thus, it follows that for  $\mu_{\mathrm{Ber}}$-a.e. ${\bs\omega}$ we must have 
   \[
   \mathrm{vol}_{\bs\omega}( (H^u_{{\bs\omega},\iota_a({\bs\omega})})^{-1}(B_{\iota_a({\bs\omega})}))=0.
   \]
   Therefore, nonsingularity of $\mathcal I_a$ follows since 
   \[
   \mathfrak m(\mathcal I_a^{-1}(B))=\int_\Omega \mathrm{vol}_{\bs\omega}( (H^u_{{\bs\omega},\iota_a({\bs\omega})})^{-1}(B_{\iota_a({\bs\omega})}))\ \mathrm d\mu_{\mathrm{Ber}}({\bs\omega})=0.
   \]
   \medskip

   \noindent\textit{Saturation and insertion invariance:} Consider the $u$-saturated set  $B^u$ constructed in  the proof of  Lemma \ref{lem:essential saturation}. For any $z\in B^u$ we have $\mathcal I_a(z)\in W^u(z)$ so $\mathcal I_a(z)\in B^u$. Therefore, we have $B^u\subset \mathcal I_a^{-1}(B^u)$. Moreover, for any $z\in \mathcal I_a^{-1}(B^u)$ we have $\mathcal I_a(z)\in B^u$, but since $z\in W^u(\mathcal I_a(z))$ it holds that $z\in B^u$. Hence, we must  have $\mathcal I_a^{-1}(B^u)\subset B^u$ and we conclude that 
   \[
   B^u=\mathcal I_a^{-1}(B^u).
   \]
   On the other hand, since $\mathcal I_a$ is non-singular
   \[
   \mathfrak m(\mathcal I_b^{-1}(B^u\triangle B))=0
   \]
   Thus, 
   \[
   \mathcal I_b^{-1}(B)=\mathcal I_b^{-1}(B^u)=B^u=B\qquad\qquad\mathrm{mod}\ \mathfrak m.
   \]
   \medskip

   \noindent\textit{Conclusion:} We denote by $G=\langle \mathcal F,\mathcal F^{-1},\{\mathcal I_b^+\}_{b}, \{\mathcal I_b^-\}_b\rangle$ and define the set
   \[
   B_{\mathrm{good}}=\bigcap_{g\in G} g^{-1}(B)
   \]
   which, clearly, satisfies items (1)-(3) in the statement.
\end{proof}
\medskip

\subsection{Ergodicity for a suitable IFS on a reference leaf}
\begin{figure}
    \centering
\includegraphics[scale=0.45]
{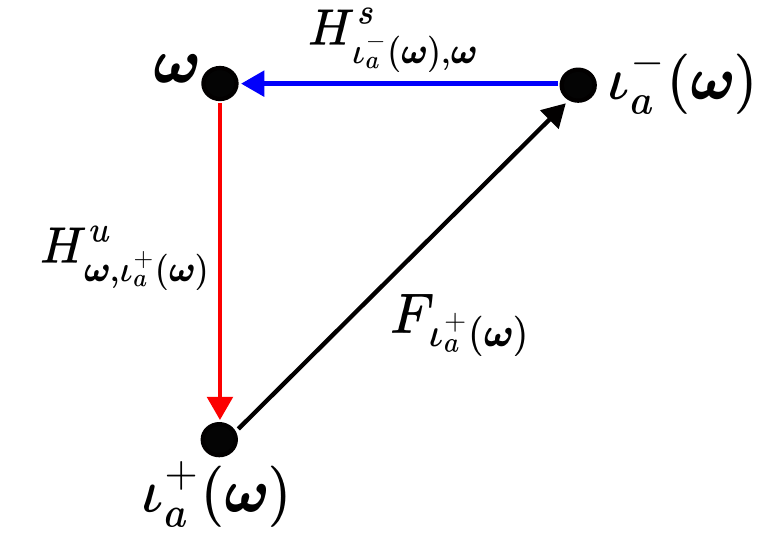}
    \caption{Definition of the map $H_a:N_{\bs\omega}\to N_{\bs\omega}$.}
    \label{fig:fig7}
\end{figure}

Given any ${\bs\omega}\in \Sigma$ and the collection of finite words $\mathcal W$ we now consider the family of diffeomorphisms (see Figure \ref{fig:fig7})
\begin{equation}\label{eq:F loop maps} 
\{H_{a,\varepsilon}\}_{a\in \mathcal W}:N_{\bs\omega}\to N_{{\bs\omega}},\qquad\qquad H_{a,\varepsilon}=H^s_{\iota^-_a({\bs\omega}),{\bs\omega},\varepsilon}\circ F_{\iota_a^+({\bs\omega})}^{|a|}\circ H^u_{{\bs\omega},\iota_a^+({\bs\omega}),\varepsilon}.
\end{equation}
To conclude the proof of Proposition \ref{prop:Reduction Skew-prod to IFS} we just need to show that given any set measurable set $E\subset N_{\bs\omega}$
\[
\mathrm{vol}_{{\bs\omega},\varepsilon}(E)>0 \quad\text{and}\quad E\text{ is invariant under all the maps }\{H_{a,\varepsilon}\}_{a\in \mathcal W}\quad\Longrightarrow\quad \mathrm{vol}_{{\bs\omega},\varepsilon}(E)=1.
\]
We will do so in two steps. 

\subsubsection*{A $\epsilon$-good family} We begin by selecting symbols $a_1,a_2$ for which the corresponding maps \eqref{eq:F loop maps} form a $\epsilon$-good family.
\begin{lem}\label{lem: epsilon good skewprod}
   Let $J\subset (1/2,1)$ be a closed interval, let $\tilde s:(1/2,1)\to (0,1/2)$ be as in \eqref{eq:reparametrization}. There exists $C>0$ such that the following holds. For any $\epsilon>0$ there exist $a_1,a_2\in \mathcal W$  such that the following holds for any 
   \begin{equation}\label{eq:Omegaepsset}
   {\bs\omega}\in \Omega_\epsilon=\{{\bs\omega}\in\Omega\colon \omega_0,{\omega}_{-1}\geq \epsilon^{-1}\}.
   \end{equation}
   If we write $\varepsilon_\epsilon(\lambda)=\epsilon^3 \tilde s(\lambda)$ the family of maps  
   \[
   \{H_{a_1,\varepsilon_\epsilon(\lambda)},H_{a_2,\varepsilon_\epsilon(\lambda)}\}_{\lambda\in J}\subset \mathrm{Diff}^3(N_{\bs\omega})
   \]
   is $C\epsilon$-good.
\end{lem}

\begin{proof}
Fix any $\epsilon>0$. On one hand, we have shown in Proposition \ref{prop: exists an epsilon good family} that there exist $n_1,n_2\in\mathcal N_\epsilon$ such that the maps $\{F_{n_1,\varepsilon_\epsilon(\lambda)},F_{n_2,\varepsilon_\epsilon(\lambda)}\}_{\lambda\in J}\subset \mathrm{Diff}^3(N)$ form a $C\epsilon$-good family. Set 
\[
a_1=(n_1,1)\qquad\qquad a_2=(n_2,1)
\]
On the other hand, for both $i=1,2$ the maps $H^s_{\iota^-_{a_i}({\bs\omega}),{\bs\omega}}$, and $H^u_{{\bs\omega},\iota_{a_i}^+({\bs\omega})}$ in $\eqref{eq:F loop maps}$ the estimates in Proposition \ref{prop: close to id holonomies} imply that (here $\bar\lambda\in (0,1)$ is as in Proposition \ref{prop: close to id holonomies} and has nothing to do with $\lambda\in J$ parametrizing the family of maps $\{F_{n_1,\varepsilon_\epsilon(\lambda)},F_{n_2,\varepsilon_\epsilon(\lambda)}\}_{\lambda\in J}$)
\[
H^u_{{\bs\omega},\iota_{a_i}^+({\bs\omega}),\varepsilon}=z+O_{C^3}(\bar\lambda^{{\omega}_{-1}+\min\{\omega_0,a\}}) \qquad\qquad H^s_{\iota^-_{a_i}({\omega}),{\omega},\varepsilon}(z)=z+O_{C^3}( \bar\lambda^{\min\{{\omega}_{-1},n_i\}}).
\]
Since ${\omega}_{-1}\geq \epsilon^{-1}$ and also $n_1,n_2\geq \epsilon^{-1}$ we deduce that indeed (uniformly for all $\lambda \in J$)
\[
H^u_{{\bs\omega},\iota_{a_i}^+({\bs\omega}),\varepsilon_\epsilon(\lambda)}=z+O_{C^3}( \bar\lambda^{C/\epsilon}) \qquad\qquad H^s_{\iota^-_{a_i}({\bs\omega}),{\bs\omega},\varepsilon_\epsilon(\lambda)}(z)=z+O_{C^3}(\bar\lambda^{C/\epsilon}).
\]
Hence, the corresponding $\{H_{a_1,\varepsilon_\epsilon(\lambda)},H_{a_2,\varepsilon_\epsilon(\lambda)}\}_{\lambda\in J}$ also form a $C\epsilon$-good family (possibly after increasing $C$).
\end{proof}
\medskip

\subsubsection*{Global transport}

We now make use of the ideas behind the proof of Proposition \ref{prop:IFS property R} to obtain the following.

\begin{lem}\label{lem:transportskewprod}
Let $N$ be a smooth compact Riemannian surface and let $\rho$ denote the associated area form. Suppose that the diffeomorphism $T_\varepsilon:N\to N$ in \eqref{eq:innerdynamicsdefn} preserves a smooth measure and that the system $\{S_{i,\varepsilon}\}_{i=1}^m:N\to N$ in  \eqref{eq:scattmapdynamicsdefn} also preserves a (possibly different) smooth measure and satisfies \ref{it:B6} for some $n\geq 1$. Choose any ${\bs\omega}\in \Omega_\epsilon$ (defined in \eqref{eq:Omegaepsset}) and given some $a\in \mathcal W$ denote by $H_{a,\varepsilon}:N_{\bs\omega}\to N_{\bs\omega}$ the diffeomorphism in \eqref{eq:F loop maps}. 

Then, there exists $K>0$  and $\varepsilon_0>0$ such that for any  $0<\varepsilon\leq \varepsilon_0$ the following holds. There exists a  subset $R_\varepsilon \subset N_{\bs\omega}$ with $\rho(N\setminus R_\varepsilon)=0$ such that given any pair of  points $z\in R_\varepsilon$,  $z_*\in N$  there exists $\bs{\tilde{\bs\omega}}\in \mathcal W$ such that 
\[
\mathrm{dist}(H_{\bs{\tilde\omega},\varepsilon}(z),z_*)\leq K \varepsilon^{n+1}.
\]
\end{lem}

The proof of Lemma \ref{lem:transportskewprod} is somewhat technical and is deferred to Appendix \ref{sec: shadowing}. Roughly speaking, the idea is that, by exploiting Poincaré recurrence, one can shadow obits of the IFS generated by $\{S_{i,\varepsilon}\}_i$ within the skew-product dynamics.

\subsection{Proof of Proposition \ref{prop:Reduction Skew-prod to IFS}}
Given some $\varepsilon_0>0$ we fix some $\epsilon\in (0,\varepsilon_0/2)$ and restrict our attention to $\varepsilon\in I_\epsilon=(\epsilon/2,2\epsilon)$. The full result follows by exhausting the interval $(0,\varepsilon_0)$ by intervals of the form $I_\epsilon$ for a sequence of values of $\epsilon$.
\medskip

We first fix any $\varepsilon\in I_\epsilon$. Let $B\subset\Omega^\mathbb Z\times N$ be a measurable set which is $\mathcal F_\varepsilon$-invariant and has positive measure $\mathfrak m_\varepsilon(B)>0$. Let $B_{\mathrm{good}}\subset \Omega^\mathbb Z\times N$ be the measurable representative constructed in Lemma \ref{lem:good representative}. Then, since the function $p_{B}:\Omega\mapsto [0,1]$ given by
\[
p_{B}(\bs\omega)=\mathrm{vol}_{\bs\omega,\varepsilon}(B_{\mathrm{good}}\cap N_{\bs\omega})
\]
is $\sigma$-invariant and the shift is ergodic it must be constant almost everywhere. Hence, from Fubini it follows that  exists a set $\Omega_*\subset \Omega$ with $\mu_{\mathrm{Ber}}(\Omega_*)=1$  such that for any ${\bs\omega}_*\in \Omega_*$
\begin{equation}\label{eq: E star set}
E_{*}=B_{\mathrm{good}}\cap N_{{\bs\omega}^*}
\end{equation}
 satisfies
\[
p_{B}(\bs\omega^*)=\mathrm{vol}_{{\bs\omega}^*,\varepsilon}(E_{*})>0.
\]
On the other hand, if we now fix any ${\bs\omega}\in \Omega_\epsilon$ (and let free $\varepsilon\in I_\epsilon$), consider the countable collection of diffeomorphisms $\{H_{a,\varepsilon}\}_{a\in\mathcal W}\subset  \mathrm{Diff}^3 (N_{{\bs\omega}})$ in \eqref{eq:F loop maps} and denote by 
\[
\mathcal G_{{\bs\omega},\varepsilon}=\langle \{H_{a,\varepsilon}\}_{a\in \mathcal W}\rangle \subset \mathrm{Diff}^3(N_{{\bs\omega}})
\]
it follows from Lemmas \ref{lem: epsilon good skewprod}-\ref{lem:transportskewprod} and  Proposition \ref{prop:global ergodicity} that the group action $(\mathcal G_{{\bs\omega},\varepsilon},N_{\bs\omega},\mathfrak m_\varepsilon)$ is ergodic for almost every $\varepsilon\in I_\epsilon$.
\medskip

It then follows from Fubini's theorem that for almost every $\varepsilon\in I_\epsilon$ there exists ${\bs\omega}^*\in \Omega_\epsilon$ such that:
\begin{itemize}
\item the set $E_*$ defined in \eqref{eq: E star set} satisfies $\mathrm{vol}_{{\bs\omega}^*,\varepsilon}(E_{*})>0$ and,
    \item the group action $(\mathcal G_{\varepsilon,{\bs\omega}^*},N_{{\bs\omega}^*},\mathrm{vol}_{\bs\omega,\varepsilon})$ is ergodic.
\end{itemize}
But then, $\mathrm{vol}_{{\bs\omega}^*,\varepsilon}(E_{*})=1$. Hence,  we must have $\mathfrak m_\varepsilon(B)=1$. This completes the proof of Proposition \ref{prop:Reduction Skew-prod to IFS}.

\section{Technical discussion of possible extensions}\label{sec: extensions section}
Here we discuss more in detail the possible extensions described in Section \ref{sec:futurework}.

\subsection{Higher dimensional center}
In order to treat the case where the dimension of $N$ is arbitrary one needs a higher dimensional generalization of the concept of Hopf box. There are two difficulties:
\begin{itemize}
    \item In the higher dimensional case $\mathrm{dim}(N)=2d$, one has $d$ stable directions and $d$ unstable directions. However, for $d\geq 2$, a pair formed by a couple of parallel stable manifolds $\gamma_d,\gamma_u$ and a couple of parallel unstable manifolds $\gamma_r,\gamma_l$ which intersect transverally pairwise ($\gamma_d\pitchfork \gamma_r\neq 0$ and $\gamma_u\pitchfork \gamma_l\neq 0$) does not bound an open region. The way to go around this obstacle is to instead define the Hopf box as the intersection $\mathcal C_\lambda=\mathcal U_\lambda\cap \mathcal S_\lambda$ with 
    \[
    \mathcal U_\lambda=\bigcup_{\omega^-} \gamma^u_{\omega^-}\qquad\qquad \mathcal S_\lambda=\bigcup_{\omega^+} \gamma^s_{\omega^+}.
    \]
    \item For the case $d=1$, in order to have good bi-covering properties using only two maps we asked the contraction rate $\lambda$ to belong to the interval $(1/2,1)$. We showed in Section \ref{sec:proof TT } that for maps satisfying \ref{it:itemB1}-\ref{it:itemB5} we can always choose a hyperbolic word with this contraction rate. 

    For $d\geq 2$ this is not possible for a generic system (see Section 4 in \cite{paradela2026robustlytransitivebehaviorsymplectic}) so we must consider different contraction rates along the different center hyperbolic directions. However, this is not an issue as it is still possible to choose a large collection of words in order to have a bi-covering property. 
\end{itemize}

\subsection{Property $K$} The current approach in Section \ref{sec: from IFS to skewprod} to establish ergodicity exploits $\mathcal F_\varepsilon$ invariance as well as invariance under un/stable holonomies (recall, in particular, the definition of the maps \eqref{eq:F loop maps}). 

To establish property $K$ we would instead need to show that any set in the Pinsker sigma algebra has measure zero or one. However such sets are not $\mathcal F_\varepsilon$ invariant. Hence, we should replace the approach in Section \ref{sec: from IFS to skewprod} and only use  $su$-loops  instead. To do so one needs to show that for a typical system one can find a collection of $su$-loops verifying hypothesis \ref{it:itemB1}-\ref{it:B6}. This seems reasonable making use of the perturbation theory developed in \cite{DelaLLaveScattmap}.

\subsection{Beyond the direct product setting} The main reason why we presented Theorem \ref{thm:mainham} in the direct product setting is that, in this case the scattering maps to the normally hyperbolic manifold $\Lambda$ are close to identity. Hence,  making use of the tools in \cite{DelaLLaveScattmap} it is rather easy to   verify that assumption \ref{it:B6} holds for a typical perturbation.

In the non-direct product setting in which we just assume that the Hamiltonian diffeomorphism $F_0$ has a compact normally hyperbolic manifold $\Lambda$ with a collection of full homoclinic channels, the scattering maps to $\Lambda$ need not to be close to the identity.  This does not seem to be an obstacle to establish ergodicity of the dynamics on the corresponding lamination as one could instead obtain a collection of close to identity maps satisfying \ref{it:B6} by looking at a suitable collection of $su$-loops.

\appendix

\medskip

\section{The real-analytic topology in the space of diffeomorphisms}\label{sec:realanalytictop}

A differentiable  manifold $X$ is said to be \textit{real-analytic} if given any two charts 
\[
\varphi_i: U_i\subset \mathbb R^d\to V_i\subset X,\qquad\qquad \varphi_j: U_j\subset \mathbb R^d\to V_j\subset X
\]
with $V_i\cap V_j\neq \emptyset$ then, the corresponding transition function
\[
\varphi_{i,j}:=\varphi_j^{-1}\circ \varphi_i:U_{i,j}\subset \mathbb R^d \to U_{j,i}\subset\mathbb R^d\qquad\qquad \text{where}\qquad U_{i,j}:=\varphi_i^{-1}(V_i\cap V_j)
\]
is real-analytic.  Provided $X$ is compact there exists a complexification $X_{\mathbb C}$, i.e. a complex manifold $X_\mathbb C$ containing $X$ as a totally real submanifold and, moreover, this complexification is unique up to biholomorphisms\footnote{If $X_{\mathbb C}$ and $X_{\mathbb C}$ are two such complexifications then there exists neighbourhoods $U\subset X_{\mathbb C}$ and $U'\subset X_{\mathbb C}'$ of $X$ in $X_{\mathbb C}$ and in $X_{\mathbb C}'$, and a biholomorphism $\varphi:U\to U'$.} (see Proposition 1 in \cite{MR102094}). 
\subsection{The space $\mathcal C^\omega(X)$}
Given a Riemannian compact manifold $X$ let 
\[
\mathcal C^\omega(X)=\{f:X\to \mathbb R\colon f \text{ is real-analytic}\}
\]
Then, the real-analytic topology on $\mathcal C^\omega(X)$ is the one generated by the open sets 
\begin{align*}
\mathcal U(f,K,\varepsilon)=&\{g\in \mathcal C^\omega(X)\colon g\text{ admits a holomorphic extension $\bar g$ to an open neighbourhood of }K\text{ and}\\
&\qquad \sup_{z\in K}|\bar g(z)-\bar f(z)|<\varepsilon\}
\end{align*}
where $X\subset K\subset X_\mathbb C$ is compact and $f\in \mathcal C^\omega (X)$ admits a holomorphic extension $\bar f$  to an open neighbourhood of $K$. Namely, two elements are close if they are close in some compact neighbourhood of $X$ in $X_{\mathbb C}$.
\subsection{The space $\mathrm{Diff}^\omega(X)$}
Given a Riemannian compact manifold $X$, the real-analytic topology on $\mathrm{Diff}^\omega(X)$ is the one generated by the open sets 
\begin{align*}
\mathcal U(F,K,\varepsilon)=&\{G\in \mathrm{Diff}^\omega(X)\colon G\text{ admits a holomorphic extension $\bar G$ to an open neighbourhood of }K\text{ and}\\
&\qquad \sup_{z\in K}\mathrm{dist}(\bar G(z),\bar F(z))<\varepsilon\}
\end{align*}
where $X\subset K\subset X_\mathbb C$ is compact and $F\in \mathrm{Diff}^\omega(X)$ admits a holomorphic extension $\bar F$  to an open neighbourhood of $K$. Namely, two elements are close if they are close in some compact neighbourhood of $X$ in $X_{\mathbb C}$.
\medskip

\section{Sobolev regularity of the slice measures}\label{sec:sobolev appendix}
In this appendix we give a proof of Theorem \ref{thm:sobolev}. In the following we fixe any value of $y\in [0,1]$ and note that all the estimates are uniform in $y$. We omit this variable from the notation. 
\medskip

We recall that for a typical $\lambda$ the density $a_\lambda(\cdot,y)=\mathrm d\eta^-_\lambda/\mathrm dx$ belongs to $L^2([0,1])$. We denote by $\tilde a_\lambda\in L^2(\mathbb R)$ its extension
\[
\tilde a_\lambda(x)=\begin{cases}
    a_\lambda(x,y) &\text{ if }x\in [0,1]\\
    0&\text{ otherwise}\\
\end{cases}
\]
and define
\[
[A_\lambda](u)=\int_{[0,1]}\tilde a_\lambda (x,y)\tilde a_\lambda(x+u,y)\mathrm dx.
\]
The first observation towards an estimate for the Sobolev norm of $\tilde a_\lambda$ is that 
\[
\lVert \tilde a_\lambda(\cdot +u,y)-\tilde a_\lambda (\cdot,y)\rVert^2_{L^2}=2([A_\lambda](0)-[A_\lambda](u)).
\]
\begin{prop}\label{prop:estimate Sobolev technical appendix}
  Choose any $\lambda_0 \in (1/2,2^{-2/3})$ and any

    \[
    0<\sigma<\frac{\log 2}{\log(1/\lambda_0)} -1
    \]
    Then, there exists $C>0$ such that  for any sufficiently small neighbourhood $U_0$ of $\lambda_0$    and any sufficiently small $\varepsilon$ we have
\begin{equation}\label{eq: estimate Sobolev technical appendix}
        \int_{U_0} ([A_\lambda](0)-[A_\lambda](u))\ \mathrm d\lambda \leq C |u|^\sigma.
    \end{equation}
\end{prop}
As we will see below, Theorem \ref{thm:sobolev} follows plainly from Proposition \ref{prop:estimate Sobolev technical appendix}. We devote the rest of the section to first establish \eqref{eq: estimate Sobolev technical appendix} and then  conclude the proof of Theorem \ref{thm:sobolev}. 

\medskip To obtain an estimate for the expression above we relate $[A_\lambda]$ to the symbolic coding. Making use of the definition of $\eta^-_\lambda$ in \eqref{eq: final projection measure} we deduce that for any $u\in \mathbb R$ we have 
\begin{equation}\label{eq:Symbolic coding difference Alambda}
[A_\lambda](0)-[A_\lambda](u)
=\lim_{r\to0}\frac{1}{2r}
\int_{\{1,2\}^\mathbb N}\int_{\{1,2\}^\mathbb N}
\Big(\mathbbm{1}_{\{|\Delta_\lambda(\omega,\eta)|<r\}}-\mathbbm{1}_{\{|\Delta_\lambda(\omega,\eta)-u|<r\}}
\Big)
\,\mathrm d\mu^-(\omega)\mathrm d\mu^-(\eta).
\end{equation}
Following \cite{MR1749437,MR4385129} we now deploy techniques similar to those in Section \ref{sec:typicalAC} to obtain an estimate for this integral. 
For any $r>0$ we write 
\[
K_r=\frac{1}{2r}\mathbbm 1_{[-r,r]}
\]
Then, given any $u\in \mathbb R$ we let 
\[
S_r(\lambda,u)=\int_{\{1,2\}^\mathbb N}\int_{\{1,2\}^\mathbb N} K_r(\Delta_\lambda (\omega,\eta)-u)\  \mathrm d\mu^-(\omega)\mathrm d\mu^-(\eta)
\]
be the integral in \eqref{eq:Symbolic coding difference Alambda}. It follows from the estimates in Section \ref{sec:typicalAC} that this integral is absolutely convergent for any $r>0$. To complete the proof of Proposition \ref{prop:estimate Sobolev technical appendix} we would like to estimate
\[
\mathcal S_r(u)=\int_{U_0} (S_r(\lambda,0)-S_r(\lambda,u))\ \mathrm d\lambda
\]
uniformly in $r$. To do so we let $U_1$ be an interval containing $U_0$ in its interior and choose a function $\chi\in C^\infty_c(U_1)$ satisfying $\chi=1$ on $U_0$. Then we define 
\[
S_r(u;\omega,\eta)=\int_{U_1}\chi(\lambda) K_r(\Delta_\lambda (\omega,\eta)-u)\  \mathrm d\lambda 
\]
and note that by Fubini we have 
\[
\mathcal S_r(u)\leq \int_{\{1,2\}^\mathbb N}\int_{\{1,2\}^\mathbb N} (\underbrace{S_r(0;\omega,\eta)-S_r(u;\omega,\eta)}_{\Delta S_r(u,\omega,\eta)})\  \mathrm d\mu^-(\omega)\mathrm d\mu^-(\eta).
\]

\subsection{Technical estimates}
Below we will make  use of the following estimates.

\begin{prop}\label{prop:transvdegree beta appendix}
    There exists $C>1$   such that for any $\epsilon\geq 0$ small enough the following holds. For any $\lambda\in U_0$ and any pair $\omega,\eta\in\{1,2\}^\mathbb N$ of the form $\omega=\tau\omega'$ and $\eta=\tau\eta'$ with $\omega_0'\neq \eta_0'$ we have
    \[
    |\partial_\lambda \Delta_{\lambda}(\omega,\eta)|\leq C (|\tau|+1) \ (\lambda+C\epsilon)^{|\tau|}\qquad\qquad   |\partial^2_{\lambda^2} \Delta_{\lambda}(\omega,\eta)|\leq C (|\tau|+1)^2 \ (\lambda+C\epsilon)^{|\tau|}.
    \]
Moreover,
    \[
     |\Delta_{\lambda}(\omega,\eta)|\leq  \frac{1}{C(|\tau|+1)}\ (\lambda-C \epsilon)^{|\tau|}\qquad \Longrightarrow\qquad  |\partial_\lambda\Delta_{\lambda}(\omega,\eta)|\geq \frac{1}{C} \ (\lambda-C \epsilon)^{|\tau|}.
    \]
\end{prop}

The estimates for the partial derivatives follow from computations analogue to those in Lemmas \ref{lem:bounded distortion} and \ref{lem:estimate for kappa dependence on lambda}. The second part is nothing but a reformulation of the property \ref{eq:transversality degree beta} and can be deduced from the proof of Proposition \ref{prop:transv degree beta}.

\begin{lem}
    Consider the setting of Proposition \ref{prop:transvdegree beta appendix} and suppose that the interval $U_0$ is small enough. Then for all $\omega,\eta\in \{1,2\}^\mathbb N$  the interval
    \begin{equation}\label{eq:monoton interval}
    U_0(\omega,\eta)=\{\lambda\in U_0\colon |\Delta_\lambda(\omega,\eta)|\leq \frac{1}{C(|\tau|+1)}\ (\lambda-C \epsilon)^{|\tau|}\}
    \end{equation}
    is either empty or it has one connected component on which the function $\lambda\mapsto \Delta_\lambda(\omega,\eta)$ is strictly monotone.
\end{lem}

This result follows plainly from the computations in the proof of Proposition \ref{prop:transv degree beta}. It will be of use when we later introduce a change of variables $t=\Delta_\lambda(\omega,\eta)$.

\subsection{The contribution of an individual pair}
We now estimate $\Delta S_r(u,\omega,\eta)$ for each pair $\omega,\eta\in\{1,2\}^\mathbb N$.

\begin{lem}\label{lem: estimate delta sr}
    Consider the setting of Proposition \ref{prop:transvdegree beta appendix} and write 
    \[
    \lambda_-=\mathrm{inf}_{\lambda\in U_1} \{\lambda-C\epsilon\}\qquad\qquad \lambda_+\mathrm{sup}_{\lambda\in U_1} \{\lambda+C\epsilon\}.
    \]
    There exists  a constant $C>1$ such that
    \[
    S_r(0;\omega,\eta)\leq C(|\tau|+1)\  \lambda_-^{-|\tau|}.
    \]
    Moreover, for any $u\in \mathbb R$, any $r<|u|$ and any  pair $\omega,\eta$ for which 
    \[
   \frac{\lambda_-^{|\tau|}}{C(|\tau|+1)}  \geq |u|
    \]
    we have
    \[
    |\Delta S_r(u;\omega,\eta)|\leq C(|\tau|^2+1) \left(\frac{\lambda_+}{\lambda_-^3}\right)^{|\tau|} |u|.
    \]
\end{lem}
\begin{proof}
    Let $U_0(\omega,\eta)$ be as in \eqref{eq:monoton interval} and introduce the change of variables $t=\Delta_\lambda(\omega,\eta)$. Writing $\lambda=\lambda(t)$ denote by 
   \[
g(t)=
\frac{\chi(\lambda(t))}
{|\partial_\lambda\Delta_{\lambda(t)}(\omega,\eta)|}.
\]
Hence, 
\[
|g'|
\leq
\frac{|\chi'|}{|\partial_\lambda\Delta(\omega,\eta)|^2}
+
\frac{|\chi| |\partial_\lambda^2\Delta(\omega,\eta)|}
{|\partial_\lambda\Delta(\omega,\eta)|^3}
\]
and, by the estimates in Proposition \ref{prop:transvdegree beta appendix}, we have 
\begin{equation}\label{eq:estimate g prime}
|g|\leq C\  \lambda_-^{|\tau|} \qquad\qquad |g'|\leq C(|\tau|^2+1)\left(\frac{\lambda_+}{\lambda_-^3}\right)^{|\tau|}.
\end{equation}
Since $\mathrm{supp}(K_r)\subset[-r,r]$, the assumption
\[
0<r\leq|u|,
\qquad
3|u|\leq\frac{\lambda_-^{|\tau|}}{C(|\tau|+1)}
\]
guarantees that whenever $K_r(t)\neq0$ both $t,t+u$ belong to the region where \eqref{eq:estimate g prime} holds. Extending $g$ by zero outside the image of the monotonicity interval  and making use of the change of variables formula
\[
S_r(0;\omega,\eta)=\int_{\mathbb R}K_r(t)g(t)\ \mathrm dt \leq |g|_{C^0} \int_{\mathbb R}K_r(t)\ \mathrm dt\leq C(|\tau|+1)\  \lambda_-^{-|\tau|}.
\]
Analogously,
\[
\Delta S_r(u;\omega,\eta)
=\int_{\mathbb R}K_r(t)(g(t+u)-g(t))\ \mathrm dt \leq |g'|_{C^0} |u|\int_{\mathbb R}K_r(t)\ \mathrm dt\leq C(|\tau|^2+1)\left(\frac{\lambda_+}{\lambda_-^3}\right)^{|\tau|}|u|.\qedhere
\]
\end{proof}

\subsection{Proof of Proposition \ref{prop:estimate Sobolev technical appendix}}

We now combine the estimates in Lemma \ref{lem: estimate delta sr} to obtain a proof of Proposition \ref{prop:estimate Sobolev technical appendix}. We proceed as follows.
\medskip

Given $u\in \mathbb R$ we let 
\[
n(u)=\max\{n\in\mathbb N\colon 3|u| \leq \frac{\lambda_-^{n}}{C(n+1)}\}.
\]
and make use of the estimates in Lemma \ref{lem: estimate delta sr} to write
\begin{align*}
\mathcal S_r(u)\leq & C|u| \sum_{n\leq n(u)}(|\tau|^2+1) \left(\frac{\lambda_+}{2\lambda_-^3}\right)^{|\tau|} + C\sum_{n>n(u)}(n+1)\left(\frac{1}{2\lambda_-}\right)^n
\end{align*}
and it follows after decreasing (if necessary) the size of the interval $U_0$ that for any 
\[
0<\sigma<\frac{\log 2}{\log(1/\lambda)}-1
\]
there exists a constant $C_\sigma$ such that for all  $\varepsilon $ small enough and all $0<r<|u|$
\[
\mathcal S_r(u)\leq C_\sigma |u|^\sigma.
\]
To complete the proof we make use of the formula \eqref{eq:Symbolic coding difference Alambda} to deduce the estimate \eqref{eq: estimate Sobolev technical appendix}.

\subsection{Proof of Theorem \ref{thm:sobolev}}
Fix any $\lambda_0\in (1/2,2^{-2/3})$ and let $U_0$ be a small neighbourhod around $\lambda_0$. Then, take any $\sigma$ as in Proposition \ref{prop:estimate Sobolev technical appendix} and let $0\leq s<\sigma/2$. Then, Proposition \ref{prop:estimate Sobolev technical appendix} implies that there exists $C$ such that for all $\varepsilon$ small enough
\[\int_{U_0}\left(\int_{\mathbb R}\int_{\mathbb R}
\frac{
|\tilde a_\lambda(x+u)-\tilde a_\lambda(x)|
}{|u|^{1+2s}}
\ \mathrm d x\mathrm du\right)\mathrm d\lambda\leq
C \int_{-1}^{1}|u|^{\sigma-1-2s}\ \mathrm du
<\infty,
\]
By Fubini we conclude that for almost every $\lambda\in U_0$ we have 
\[
\int_{\mathbb R}\int_{\mathbb R}
\frac{
|\tilde a_\lambda(x+u)-\tilde a_\lambda(x)|
}{|u|^{1+2s}}
\ \mathrm d x\mathrm du<\infty
\]
as we wanted to show.
\medskip

\section{The average growth of cylinder overlaps}\label{sec: average growth proof appendix}

In this appendix we study the average growth of cylinder overlaps. Our first observation is that for $\lambda$ close to one the growth of cylinder overlaps is not fast enough to verify \eqref{eq:heuristic rho with slambda}. To do so observe that for the affine case
\[
N_{\lambda,n}^{\mathrm{aff}}(z)=1\qquad\qquad \text{for }0<x<(1-\lambda)\lambda^n.
\]
Hence, 
\begin{equation}\label{eq:lowerbound affine growth average counting}
\int_{D} N_{\lambda,n}^{\mathrm{aff}}(z)^{-1} \mathrm dz\geq (1-\lambda)\lambda^n.
\end{equation}
Combined, \eqref{eq:heuristic rho with slambda} and \eqref{eq:lowerbound affine growth average counting} cannot hold simultaneously for $\lambda \geq 2^{-1/3}$.

\medskip
Having understood the existence of limitations to the average growth of cylinder overlaps we now give a proof of Proposition \ref{prop:averagegrowth}. It is divided in two steps. First, in Section \ref{sec:affine lambda star}, we control the growth rate for the affine system and extract a suitable interval for which the associated rate verifies \eqref{eq:heuristic rho with slambda}. Then, we handle the non-linear case in Section \ref{sec: nonlinear growth rate}.

\subsection{The affine case}\label{sec:affine lambda star}
In the affine case $N_{n,\lambda}^{\mathrm{aff}}$ only depends on $x$. Moreover, it is easy to check that (we write $f_{i,\lambda}$ for the horizontal component of the affine map $F_{i,\lambda}^{\mathrm{aff}}$ in \eqref{eq:affinemaps})
\[
N_{\lambda,n+1}(x)=\sum_{i\colon x\in f_{i,\lambda}([0,1])} N_{\lambda,n}(f_{i,\lambda}^{-1}(x)).
\]
Write $h_{\lambda,n}(x)=N_{\lambda,n+1}(x)^{-1}$. Outside the overlapping region $O_\lambda=(1-\lambda,\lambda)$ there is only one term in the sum so 
\[
h_{\lambda,n+1}(x)=h_n(f_{i,\lambda}^{-1}(x))\qquad\qquad\text{for }x\notin O_\lambda.
\]
Inside the overlapping region the inequality $(a+b)^{-1}\leq 4^{-1}(a^{-1}+b^{-1})$ shows that 
\[
h_{\lambda,n+1}(x)\leq \frac 14 \sum_{i\colon x\in f_{i,\lambda}([0,1])} h_{\lambda,n}(f_{i,\lambda}^{-1}(x))\qquad\qquad\text{for }x\in O_\lambda.
\]
Setting $q_\lambda(x)=1$ for $x\notin O_\lambda$ and $q_\lambda(x)=1/4$ for $x\in O_\lambda$ we obtain 
\begin{equation}\label{eq:recursion for h}
h_{\lambda,n+1}(x)\leq q_\lambda(x) \sum_{i\colon x\in f_{i,\lambda}([0,1])} h_{\lambda,n}(f_{i,\lambda}^{-1}(x))\qquad\qquad\text{for }x\in [0,1].
\end{equation}
We now use the recursion \eqref{eq:recursion for h} to obtain an estimate for $\int_0^1 h_{\lambda,n}(x)\ \mathrm dx$. 

\begin{lem}\label{lem:weigh integration estimate growth affine}
    For any $\lambda\in ((\sqrt 5-1)/2,1)$ we have 
    \[
    \int_0^1 h_{\lambda,n}(x)\ \mathrm dx\leq \left(\lambda\frac{1+\sqrt 3}{2}\right)^n.
    \]
\end{lem}
\begin{rem}
    The restriction $\lambda>(\sqrt 5-1)/2$ is only a limitation of our argument.  
\end{rem}

\begin{proof}
We consider a non-negative weight $w:[0,1]\to (0,\infty)$ and estimate instead
\[
I_{\lambda,n}=\int_0^1 w(x) h_{\lambda,n}(x)\mathrm dx.
\]
Using the recursion \eqref{eq:recursion for h}
\begin{equation}\label{eq:recursion for h integrated}
I_{\lambda,n+1}=\int_0^1 w(x) h_{\lambda,n+1}(x)\mathrm dx\leq \lambda \int_0^1 h_{\lambda,n}(x) \ ( \underbrace{(w q_\lambda)(f_{1,\lambda}(x)) + (w q_\lambda)(f_{2,\lambda}(x))}_{R_{\lambda}(x)}  )\ \mathrm dx.
\end{equation}
The first observation is that for $\lambda>(\sqrt 5-1)/2$ we have 
\begin{equation}\label{eq:observation for inclusions in overlap}
f_1((\lambda,1)), f_2((0,1-\lambda))\subset O_\lambda.
\end{equation}
Choose now $w(x)=1$ for $x\notin O_\lambda$, $w(x)=c$ for $x\in O_\lambda$ (for some $c\in(0,\infty)$ to be chosen) and observe that (in view of the inclusions \eqref{eq:observation for inclusions in overlap})
\[
R_\lambda(x)\leq (1+c/4)\  w(x)\qquad\qquad \text{for }x\notin O_\lambda
\]
and, trivially
\[
R_{\lambda(x)}\leq 2 = (2/c)\ w(x)\qquad\qquad \text{for }x\in O_\lambda.
\]
Hence, for any $c\in (0,1)$ \eqref{eq:recursion for h integrated} yields
\[
I_{\lambda,n+1}\leq \lambda \min\{1+c/4,2/c\}\  I_{\lambda,n}.
\]
Optimizing $c$ at $(1+\sqrt3)/2$ and using that $I_{\lambda,0}=1$ we deduce that for all $n\in \mathbb N$
\[
I_{\lambda,n}\leq \left(\lambda\frac{1+\sqrt 3}{2}\right)^n
\]
and the desired conclusion follows.
\end{proof}
\subsection{The nonlinear case: proof of Proposition \ref{prop:averagegrowth}}\label{sec: nonlinear growth rate}

The proof for the general case follows essentially from the fact that the nonlinear dynamics stays $\epsilon$-close to the affine dynamics.

\begin{lem}
    For a word $(\tau=(\tau_0,\ldots,\tau_{n-1})\in\mathcal W$, write
\[
x_0^\tau(t)=t,\qquad
x_{j+1}^\tau(t)=f_{\tau_j,\lambda}(x_j^\tau(t)).
\]
Fix $\kappa>0$ and set
\[
q_{\lambda,\kappa}(x)=
\begin{cases}
1/4 &1-\lambda+\kappa<x<\lambda-\kappa,\\
1 &\text{otherwise}.
\end{cases}
\]
Then, there exists $C$ such that for all $y\in [0,1]$ and for $\epsilon$ sufficiently small  (depending only on $\kappa$),
\begin{equation}\label{eq: nonlinear growth technical one}
\int_{\gamma_y} N_{\lambda,n}(z)^{-1}\ \mathrm dz
\leq (\lambda+C\epsilon)^n
\int_0^1\sum_{|\tau|=n}
\prod_{j=1}^n q_{\lambda,\kappa}(x_j^\tau(t))\ dt.
\end{equation}
\end{lem}

\begin{proof}

Fix any $y\in [0,1]$. To handle the non-linear dynamics we start by considering the curve 
\[
\gamma_{\tau,y,\lambda}=(\mathtt F^{|\tau|}_\tau)^{-1}(\gamma_y\cap V_{\tau,\lambda})
\]
and let $\{z(t)\}_{t\in[0,1]}$ stand for its parametrization over the horizontal coordinate. Then, we let 
\[
z_{j}^\tau(t)=\mathtt F_{\tau_j-1,\lambda}(z_{j-1}^\tau(t)),\qquad\qquad z_0^\tau(t)=z(t).
\]
We have shown in Section \ref{sec:geometric preli} that there exists $C$ such that uniformly in $\tau$ and for any $\epsilon$ small enough (here $\{x_j^\tau(t)\}_j$ is the horizontal orbit of the associated affine dynamics)
\begin{equation}\label{eq:close x dinamycs}
|\pi_x z_j^\tau(t)-x_j^\tau(t)|\leq C\epsilon.
\end{equation}
That is, the non-linear dynamics stays uniformly close to that of the affine system.
\medskip

    Write $h_{\lambda,n}=N_{\lambda,n}^{-1}$, and let $\tilde q_\lambda$ equal $1/4$ on $V_{1,\lambda}\cap V_{2,\lambda}$ and $1$ elsewhere. The same reciprocal inequality as in the affine case gives
\[
h_{\lambda,n+1}(z)\leq \tilde q_\lambda(z)
\sum_{i:\,z\in V_{i,\lambda}}
h_{\lambda,n}(\mathtt F_{i,\lambda}^{-1}z).
\]
Iterating, the contribution of a word $\tau$  is bounded by $\prod_{j=1}^n \tilde q_\lambda(z_j^\tau)$. Then, it follows from the change of variables formula and Lemma \ref{lem:bounded distortion} that 
\[
\int_{\gamma_y} N_{\lambda,n}(z)^{-1}\ \mathrm dz\leq (\lambda+C\epsilon)^n \int_0^1 \sum_{|\tau|=n}\prod_{j=1}^n \tilde q_\lambda
(z_j^\tau(t))\ \mathrm dt.
\]
But then, if $x_j^\tau(t)\in(1-\lambda+\kappa,\lambda-\kappa)$, then $z_j^\tau(t)\in V_{1,\lambda}\cap V_{2,\lambda}$. In particular,
\[
\tilde q_\lambda(z_j^\tau(t))
\leq q_{\lambda,\kappa}(x_j^\tau(t))
\]
and we obtain the estimate \eqref{eq: nonlinear growth technical one}.
\end{proof}

Let 
\[
Z_{\lambda,\kappa,n}(t)=\lambda^n\sum_{|\tau|=n}
\prod_{j=1}^n q_{\lambda,\kappa}(x_j^\tau(t)).
\]
and define
\[
\rho_\lambda= \lambda \frac{1+\sqrt 3}{2}.
\]
\begin{lem}\label{lem:bufferaffine}
    For $\lambda\in ((\sqrt 5-1)/2,1)$ fix  any $r_\lambda>\rho_\lambda$, any $\delta>0$ sufficiently small  and any $\kappa>0$ small enough. There exists $C$ such that for all $\epsilon$ small enough, all $\lambda\in ( (\sqrt 5-1)/2 +\delta, 1-\delta)$ and all $n\in\mathbb N$ we have 
    \[
    Z_{n,\lambda,\kappa}(t)\leq C r^n_\lambda\qquad\qquad\text{for all }t\in [0,1].
    \]
\end{lem}

\begin{proof}
    Arguing as in the proof of Lemma \ref{lem:weigh integration estimate growth affine} it is not difficult to check that  there exists  $m\in \mathbb N$ large enough such that for all $\epsilon$ small enough, all $\lambda\in ( (\sqrt 5-1)/2 +\delta, 1-\delta)$ we have 
    \[
    Z_{m,\lambda,\kappa}(t)\leq r^m_\lambda\qquad\qquad\text{for all }t\in [0,1].
    \]
    Writing $n=qm+l$
\[
\sup_t Z_{n,\lambda,\kappa}(t)
\leq (\sup_t Z_{m,\lambda,\kappa}(t))^q\,
      \sup_t Z_{l,\lambda,\kappa}(t)\leq C r_\lambda^{qm}.
\]
gives the desired conclusion.
\end{proof}

We are ready to complete the proof of Proposition \ref{prop:averagegrowth}.

\begin{proof}[Proof of Proposition \ref{prop:averagegrowth}]

We observe that for any $\lambda< 2^{-2/3}$ we have  
\[
\rho_\lambda \sqrt{2\lambda}<1.
\]
Thus, if we fix  any $\delta>0$ and any $r_\lambda<(2\lambda)^{-1/2}$,  Lemma \ref{lem:bufferaffine}  shows that there exist $C>0$ and $\epsilon_0>0$ such that for any $0\leq \epsilon\leq \epsilon_0$  and  any $\lambda\in ( (\sqrt 5-1)/2 +\delta, \ 2^{-2/3}-\delta)$ we have 
\[
\int_D N_{\lambda,n}(z)^{-1}\ \mathrm dz\leq C r_\lambda^n.
\]
The proof of Proposition \ref{prop:averagegrowth} is complete.

\end{proof}
\medskip

\section{Shadowing orbits of the IFS}\label{sec: shadowing}

In this appendix we give the proof of Lemma \ref{lem:transportskewprod}. The proof is divided in a number of steps:
\begin{enumerate}
    \item First, in Section \ref{sec: results from old} we gather a few results obtained in \cite{paradela2026robustlytransitivebehaviorsymplectic} concerning the IFS generated by the maps $\{T_\varepsilon,\{S_{i,\varepsilon}\}_i\}$ acting on $N$.
    \item Second, in Section \ref{sec:compactification} we introduce a suitable compactification of the lamination $\mathcal X_\varepsilon$ in Proposition \ref{prop:mainnhim}. The reason is that to shadow the orbits of the IFSs generated by the maps $\{T_\varepsilon,\{S_{i,\varepsilon}\}_i\}$ we will need to approach arbitrarily close to the invariant sets $W^{u,s}(\Lambda_\varepsilon)$. The purpose of this compactification  is to give an extension of the center dynamics for orbits in $W^{u,s}(\Lambda_\varepsilon)$. 

    \item Finally, in Section \ref{sec: proof shadowing} we complete the proof of Lemma \ref{lem:transportskewprod}.
\end{enumerate}

\subsection{Transport in the IFS}\label{sec: results from old}

The following result was obtained in \cite{paradela2026robustlytransitivebehaviorsymplectic} (see Proposition 5.3 therein).

\begin{prop}\label{prop:epsaccesibility}

  Let $N$ be a smooth compact Riemannian surface and denote by $\rho$ the associated area form.   Let $T_\varepsilon:N\to N$ be a diffeomorphism which preserves a smooth measure and let  $\{S_{i,\varepsilon}\}_{i=1}^m:N\to N$ be diffeomorphisms preserving a (possibly different) smooth measure  and satisfying \ref{it:B6} for some $n\geq 0$.  Then, there exists $K>0$ and $\varepsilon_0>0$ such that for any $0<\varepsilon\leq \varepsilon_0$, there exists a  subset $R_\varepsilon\subset N$ with $\rho(N\setminus R_\varepsilon)=0$ such that the following holds.
  
  For   any two points $z,z_*\in N$ with $z\in R_\varepsilon$ there exists $L\in \mathbb N$,   $\iota\in\{1,2\}^L$ such that for any sequence $\{n_r\}_{r=1}^{L-1}$ satisfying $ n_1\gg \cdots \gg n_{L-1}\gg 1$ there exists  $\omega\in \mathbb N^L$ satisfying 
    \begin{equation}\label{eq:increasingtimes}
\omega_{r-1}-\omega_{r}\geq n_r \qquad\qquad \text{for all }r=1,\dots, L-1\quad\text{and } \omega_{L-1}\geq n_*
    \end{equation}
    and such that 
    \[
    \mathrm{dist}(T^{\omega_{L-1}}\circ S_{\iota_{L-1}}\circ\cdots\circ T^{ \omega_0}\circ S_{\iota_0}(z),\ z_*)\leq K\varepsilon^{n+1}.
    \]
\end{prop}

We refer the interested reader to \cite{paradela2026robustlytransitivebehaviorsymplectic} for details on the proof. Here we just outline the main idea. 
\begin{rem}
    To be precise, in \cite{paradela2026robustlytransitivebehaviorsymplectic}  Proposition \ref{prop:epsaccesibility} was only established for the case $n=0$. The reader can easily check that the very same argument gives the general case. 
\end{rem}

The first step is to show that, Assumption \ref{it:B6} guarantees that for any two points $z,z_*\in N$ there exists an orbit of the IFS generated by $\{\{S_{i,\varepsilon}\}_i\}$ which connects $z$ to a $O(\varepsilon^{n+1})$-neighbourhood of $z_*$. The proof of this result follows from a trivial adaptaition of that of the Ball-box Theorem (see Proposition 5.1. in \cite{paradela2026robustlytransitivebehaviorsymplectic}). Since the maps $T,\{S_{i,\varepsilon}\}_i$ preserve a smooth measure, making use of the Poincaré recurrence theorem (hence the restriction to a full measure subset $R_\varepsilon$), we deduce that the same conclusion holds for orbits of the IFS generated by $\{T^n\circ S_{i,\varepsilon}\}_{i,n}$. Moreover, it is not difficult to check that the orbits realizing the construction can be chosen to have recurrence times with arbitrarily large spacings (see \eqref{eq:increasingtimes}). This last property will prove crucial in Section \ref{sec: proof shadowing} to implement a shadowing argument and substitute the IFS dynamics by the true skew-product dynamics.

\medskip

\subsection{Compactification}\label{sec:compactification}
In  Proposition \ref{prop:mainnhim} we showed that the dynamics on the lamination $\mathcal X_\varepsilon$  is driven by the skew-product 
  \begin{align*}
    \mathcal F_\varepsilon: (\mathbb N\times\{1,\dots, m\})^\mathbb Z\times N&\to (\mathbb N\times\{1,\dots,m\})^\mathbb Z\times N\\
    (\bs \omega,z)&\mapsto (\sigma(\bs\omega), F_{\bs\omega,\varepsilon}(z))
    \end{align*}
   where 
   \[
   F_{\bs\omega,\varepsilon}(z)=T_\varepsilon^{\omega_0}\circ (S_{\iota_0,\varepsilon}+R_{\bs\omega,\varepsilon})(z).
   \]
   with $T_\varepsilon:N\to N$ being the inner map \eqref{eq:innerdynamicsdefn} and $\{S_{i,\varepsilon}\}_{i=1}^m:N\to N$ being the system of scattering maps in \eqref{eq:scattmapdynamicsdefn}. Following Moser \cite{MR442980} we introduce the compactification of  \[
 \Omega=(\mathbb N\times\{1,\dots,m\})^\mathbb Z
 \]
 obtained by including also sequences of the form $\bs\omega=(\omega,\iota)$ with 
    \begin{enumerate}
        
         \item (right-finite) $\omega=(\dots\omega_{-l},\dots;\omega_0,\dots,\omega_{r},\infty)$ and $\iota=(\dots,-\iota_{-l},\dots;\iota_0,\dots,\iota_r)$ for some $r\in\mathbb N$,
          \item (left-finite) $\omega=(\infty,\omega_{-l},\dots;\omega_0,\dots,\omega_{r},\dots)$ and $\iota=(-\iota_{-l},\dots;\iota_0,\dots,\iota_r,\dots)$ for some $l\in\mathbb N$,
         \item (finite) $\omega=(\infty,\omega_{-l},\dots;\omega_0,\dots,\omega_{r},\infty)$ and $\iota=(-\iota_{-l},\dots;\iota_0,\dots,\iota_r)$ for some $l,r\in\mathbb N$.
    \end{enumerate}
  We denote the corresponding space by $\Omega_c$ and introduce as neighbourhoods of an element $\bs\omega^*$ as in (1), the sets of $\bs \omega$ such that 
  \[
  \omega^*_k=\omega_k\qquad\text{for }-n\leq k\leq r,\qquad\qquad \omega_{r+1}\geq n
  \]
  and
  \[
  \iota^*_k=\iota_k\qquad\text{for }-n\leq k\leq r
  \]
  for $n\in\mathbb N$. Similarly we introduce the set of neighbourhoods of an element $\bs\omega^*$ as in (2) or (3). In this way we obtain a topology on $\Omega_c$ extending that of $\Omega$ and for which $\Omega_c$ is compact. 
  \medskip

  The key point in introducing this extension is that, if we fix any $\bs\omega^*$ as in (1), for any sequence $\{\bs\omega^{(k)}\}_k\in \Omega_c$ such that $\bs\omega^{(k)}\to \bs\omega^*$ the sequence of functions
  \[
R_{\bs\omega^{(k)},\varepsilon}:N\to N
  \]
  converges uniformly (in the $C^{2}$ topology) to a function
  \[
  R_{\bs\omega^*,\varepsilon}:N\to N.
  \]
Similar (suitably modified) statements hold sequences $\bs\omega^*$ of types (2) or (3). The skew-product dynamics is only defined on  
\[
\Omega_{\mathrm{adm}}=\{(\bs\omega,z)\in \Omega_c\times N\colon \omega_0<\infty\}.
\]
However, for $\bs\omega^*\in \Omega$ such that $\omega_0^*=\infty$ we can define a map
\begin{align*}
F_{\bs\omega^*,\varepsilon}(z):N&\to N\\z&\mapsto (S_{\iota_0,\varepsilon}+R_{\bs\omega^*,\varepsilon})(z).
\end{align*}
It can also be shown (using the estimates obtained in the proof of  Lemma \ref{lem:formulas holonomies}) that, given any $\bs\omega\in \Omega$, any $\bs\omega^*\in \Omega_c$ of type (1) and any sequence  $\{\bs\omega^{(k)}\}_k\in \Omega_c$ such that $\bs\omega^{(k)}\to \bs\omega^*$ the sequence of maps 
\[
H^u_{\bs\omega,\bs\omega^{(k)}}:N_{\bs\omega}\mapsto N_{\bs\omega^{(k)}}
\]
converges uniformly in the $C^{2}$ topology. Analogous statements could be formulated for $\bs\omega^*$ of types (2)-(3) but we omit them since we will not need them.
\medskip

\subsection{Proof of Lemma \ref{lem:transportskewprod}}\label{sec: proof shadowing}
We omit the dependence of $T_\varepsilon, S_{i,\varepsilon}$ and $R_{\bs\omega,\varepsilon}$ in $\varepsilon$ in order to alleviate the notation. Consider the right-finite sequence 
\[
\bs\omega^*=(\omega^*,\iota^*)\qquad\qquad\text{with}\qquad\qquad \omega^*=(\dots,\omega_{-1};\infty)\qquad \iota^*=(\dots,\iota_{-1};\iota_0)
\]
 and define the point
    \[
\tilde z(\bs\omega^*)=(S_{\iota_0}+R_{\bs\omega^*})\circ H^u_{\bs\omega,\bs\omega^*}(z).
    \]
    In view of \ref{eq:Romega estimate C1} and \ref{eq:comparing holonomies} we have that for any $\bs\omega'=(\omega',\iota')$ with $\omega'=(\dots,\omega_{-1};\omega'_0,\dots)$ and $\iota'=(\dots,\iota_{-1};\iota_0,\iota_1',\dots)$
\[
\lVert R_{\bs\omega'}-R_{\bs\omega^*}\rVert_{C^2}\leq C \bar\lambda^{\omega_0'},\qquad\qquad \lVert H^u_{\bs\omega,\bs\omega'}-H^u_{\bs\omega,\bs\omega^*}\rVert_{C^2}\leq C\bar\lambda^{\omega_{-1}+\omega_0'}
\]
    so
    \[
    \mathrm{dist}((S_{\iota_0}+R_{\bs\omega'})\circ H^u_{\bs\omega,\bs\omega'}(z),\tilde z(\bs\omega^*))\leq C \bar\lambda^{\omega_0'}.
    \]
    Moreover, since $z\in R_\varepsilon$ and  the transformation $z\in N_{\bs\omega}\mapsto (S_{\iota_0}+R_{\bs\omega^*})\circ H^u_{\bs\omega,\bs\omega^*} (z)\in N_{\omega^*}$ is nonsingular, we can assume that $\tilde z(\bs\omega^*)\in R_\varepsilon$ (after deleting, if necessary, a zero measure subset).
    \medskip
    
   Let $n$ be as in the statement of Lemma \ref{lem:transportskewprod} and denote by 
   \[
   \delta=\varepsilon^{n+1}.
   \]
   We now proceed as in the proof of Lemma 6.1. in \cite{paradela2026robustlytransitivebehaviorsymplectic}, taking $\tilde z(\bs\omega^*)$ as starting point. Applying Proposition 5.3. in \cite{paradela2026robustlytransitivebehaviorsymplectic} we find $L\in\mathbb N$ and $\bs{\hat\omega}\in(\mathbb N\times A)^L$ such that
    \[
    \mathrm{dist}(T^{\hat\omega_{L-1}}\circ S_{\hat\iota_{L-1}}\circ\cdots T^{ \hat\omega_0}\circ S_{\hat\iota_0}(\tilde z(\bs\omega^*),\ z_*)\leq C\delta
    \]
    and, for any $r=0,\dots,L-1$
    \[
    \bar\lambda^{\hat\omega_r}\leq \frac{\delta}{LK_r}\qquad\qquad\text{where}\qquad K_r=\lVert T^{\hat\omega_{L-1}}\circ S_{\hat\iota_{L-1}}\circ\cdots T^{\hat\omega_r}\circ S_{\hat\iota_r}\rVert_{C^1}.
    \]
 Fix now any $\tilde\omega_0\gg \hat\omega_0$ such that 
    \[
  \bar\lambda^{\tilde\omega_0}\leq\frac{\delta}{K_0}\qquad\qquad\text{and}\qquad\qquad \mathrm{dist}( T^{\tilde\omega_0}\tilde z(\bs\omega^*),\ \tilde z(\bs\omega^*))\leq \frac{\delta}{K_0}
    \]
We claim that $\tilde{\bs\omega}=(\tilde\omega,\tilde\iota)\in (\mathbb N\times\{1,\dots,m\})^{L+1}$ with 
\[
\tilde\omega=(\tilde\omega_0,\hat\omega_0,\dots,\hat\omega_{L-1}),\qquad\qquad \tilde\iota=(\iota_0,\hat\iota_0,\dots,\hat\iota_{L-1})
\]
verifies that 
\[
\mathrm{dist}(F_{\iota_{\bs{\tilde\omega}(\bs\omega)}}^{L+1}\circ H^u_{\bs\omega,\iota_{\bs{\tilde\omega}(\bs\omega)}}(z),\ z_*)\leq C \delta.
\]
This claim can be verified  by making use of the above estimates and the very same  telescoping argument used in the proof of Lemma 6.1. of \cite{paradela2026robustlytransitivebehaviorsymplectic}.
\medskip

Finally, since the stable holonomy $H^s_{\iota^-_{\bs{\tilde\omega}(\bs\omega)},\bs\omega}$ satisfies (see Proposition \ref{prop: close to id holonomies})
\[
\lVert H^s_{\iota^{-}_{\bs{\tilde\omega}(\bs\omega)},\bs\omega}-\mathrm{id}\rVert_{C^2}\leq C\bar\lambda^{\min\{\omega_{-1},\tilde\omega_{L-1}\}}\leq C\delta
\]
and 
\[
\mathcal H_{\bs{\tilde\omega}}=H^s_{\iota^{-}_{\bs{\tilde\omega}(\omega)},\omega}\circ F_{\iota_{\bs{\tilde\omega}(\omega)}}^{L+1}\circ H^u_{\omega,\iota_{\bs{\tilde\omega}(\omega)}}
\]
we obtain that  
\[
\mathrm{dist}(\mathcal H_{\bs{\tilde\omega}}(z), z_*)\leq K \delta = K\varepsilon^{n+1}.
\]

\bibliography{biblio}
\bibliographystyle{alpha}
\end{document}